\documentclass[UTF-8,reqno]{amsart}
\usepackage{enumerate, bbm}
\usepackage{amssymb,url,color, booktabs}
\usepackage{mathrsfs}

\usepackage{color}
\usepackage[colorlinks=true]{hyperref}
\hypersetup{
    linkcolor=blue,          
    citecolor=red,        
    filecolor=blue,      
    urlcolor=cyan
}

\usepackage{color}
\definecolor{MyDarkBlue}{cmyk}{0.8,0.3,0.8,0.4}
\definecolor{yellow}{rgb}{0.99,0.99,0.70}
\definecolor{white}{rgb}{1.0,1.0,1.0}
\definecolor{black}{rgb}{0.00,0.00,0.00}

\numberwithin{equation}{section}

\newcommand{\be}{\begin{eqnarray}}
\newcommand{\ee}{\end{eqnarray}}
\newcommand{\ce}{\begin{eqnarray*}}
\newcommand{\de}{\end{eqnarray*}}
\newtheorem{theorem}{Theorem}[section]
\newtheorem{lemma}[theorem]{Lemma}
\newtheorem{remark}[theorem]{Remark}
\newtheorem{definition}[theorem]{Definition}
\newtheorem{proposition}[theorem]{Proposition}

\newtheorem{corollary}[theorem]{Corollary}

\def\nor{|\mspace{-3mu}|\mspace{-3mu}|}

\def\eps{\varepsilon}

\def\e{\mathrm{e}}

\def\u{\mathbf{u}}

\def\p{\partial}

\def\[{{\Big[}}
\def\]{{\Big]}}
\def\<{{\langle}}
\def\>{{\rangle}}
\def\({{\big(}}
\def\){{\big)}}

\def\tr{\mathrm {tr}}
\def\W{{\mathcal W}}

\def\dif{{\mathord{{\rm d}}}}

\def\min{{\mathord{{\rm min}}}}

\def\bbp{{\boldsymbol{p}}}
\def\bbr{{\boldsymbol{r}}}
\def\bbq{{\boldsymbol{q}}}

\def\bb2{{\boldsymbol{2}}}
\def\no{\nonumber}
\def\={&\!\!=\!\!&}

\def\cB{{\mathcal B}}
\def\cC{{\mathcal C}}

\def\cH{{\mathcal H}}
\def\cI{{\mathcal I}}
\def\cJ{{\mathcal J}}

\def\cL{{\mathcal L}}
\def\cM{{\mathcal M}}
\def\cN{{\mathcal N}}

\def\cP{{\mathcal P}}

\def\mB{{\mathbb B}}
\def\mC{{\mathbb C}}

\def\mE{{\mathbb E}}

\def\mI{{\mathbb I}}

\def\mL{{\mathbb L}}

\def\mN{{\mathbb N}}

\def\mP{{\mathbb P}}
\def\mQ{{\mathbb Q}}
\def\mR{{\mathbb R}}

\def\bP{{\mathbf P}}

\def\bJ{{\mathbf J}}

\def\b1{{\mathbbm 1}}

\def\sE{{\mathscr E}}
\def\sF{{\mathscr F}}

\def\sI{{\mathscr I}}

\def\sL{{\mathscr L}}
\def\sM{{\mathscr M}}

\def\sU{{\mathscr U}}

\def\sX{{\mathscr X}}

\def\geq{\geqslant}
\def\leq{\leqslant}
\def\ge{\geqslant}
\def\le{\leqslant}

\def\div{\mathord{{\rm div}}}

\def\eps{\varepsilon}

\def\e{\mathrm{e}}

\def\u{\mathbf{u}}

\def\p{\partial}

\def\[{{\Big[}}
\def\]{{\Big]}}
\def\<{{\langle}}
\def\>{{\rangle}}

\def\bx{{\mathbf{x}}}
\def\tr{\mathrm {tr}}
\def\W{{\mathcal W}}

\def\dif{{\mathord{{\rm d}}}}

\def\min{{\mathord{{\rm min}}}}

\def\no{\nonumber}
\def\={&\!\!=\!\!&}
\def\bt{\begin{theorem}}
\def\et{\end{theorem}}
\def\bl{\begin{lemma}}
\def\el{\end{lemma}}
\def\br{\begin{remark}}
\def\er{\end{remark}}
\def\bd{\begin{definition}}
\def\ed{\end{definition}}
\def\bp{\begin{proposition}}
\def\ep{\end{proposition}}
\def\bc{\begin{corollary}}
\def\ec{\end{corollary}}

\def\geq{\geqslant}
\def\leq{\leqslant}
\def\ge{\geqslant}
\def\le{\leqslant}

\def\div{\mathord{{\rm div}}}

\def\bW{{\mathbf W}}
\def\bw{{\boldsymbol w}}
\def\bX{{\mathbf X}}

 \def\R{\mathbb R}
 \def\R{\mathbb R}

\def\<{\langle} \def\>{\rangle}
\def\W{{\widetilde W}}
\def\X{{\widetilde X}}
\def\x{{\mathbf x}}
\def\wt{\widetilde}

\allowdisplaybreaks

\begin{document}

\title[Propagation of chaos of McKean-Vlasov SDEs]
{Strong convergence of propagation of chaos for McKean-Vlasov SDEs with singular interactions}

\author{Zimo Hao, Michael R\"ockner, Xicheng Zhang}

\address{Zimo Hao:
Fakult\"at f\"ur Mathematik, Universit\"at Bielefeld,
33615, Bielefeld, Germany\\ 
and School of Mathematics and Statistics, Wuhan University,
Wuhan, Hubei 430072, P.R.China\\
Email: zhao@math.uni-bielefeld.de}

\address{Michael R\"ockner:
Fakult\"at f\"ur Mathematik, Universit\"at Bielefeld,
33615, Bielefeld, Germany\\
and Academy of Mathematics and Systems Science, CAS, Beijing, China\\
Email: roeckner@math.uni-bielefeld.de
 }

\address{Xicheng Zhang:
School of Mathematics and Statistics, Wuhan University, Wuhan, Hubei 430072, China\\
School of Mathematics and Statistics, Beijing Institute of Technology, Beijing 100081, China\\
Email: XichengZhang@gmail.com
 }

\thanks{
This work is partially supported by NNSFC grants of China (Nos. 12131019, 11731009), and the German Research Foundation (DFG) through the 
Collaborative Research Centre(CRC) 1283/2 2021 - 317210226 ``Taming uncertainty and profiting from randomness and low regularity in analysis, stochastics and their applications".}

\begin{abstract}
In this work we show the strong convergence of propagation of chaos for the particle approximation of McKean-Vlasov SDEs 
with singular $L^p$-interactions as well as for the moderate interaction particle systems on the level of particle trajectories. One of the main obstacles 
is to establish the strong well-posedness of the SDEs for particle systems with singular interaction. To this end, we extend the results on strong well-posedness of Krylov and R\"ockner \cite{Kr-Ro} to  
the case of mixed $L^\bbp$-drifts, where the heat kernel estimates play a crucial role. Moreover, when the interaction kernel is bounded measurable, 
we also obtain the optimal rate of strong convergence, 
which is partially based on Jabin and Wang's entropy method \cite{JW16} and Zvonkin's transformation.


\end{abstract}
\maketitle \rm

\tableofcontents

\section{Introduction}
Let $\phi:\mR_+\times\mR^d\times\mR^d\to\mR^m$, $F:\mR_+\times\mR^d\times\mR^m\to\mR^d$ 
and $\sigma:\mR_+\times\mR^d\to\mR^d\otimes\mR^d$ be Borel measurable functions.
For a (sub)-probability measure $\mu$ over $\mR^d$, we define 
$$
b(t,x,\mu):=F(t,x,(\phi_t\circledast\mu)(x)),
$$
where { $\phi_t(x,y):=\phi(t,x,y)$} and
$$
(\phi_t\circledast\mu)(x):=\int_{\mR^d}\phi_t(x,y)\mu(\dif y).
$$
Consider the following interacting system of  $N$-particles,
\begin{align}\label{MV2}
\dif X^{N,i}_t=b\(t,X^{N,i}_t,\eta_{\bX^N_t}\)\dif t+\sigma\(t,X^{N,i}_t\)\dif W^i_t,\ \ i=1,\cdots,N,
\end{align}
where $\eta_{\bX^N_t}$ stands for the empirical distribution measure of $N$-particles $\bX^N_t:=(X^{N,1}_t,\cdots,X^{N,N}_t)$,
$$
\eta_{\bX^N_t}(\dif y):=\frac1{N}\sum_{j=1}^N\delta_{X^{N,j}_t}(\dif y),
$$
and $\{W^i,i\in N\}$ is a sequence of independent standard Brownian motions on some stochastic basis $(\Omega,\sF,\mP, (\sF_t)_{t\geq 0})$. 
The infinitesimal generator of the above system is given by
$$
\cL^N_t\varphi(\bx)=\tr\big(a(t,x^i)\cdot\nabla^2_{x^i}\varphi(\bx)\big)+F\Big(t,x^i, \frac{1}{N}\sum_{j=1}^N\phi_t(x^i,x^j)\Big)\cdot\nabla_{x^i}\varphi(\bx),
$$
where $\x=(x^1,\cdots,x^N)\in(\mR^{d})^N$ and $a=\frac12\sigma\sigma^*$. Here and below we use Einstein's convention for summation.

\medskip

In this paper we are mainly concerned with the weak and strong convergence of {the solutions to \eqref{MV2}}
with general $L^p$-{\it singular} interaction $\phi_t(x,y)$ to {the solution of} the following distribution-dependent 
(or McKean-Vlasov) SDE  (abbreviated as DDSDE)
when $N\to\infty$:
\begin{align}\label{MV1}
\dif X_t=b(t,X_t,\mu_{X_t})\dif t+\sigma(t,X_t)\dif W^1_t,
\end{align}
where $\mu_{X_t}$ denotes the distribution of $X_t$. In particular, 
$\mu:=(\mu_{X_t})_{t\geq 0}$ solves the following nonlinear Fokker-Planck equation in distributional sense:
$$
\p_t\mu= \p_i\p_j(a_{ij}\mu)+\div (b(\mu)\mu),
$$
Moreover, we are also interested in the so called moderately interacting kernel $\phi_t(x,y)=\phi_{\eps_N}(x-y)$,
where $\phi_{\eps_N}$ is a family of mollifiers and $\eps_N\to 0$ as $N\to\infty$. In this case, the {solution to the} interacting particle system
\begin{align}\label{MV30}
\dif X^{N,i}_t=F\(t,X^{N,i}_t,(\phi_{\eps_N}\circledast\eta_{\bX^N_t})(X^{N,i}_t)\)\dif t+\sigma\(t,X^{N,i}_t\)\dif W^i_t,\ \ i=1,\cdots,N,
\end{align}
is expected to converge to {the solution of} the following density-dependent SDE (see \cite{Oe85,JM98}):
\begin{align}\label{MV3}
\dif X_t=F(t,X_t,\rho_{X_t}(X_t))\dif t+\sigma(t,X_t)\dif W_t,
\end{align}
where $\rho_{X_t}$ stands for the density of $X_t$. Here $\rho:=(\rho_{X_t})_{t\geq 0}$ solves the following nonlinear and {\it local (or Nemytskii-type)} Fokker-Planck equation:
\begin{align}\label{MV4}
\p_t\rho= \p_i\p_j(a_{ij}\rho)+\div (F(\rho)\rho).
\end{align}
It should be kept in mind that for $d=1$ and $F(\rho)=\rho$, this is Burgers-type equation.

\medskip

For the motion of a single particle, when $\phi\in L^q_t(L^p_x)$ with $\frac{d}p+\frac2q<1$, Krylov and R\"ockner \cite{Kr-Ro}  
showed the existence and uniqueness of strong solutions to the following SDE by Girsanov's transformation:
$$
\dif X_t=\phi_t\(X_t\)\dif t+\dif W_t.
$$
Later,  Zhang \cite{Z11} extended their result to the multiplicative noise case by Zvonkin's transformation of \cite{Zv}
(see also \cite{Zh-Zh1, Xi-Xi-Zh-Zh}).
However, for $N$-particle 
system
\eqref{MV2} with $F(t,x,r)=r$ and $\phi_t(x,y)=\phi_t(x-y)$, where $\phi$ is as above,
one can not use these well-known results {for $L^q_t(L^p_x)$ drifts} to derive the well-posedness by considering \eqref{MV2} as an  SDE in $\mR^{Nd}$.
For instance, when {$N=3$}, consider the following SDE in {$\mR^{3d}$}:
\begin{align}\label{AS09}
\left\{
\begin{aligned}
&\dif X^{1}_t=\[\phi\(X^{1}_t,X^{2}_t\)+\phi\(X^{1}_t,X^{3}_t\)\]\dif t+\dif W^1_t,\\
&\dif X^{2}_t=\[\phi\(X^{2}_t,X^{1}_t\)+\phi\(X^{2}_t,X^{3}_t\)\]\dif t+\dif W^2_t,\\
&\dif X^{3}_t=\[\phi\(X^{3}_t,X^{1}_t\)+\phi\(X^{3}_t,X^{2}_t\)\]\dif t+\dif W^3_t,
\end{aligned}
\right.
\end{align}
where $|\phi(x,y)|\leq h(x-y)$ and $h\in L^p$ with $p>d$. For $i=1,2,3$, let $\phi_i(x_1,x_2,x_3):=\sum_{j\not=i}\phi(x_i,x_j)$.
{As a function of {$(x_1,x_2,x_3)$} in {$\mR^{3d}$}, one only has 
\begin{align}\label{HZA00}
\phi_i\in L^\infty_{x^*_i}L^p_{x_i},\quad i=1,2,3,
\end{align}
where $x^*_i$ stands for the remaining variables except for $x_i$. It does not satisfy the conditions in \cite{Kr-Ro}.
Note that in the same work \cite{Kr-Ro}, Krylov and R\"ockner also showed the strong well-posedness for a class of special stochastic particle system
with singular gradient interaction $\phi=\nabla V$, where $V$ is continuously differentiable on $\mR^d\backslash\{0\}$ and satisfies some other conditions (see Section 9 in \cite{Kr-Ro}). Moreover, the strong well-posedness for particle system with Biot-Savart law interaction kernel 
$\phi(x)=(-x_2,x_1)/|x|^2$ was established in \cite{Os85} and \cite{FM07}, 
which is related to the random point vortex approximation for two dimensional Navier-Stokes equations. 
In the above well-known works, the key point of establishing the strong well-posedness is to prove that the process $X^i_t-X^j_t$ for $i\not=j$ 
does not touch the singular point $0$, i.e. the state space is $\mR^{Nd}$ 
``without diagonals''.
However, the strong well-posedness for particle systems \eqref{MV2} with general $L^p$-interaction kernels on all of $\mR^{Nd}$ has still been open.
}

\medskip

Therefore, our first task is to extend \cite{Kr-Ro, Z11} to the case of mixed $L^p$-spaces. We mention here that although Ling and Xie \cite{LX21} have already
considered singular SDEs in mixed $L^p$-spaces, their result cannot be applied to equation \eqref{AS09}
due to the new feature that we need to consider the order of the integral in {$x_1,x_2,x_3$} as well as the different integrability indices. Notice
that each $\phi_i$ belongs to a different mixed $L^p$-space.
For DDSDE \eqref{MV1}, in \cite{RZ21}, the last two authors of the present paper 
have already shown the weak and strong well-posedness 
 (see also \cite{Li-Mi} and \cite{Mi-Ve} for bounded measurable interaction kernel). 
Furthermore, weak solutions to the distribution density-dependent SDE \eqref{MV3}, were constructed in \cite{Ba-Ro1}, first solving the corresponding Fokker-Planck-Kolmogorov equation 
and using the superposition principle, and strong solutions were constructed in \cite{HRZ21} by directly using Euler's scheme. Recently, Wang \cite{Wa21} studied the weak and strong well-posedness for
more general distribution density-dependent SDEs with singular coefficients by a fixed point argument, but not for mixed $L^\bbp$-drifts. Nowadays, there is a vast literature about McKean-Vlasov or mean-field SDEs. 
We do not intend to list all the papers here. The interesting reader is referred to the references in the already mentioned papers.

\subsection{Propagation of chaos}In this subsection we recall some notions and well-known results about the propagation of chaos.

{\bf Kac's chaos:} Let $E$ be a Polish space and $\mu\in\cP(E)$ a probability measure on $E$. Let $(\mu^N)_{N\in\mN}$ be a sequence of symmetric probability measures on 
the respective product space $E^N$, where symmetric means that for any permutation $(x_{i_1},\cdots, x_{i_N})$ of $(x_1,\cdots,x_N)$,
$$
\mu^N(\dif x_{i_1},\cdots,\dif x_{i_N})=\mu^N(\dif x_1,\cdots,x_N).
$$
In particular, $\mu^N$ has a common $1$-marginal distribution. 
One says that $(\mu^N)_{N\in\mN}$ is {\it $\mu$-chaotic} if for any $k\in\mN$ (see \cite{Kac}),
\begin{align}\label{Kac}
\mu^{N,k}\mbox{ weakly converges to $\mu^{\otimes k}$ as $k\leq N\to\infty$},
\end{align}
where $\mu^{N,k}(\dif x_1,\cdots,\dif x_k)=\mu^N(\dif x_1,\cdots,\dif x_k, E,\cdots, E)$ is the $k$-fold marginal distribution of $\mu^N$.
It is well known that \eqref{Kac} holds if and only if \eqref{Kac} holds for only $k=2$ (see \cite[(i) of Proposition 2.2]{Sz}). In the language of random variables, Kac's chaos can be restated as follows:
Let ${\boldsymbol\xi^N}:=(\xi^{N,1},\cdots,\xi^{N,N})$ be a family of $E$-valued random variables. If the law of ${\boldsymbol\xi^N}$ is symmetric and $\mu$-chaotic, one says that
${\boldsymbol\xi^N}$ is $\mu$-chaotic. It is also equivalent to  (see \cite[(ii) of Proposition 2.2]{Sz})
\begin{align}\label{LAR}
\mbox{the empirical measure $\eta_{{\boldsymbol\xi}^N}(\dif y):=\frac1N\sum_{j=1}^N\delta_{\xi^{N,j}}(\dif y)\in\cP(E)$ converges to $\mu$ in law.}
\end{align}
Note that ${\boldsymbol\xi^N}$ can be regarded as $N$-random particles in state space $E$. From this viewpoint, Kac's chaos means that if one observes the distribution of any $k$-particles,
then they become statistically independent as $N$ goes to infinity. Indeed, \eqref{LAR} is a law of large numbers, i.e., for any $\varphi\in C_b(E)$,
$$
\eta_{{\boldsymbol\xi}^N}(\varphi):=\frac1N\sum_{j=1}^N\varphi(\xi^{N,j}) \to \mu(\varphi):=\int_E\varphi(x)\mu(\dif x),\ \ \mbox{ in law}.
$$
In Hauray and Mischler's work \cite{HM14}, various quantitative and qualitative estimates related to the chaos are obtained for different notions such as 
Kac's chaos, entropy chaos and Fisher information chaos. More references about Kac's chaos can be also found in \cite{HM14}.

\medskip

{\bf Propagation of chaos:} If one considers Kac's chaos as a static version of chaos, then propagation of chaos is usually understood as a dynamical version of Kac's chaos.
More precisely, let $({\boldsymbol\xi^N_t})_{t\geq 0}:=(\xi^{N,1}_t,\cdots,\xi^{N,N}_t)_{t\geq 0}$ be a family of $E^N$-valued continuous stochastic processes, 
which can be thought of as the evolution of $N$-particles. Let $(\xi_t)_{t\geq 0}$ be a limit $E$-valued continuous  stochastic process defined on the same probability space.
Let $\mu^N_t$ be the law of ${\boldsymbol\xi^N_t}$ in $E^N$ and $\mu_t$ be the law of $\xi_t$ in $E$. Suppose that $\mu^N_0$ is $\mu_0$-chaotic at time $0$. One says that
{\it propagation of chaos} holds if for any time $t>0$, $\mu^N_t$ is $\mu_t$-chaotic. Usually, as the evolution of particle distributions, the probability measures
$\mu^N_t$ and $\mu_t$ satisfy some Fokker-Planck equation in the weak sense. Therefore, it can be studied by purely PDE's method.
However, as stochastic processes, one would like to have the following stronger convergence in a probabilistic sense: for each $t>0$,
$$
\lim_{N\to\infty}\mE|\xi^{N,1}_t-\xi_t|=0,
$$
or in the functional path sense
\begin{align}\label{DZ90}
\lim_{N\to\infty}\mE\left(\sup_{s\in[0,t]}|\xi^{N,1}_s-\xi_s|\right)=0.
\end{align}
In fact, when $F$, $\phi$ and $\sigma$ are globally Lipschitz continuous in $x,r$ and uniformly in $t$, McKean \cite{Mc68}  firstly 
established the following result for \eqref{MV2} and \eqref{MV1}:  
for any $T>0$,
\begin{align}\label{DZ91}
\mE\left(\sup_{s\in[0,T]}|X^{N,1}_s-X_s|^2\right)\leq \frac{C(b,\sigma, T)}{N},
\end{align}
where the constant $C(b,\sigma, T)>0$ can be estimated explicitly. We note that the power of convergence rate $N$ is sharp.
The above estimate is also reproven by Sznitman \cite{Sz} by more direct synchronous coupling methods. 
Since then, propagation of chaos has undergone an enormous development in mathematical kinetic theory (see \cite{TT96, MM13, HM14, JW18}). 
Moreover, propagation of chaos  also appears in many other disciplines including data science \cite{Del13}, mean-field games \cite{CD18a, CD18b} 
and the training of neural networks \cite{RV19}, etc.
In a recent paper \cite{CD21}, Chaintron and Diez reviews various models, methods as well as applications for propagation of chaos.

\medskip

Obviously, Lipschitz assumptions on $F,\phi$ and $\sigma$ are too strong in practice. In fact, 
most of the interesting physical models have bounded measurable or even singular interaction kernels.
For examples, the rank-based interaction diffusion studied in \cite{Sh12, La18} has a discontinuous interaction kernel (see \eqref{AB1} below), and
the Biot-Savart law appearing in the vortex description of 2d imcompressible Navier-Stokes equations has a singular kernel like $x^\perp/|x|^2$.
For this type of singular kernels, Osada \cite{Os87} firstly showed the propagation of chaos for the point vortices associated with the 2d Navier-Stokes equation with large viscosity. 
Recently, in \cite{FHM14}, Fournier, Hauray, and Mischler dropped the assumption of large viscosity  by the classical martingale method. More recently,
 Jabin and Wang \cite{JW18} obtained a first quantitative convergence rate about the relative entropy between the law of particle system and the tensorized limit law, where the key point is an estimate for the entropy and a large deviation type exponential functional.
In fact, the results in \cite{JW18} can be applied to a large class of singular kernels $K$ in $W^{-1,\infty}$ with $K(x)=-K(-x)$, 
as well as the bounded measurable interaction kernel (see Section \ref{55}).
We note that the proof in \cite{JW18} strongly depends on the symmetry of the kernel $K(x)$, not valid for general $L^p$-singular kernel.

\medskip

For general $L^p$-singular interaction kernels, in \cite{To20}, Toma$\rm\check{s}$evi\'c uses the partial Girsanov transform as in \cite{JTT18} to derive the
propagation of chaos under the extra assumption that the set of discontinuous points of the interaction kernel has Lebesgue measure zero. In \cite{HHMT}, 
Hoeksema, Holding, Maurelli and Tse showed a large deviation result for a particle system with $L^p$-singular interaction kernels.
As a byproduct, they also obtained a result of propagation of chaos (see also \cite{La18}). However, in \cite{To20} and \cite{HHMT},  both of them assume 
the initial distributions of the particle system are i.i.d, that is, the initial distributions are not really {\it chaotic}. This assumption is crucial for them to construct a weak solution
for the interaction particle system by Girsanov's transform. 
In the present paper we overcome this difficulty by showing the existence of strong solutions for the particle system (see Lemma \ref{Le52} below), and then
obtain the strong convergence as in \eqref{DZ90} for singular interaction kernels and the quantitative convergence \eqref{DZ91} for bounded measurable kernels
by Zvonkin's transformation. Note that Bao and Huang \cite{BH19} have already used the Zvonkin transformation to obtain propagation of chaos for H\"older interaction kernels
with non-optimal rate $N^{-1/4}$.

\subsection{Main results}
Before stating our main assumptions, we introduce the following index sets:
\begin{align}\label{Indo}
\sI^o:=\big\{(q,\bbp)\in(2,\infty)^{1+d}: |\tfrac1\bbp|+\tfrac 2q<1\big\}
\end{align}
and 
$$
\sX:=\big\{\x=(x_{i_1},\cdots,x_{i_d}):\mbox{any permutation of $(x_1,\cdots,x_d)$}\big\}.
$$
Now we make the following main assumptions:
\begin{enumerate}[{\bf (H$^{\sigma}$)}]
\item There are $\kappa_0\geq 1$ and $\gamma_0\in(0,1]$ such that for all $t\geq 0$ and $x,x',\xi\in\mR^d$,
\begin{align}\label{SIG}
\qquad\quad \kappa_0^{-1}|\xi|\leq |\sigma(t,x)\xi|\leq \kappa_0|\xi|,\ 
\|\sigma(t,x)-\sigma(t,x')\|_{HS}\leq \kappa_0|x-x'|^{\gamma_0},
\end{align}
where $\|\cdot\|_{HS}$ is the usual Hilbert-Schmidt norm of a matrix.
Moreover, for some $(q_0,\bbp_0)\in\sI^o$ and $\x_0\in\sX$ and any $T>0$,
\begin{align}\label{SIG0}
\|\nabla\sigma\|_{\mL^{q_0}_T(\wt\mL^{\bbp_0}_{\x_0})}\leq\kappa_0,
\end{align}
where the localized space $\wt\mL^{\bbp}_{\x}$ is defined in Subsection \ref{Sec21} below.
\end{enumerate}
\begin{enumerate}[{\bf (H$^b$)}]
\item 
Suppose that $\phi_t(x,x)=0$ and for some measurable $h:\mR_+\times\mR^d\to\mR_+$ and $\kappa_1>0$,
\begin{align}\label{CC1}
|F(t,x,r)|\leq h(t,x)+\kappa_1|r|,\ \ |F(t,x,r)-F(t,x,r')|\leq \kappa_1|r-r'|,
\end{align}
and for some $(q,\bbp)\in\sI^o$ and $\x\in\sX$ and for any $T>0$, 
\begin{align}\label{CC2}
\nor h\nor_{\mL^{q}_T(\wt\mL^{\bbp}_{\x})}
+\left[\int^T_0\sup_{y\in\mR^d}\left(\nor\phi_t(\cdot,y)\nor^{q}_{\wt\mL^{\bbp}_{\x}}
+\nor\phi_t(y,\cdot)\nor^{q}_{\wt\mL^{\bbp}_{\x}}\right)\dif t\right]^{\frac1q}\leq\kappa_1.
\end{align}
\end{enumerate}

\noindent {\bf Example 1.} We provide two examples to illustrate condition \eqref{CC2}.
\begin{enumerate}[(i)]
\item Let $d\geq 2$ and $\phi_t(x,y)=c_t(x,y)/|x-y|^\alpha$, where $c_t(x,y)$ is bounded measurable and $\alpha\in(0,1)$.
It is easy to see that \eqref{CC2} holds for $q$ close to $\infty$ and $p\in(d,\frac d\alpha)$ with $\frac dp+\frac2q<1$.

\item Let $d\geq 1$ and $\phi_t(x,y)=c_t(x,y)/\Pi_{i=1}^d|x_i-y_i|^{\alpha_i}$, where $\alpha_i\in(0,\frac12)$ satisfies $\alpha_1+\cdots+\alpha_d<1$ 
and $c_t(x,y)$ is bounded measurable. Note that one can choose $q$ close to $\infty$ and $p_i>2$ close to $1/\alpha_i$ 
so that $|\frac 1\bbp|+\frac2q<1$ and \eqref{CC2} holds. In this case, the kernel is allowed to have singularities along each axis.
\end{enumerate}

Throughout this paper we use $\Theta$ to denote the set of parameters that a constant may depend on. $\Theta$ may have different parameters in different occasions, which should be clear from the context, e.g.,
$$
\Theta=(m,d,\gamma_0,\kappa_0,\kappa_1, q_0,\bbp_0,q,\bbp,\cdots).
$$
The aim of this paper is to show the following strong convergence of the particle approximation.
\bt\label{Th11}
Let $T>0$. Under {\bf (H$^{\sigma}$)} and {\bf (H$^b$)}, for any  initial values $\bX^N_0$ and $X_0$, 
there are unique strong solutions $\bX^N_t$ and $X_t$ to particle system \eqref{MV2} and DDSDE \eqref{MV1}, respectively. 
Moreover, letting $\mu^N_0$ be the law of  $\bX^N_0$ in $\mR^{dN}$ and $\mu_0$  the law of $X_0$ in $\mR^d$,
we have the following strong convergence results:
\begin{enumerate}[(i)]
\item {\bf (Singular kernel)} Suppose that $\mu^N_0$ is symmetric and $\mu_0$-chaotic, and
$$
\lim_{N\to\infty}\mE|X^{N,1}_0-X_0|^2=0. 
$$ 
Then for any $\gamma\in(0,1)$,
\begin{align}\label{LIM1}
\lim_{N\to\infty}\mE\left(\sup_{t\in[0,T]}|X^{N,1}_t-X_t|^{2\gamma}\right)=0.
\end{align}
\item {\bf (Bounded kernel)} If $h$ and $\phi$ in  {\bf (H$^b$)} are bounded measurable and
\begin{align}\label{EN1}
\kappa_2:=\sup_N\cH\(\mu^N_0|\mu^{\otimes N}_0\)<\infty,
\end{align}
where $\mu^{\otimes N}_0\in\cP((\mR^d)^N)$ is the $N$-tensor of $\mu_0$ 
and $\cH$ stands 
for the relative entropy (see \eqref{Rela} below),
then for any $\delta>2$ and $\gamma\in(0,1)$, there are constants $C_i=C_i(T,\gamma,\delta,\Theta)>0$, $i=1,2$ independent of $\phi$ and $\kappa_2$ such that
\begin{align}\label{LIM2}
\mE\left(\sup_{t\in[0,T]}|X^{N,1}_t-X_t|^{2\gamma}\right)\leq C_1\e^{C_2\|\phi\|_\infty^\delta}\left(\mE|X^{N,1}_0-X_0|^2+
\frac{\kappa_2+1}{N}\right)^\gamma.
\end{align}
\end{enumerate}
\et
\br\rm
If $\sup_N\mE|X^{N,1}_0|^p<\infty$ for some $p>2$, then by interpolation one in fact has
$$
\lim_{N\to\infty}\mE\left(\sup_{t\in[0,T]}|X^{N,1}_t-X_t|^{p\gamma}\right)=0,\ \ \gamma\in(0,1).
$$
The Euler approximation for particle system \eqref{MV2} with bounded interaction kernel was studied in \cite{Zh1}, which combined with
\eqref{LIM2} implies the full discretization approximation for DDSDE \eqref{MV1}.
\er
\noindent {\bf Example 2.}
Let $d=1$. Consider the following rank-based interaction:
\begin{align}\label{AB1}
b(t,x,\mu)=F(t,x,\mu(-\infty,x]).
\end{align}
In this case, the interaction kernel is $\phi(x,y)=\b1_{(-\infty,x]}(y)=\b1_{x-y\geq 0}$, which is bounded and discontinuous.
Thus, by \eqref{LIM2} we have the strong convergence rate of the particle approximation. 
In particular, if we let $V(x):=\mu((-\infty,x])$, $\sigma(t,x)=\sqrt{2}$ and $F(t,x,r)=g(r)$, then $V$ solves the following Burgers type equation:
$$
\p_t V=\Delta V+\left(\int^V_0 g(r)\dif r\right)'.
$$
For $g(r)=r$, this is the classical Burgers equation. {In this way}, the above Burgers type equation has been studied in 
\cite{BT96, Jo97, La18}. In the following Example 3, we have another way to simulate Burgers equation via moderate interaction particle system.

\medskip

Next we turn to the moderate interaction system \eqref{MV30} and have the following result.
\bt\label{Th12}
Let $T>0$. Suppose that {\bf (H$^{\sigma}$)} holds, and 
\begin{align}\label{CN1}
|F(t,x,r)|\leq\kappa_1,\ \ |F(t,x,r)-F(t,x,r')|\leq \kappa_1|r-r'|,
\end{align}
and for $\eps_N\in(0,1)$ with $\eps\to 0$ as $N\to\infty$,
$$
\phi_t(x,y)=\phi_{\eps_N}(x-y)=\eps_N^{-d}\phi((x-y)/\eps_N),
$$ 
where $\phi$ is a bounded probability density function in $\mR^d$ with support in the unit ball. Then for any initial value $X_0$ with bounded density $\rho_0$,
there is a unique strong solution $X$ to density-dependent SDE \eqref{MV3} such that for each $t>0$, $X_t$ admits a density $\rho_t$ with
\begin{align}\label{XM1}
\|\rho_t\|_\infty\leq C(T,\Theta)\|\rho_0\|_\infty,\ t\in[0,T]. 
\end{align}
Moreover, 
under \eqref{EN1}, for any $T>0$, $\beta\in(0,\gamma_0)$, $\gamma\in(0,1)$ and $\delta>2$,
there are constants $C_i=C_i(T,\beta,\gamma,\delta,\Theta)>0$, $i=1,2,3$ such that for all $N\geq 2$,
\begin{align}\label{ES08}
\mE\left(\sup_{t\in[0,T]}|X^{N,1}_t-X_t|^{2\gamma}\right)\leq C_1\e^{C_2\eps_N^{-\delta d}}\left(\mE|X^{N,1}_0-X_0|^2
+\frac{\kappa_2 +1}{N}\right)^\gamma+C_3\eps_N^{2\beta\gamma}.
\end{align}
\et
\br\rm
Suppose that for some $C>0$,
$$
\mE|X^{N,1}_0-X_0|^2\leq C/N.
$$ 
If one chooses $\eps_N=C_4/(\ln N)^{1/(\delta d)}$ with $C_4$ being large enough, then by \eqref{ES08}, for some $C>0$,
$$
\mE\left(\sup_{t\in[0,T]}|X^{N,1}_t-X_t|^{2\gamma}\right)\leq \frac{C}{(\ln N)^{(2\beta\gamma)/(\delta d)}}.
$$
In \cite{JM98}, under smoothness assumptions on  $F$, $\phi$ and the initial density $\rho_0$,  Jourdain and M\'el\'eard \cite[Theorem 2.7]{JM98} have proven a similar estimate as \eqref{ES08}.
We note that the concept of moderately interacting particles was introduced by Oelschl\"ager  in \cite{Oe85}. Therein, $\eps_N=N^{-\beta/d}$ and $\beta\in(0,1)$.
For $\beta=0$ and $\beta=1$, they are called weakly and strongly interacting, since they correspond to the scaling order $1/N$ and $1$, respectively. While, the moderate interaction refers to any choice of $\eps_N$ with that $\eps_N\to 0$ and $\eps_N^{-d}/N=o(1)$.
\er

Although we assume that $F$ is bounded in \eqref{CN1}, once we can establish the existence of bounded solutions to the Fokker-Planck equation \eqref{MV4} under linear growth assumptions of $F$ in $r$,
then the boundedness of $F$ in \eqref{CN1} is no longer a restriction. We illustrate this in the following example.

\medskip

\noindent{\bf Example 3.} Consider the following special case:
$$
\p_t\rho=\Delta\rho+\div (F(\rho)\rho),
$$
where $F:\mR_+\to\mR^d$ satisfies $\sum_{i=1}^d|F'_i(r)|\leq\kappa_1$. 
Since the above equation can be written in the following transport form:
$$
\p_t\rho=\Delta\rho+ (F(\rho)+F'(\rho)\rho)\cdot\nabla\rho,
$$
it is easy to see that by the maximum principle,
$$
\|\rho_t\|_\infty\leq\|\rho_0\|_\infty.
$$
This can be established rigorously by considering the truncated $F$ as 
$F_n(r)=F(r\wedge n)$, where $n>\|\rho_0\|_\infty$.
In particular, the above example covers the one dimensional Burgers equation, i.e., $F(r)=r$.
In this case, if one takes $\phi(x)=\b1_{[-1,1]}(x)/2$ in \eqref{MV30}, then 
$$
(\phi_{\eps_N}\circledast\eta_{\bX^N_t})(X^{N,i}_t)=\frac1{2N\eps_N}\sum_{j=1}^N\b1_{|X^{N,i}_t-X^{N,j}_t|\leq\eps_N}.
$$
We believe that this is useful for numerical experiments.

\subsection{Structure of the paper}
In Section 2 we first introduce the localized mixed $L^\bbp$-spaces and its basic properties used in this paper. Then we study second order parabolic PDEs with mixed $L^\bbp$-drifts
and show the unique existence of strong solutions. Since each component of the drift may be in a different mixed $L^\bbp$-space, 
the new point here is that the second order derivative of the solution shall stay in a direct sum space (see Theorem \ref{Pre3:Sch0}).

In Section 3, we show the weak and strong well-posedness for stochastic differential equations with mixed $L^\bbp$-drifts. 
As usual, we need to prove a priori Krylov estimates based on the PDE estimates obtained in Section 2, and then show that we can perform the Zvonkin transformation.
Since Zvonkin's transformation is a $C^1$-diffeomorphism and reduces the original singular SDE to an equivalent regular SDE, one can use well-known 
results such as the heat kernel estimates to derive some apriori estimates for the original SDE, and then show our main results. 
We emphasize that the mixed $L^\bbp$-space is not invariant under $C^1$-diffeomorphism transformation. Thus one can not obtain the Krylov estimate 
directly through the transformed equation. Instead, we use the heat kernel estimates to show the Krylov estimate for the indices $(q,\bbp)\in\sI_2$.

In Section 4, by Picard's iteration, we show the weak and strong well-posedness for {distribution density-distribution dependent SDEs} with mixed  $L^\bbp$-drifts,
where we use the entropy formula, Pinsker's inequality and the Fokker-Planck equation to show that the Picard iteration of the density is a Cauchy sequence
in $L^1\cap L^\infty$.

In Section 5, by the classical martingale method we show that the propagation of chaos for systems as in \eqref{MV2} with singular kernels
holds in the weak sense, where the key point is to use the partial Girsanov transform used in \cite{JTT18, To20} to 
derive some uniform estimate for the exponential functional. In particular, the strong solution is used to treat the chaos of the initial distributions.
Moreover, we also provide a detailed proof for Jabin and Wang's quantitative result \cite{JW18} for bounded interaction kernels. This is not new and only for  the readers' convenience.

In Sections 6 and 7, we give the proofs of Theorems \ref{Th11} and \ref{Th12}, and
show how to use Zvonkin's transformation again to derive the strong convergence from the weak convergence obtained in Section 5, where the key point is Lemma \ref{Le61}.

We conclude this introduction by introducing the following convention: 
Throughout this paper, we use $C$ with or without subscripts to denote constants, whose values
may change from line to line. We also use $:=$ to indicate a definition and $a^+:=0\vee a$. By $A\lesssim_C B$ and $A\asymp_C B$
or simply $A\lesssim B$ and $A\asymp B$, we mean that for some constant $C\geq 1$,
$$
A\leq C B,\ \ C^{-1} B\leq A\leq CB.
$$

\section{Preliminaries}

\subsection{Mixed $L^p$-spaces}\label{Sec21}
In this section we recall the definition of localized mixed $L^\bbp$-spaces, which was originally introduced in \cite{BP61}.
As we have seen in the introduction, these are very suitable for singular interacting particle system (see also \cite{HHMT}).
Let $d\in\mN$. For a multi-index $\bbp=(p_1,\cdots,p_d)\in(0,\infty]^d$
and any  permutation $\x\in\sX$, the mixed $\mL^\bbp_\x$-space is defined by
\begin{align}\label{AM9}
\|f\|_{\mL^\bbp_\x}:=\Bigg(\int_\mR\Bigg(\int_\mR\cdots\left(\int_\mR |f(x_1,\cdots,x_d)|^{p_d}\dif x_{i_d}\right)^{\frac{p_{d-1}}{p_d}}\cdots
\dif x_{i_2}\Bigg)^{\frac{p_1}{p_2}}\dif x_{i_1}\Bigg)^{\frac{1}{p_1}}.
\end{align}
When $\bbp=(p,\cdots,p)\in(0,\infty]^d$, the mixed
$\mL^\bbp_\x$-space is the usual $L^p (\R^d)$-space, simply denoted by $\mL^p$. Note that for general $\x\not=\x'$ and $\bbp\not=\bbp'$,
$$
\mL^{\bbp'}_\x\not=\mL^{\bbp}_\x\not=\mL^{\bbp}_{\x'}.
$$
For multi-indices $\bbp, \bbq\in(0,\infty]^d$, we shall use the following notations:
$$
\frac1{\bbp}:=\Big(\frac1{p_1},\cdots,\frac1{p_d}\Big),
\quad  \bbp\cdot\bbq:=\sum_{i=1}^dp_iq_i,\ \ \Big|\frac1\bbp\Big|=\sum_{i=1}^d\frac1{p_i},
$$
and 
$$
\bbp>\bbq\ \ (\hbox{\rm resp.  } \bbp\geq\bbq;\ \bbp=\bbq)\Longleftrightarrow p_i>q_i
 \ (\hbox{\rm resp.  }  p_i\geq q_i;\ p_i=q_i) \ \hbox{ for all } i=1,\cdots,d.
$$
Moreover, we use bold numbers to denote constant vectors in $\mR^d$, for example, 
$$
\boldsymbol{1}=(1,\cdots,1),\ \ \boldsymbol{2}=(2,\cdots,2).
$$
For multi-indices $\bbp,\bbq,\bbr\in(0,\infty]^d$ with $\tfrac1\bbp+\tfrac1\bbr=\tfrac1\bbq$, the following H\"older inequality holds
\begin{align}\label{HH1}
\|fg\|_{\mL^\bbq_\x}\leq \|f\|_{\mL^\bbp_\x}\|g\|_{\mL^\bbr_\x}.
\end{align}
For any multi-indices $\bbp,\bbq,\bbr\in[1,\infty]^d$ with $\tfrac1\bbp+\tfrac1\bbr=\boldsymbol{1}+\tfrac1\bbq$, the following Young inequality holds
\begin{align}\label{HH2}
\|f*g\|_{\mL^\bbq_\x}\leq \|f\|_{\mL^\bbp_\x}\|g\|_{\mL^\bbr_\x}.
\end{align}

For any $r>0$, let $B_z^r$ be the ball in $\mR^d$ with radius $r$ and center $z$.
Let $\chi:\mR^d\to[0,1]$ be a smooth cutoff function with $\chi|_{B_1}=1$
and $\chi|_{B_2^c}=0$. For fixed $r>0$, we set
\begin{align}\label{Cut}
\chi^r_z(x):=\chi((x-z)/r),\ \ x,z\in\mR^d.
\end{align}
For $\bbp\in[1,\infty]^d$, we introduce the following localized $L^\bbp$-space (see \cite{Zh-Zh2}):
\begin{align*}
\widetilde{\mL}^\bbp_\x:=\Big\{f\in L^1_\text{\rm loc}(\mR^d),
\nor f\nor_{\wt\mL^\bbp_\x}:=\sup_z\|\chi^r_zf\|_{\mL^\bbp_\x}<\infty\Big\},
\end{align*}
and for a finite time interval ${\rm I}\subset\mR$ and $q\in[1,\infty]$,
\begin{align}\label{LOC1}
\wt\mL^q_{\rm I}(\wt\mL^\bbp_\x):=\Big\{f\in L^1_\text{\rm loc}({\rm I}\times\mR^d),\nor f\nor_{\wt\mL^q_{\rm I}(\wt\mL^\bbp_\x)}
:=\sup_z\|\chi^r_zf\|_{\mL^q_{\rm I}(\mL^\bbp_\x)}<\infty\Big\},
\end{align}
where for a Banach space $\mB$ we set
$$
\mL^q_{{\rm I}}(\mB):=L^q({\rm I};\mB).
$$ 
By a finitely covering technique, it is easy to see that the definitions of $\widetilde{\mL}^\bbp_\x$ and $\wt\mL^q_{\rm I}(\wt\mL^\bbp_\x)$ do not depend on the choice of 
$r$ (see \cite{Zh-Zh2}), and for any $1\leq q_2\leq q_1\leq\infty$ and $1\le \bbp_2\le \bbp_1\le\infty$,
\begin{align}\label{Pre1:emb}
\widetilde{\mL}^{\bbp_1}_\x\subset\widetilde{\mL}^{\bbp_2}_\x,\ \ \wt\mL^{q_1}_{\rm I}(\wt\mL^{\bbp_1}_\x)\subset\wt\mL^{q_2}_{\rm I}(\widetilde{\mL}^{\bbp_2}_\x).
\end{align}
This property is the main advantage of using localized spaces.
Since the supremum $z$  in the definition of $\wt\mL^q_{\rm I}(\wt\mL^\bbp_\x)$ is taken outside the time integral, we obviously have
$$
\mL^q_{\rm I}(\wt\mL^\bbp_\x)\subset\wt\mL^q_{\rm I}(\wt\mL^\bbp_\x).
$$
Moreover, for $\alpha\geq 0$, let $\cC^\alpha$ be the usual H\"older space with norm:
$$
\|f\|_{\cC^\alpha}:=\sum_{j=0}^{[\alpha]}\|\nabla^jf\|_\infty+\sup_{x\not=y\in\mR^d}\frac{|\nabla^{[\alpha]}f(x)-\nabla^{[\alpha]}f(y)|}{|x-y|^{\alpha-[\alpha]}},
$$
where $\nabla^j$ stands for the $j$-order gradient and $[\alpha]$ stands for the integer part of $\alpha$.
For simplicity we write
\begin{align*}
\wt\mL^q_T(\wt\mL^\bbp_\x):=\wt\mL^q_{[0,T]}(\wt\mL^\bbp_\x),~~\mL^p_T:=\mL^p_{[0,T]}(\mL^p),\
\mL^\infty_T(\cC^\alpha):=\mL^\infty_{[0,T]}(\cC^\alpha).
\end{align*}
{\bf Example.}
For $i=1,\cdots,d$ and $\alpha\in(0,1)$, let $f_i(x)=b(x)|x_i|^{-\alpha}$, where $b(x)$ is a bounded measurable function. It is easy to see that
$f_i\in\wt\mL^\bbp_{\x_i}$, where $\x_i=(x_1,\cdots,x_{i-1},x_{i+1}, \cdots,x_d,x_i)$ and $\bbp=(\infty,\cdots,\infty,p)$ with $p\in(1,\frac1\alpha)$.
From this example, one sees that for a $C^1$-diffeomorphism $\Phi$ from $\mR^d$ to $\mR^d$,
say $\Phi(x)=(x_i,x_1,\cdots,x_{i-1},x_{i+1},\cdots,x_d)$, it may happen that
$$
f_i\circ\Phi\notin\wt\mL^\bbp_{\x_i}.
$$

Throughout this paper, we shall use the same notation $\varGamma_\eps$ to denote mollifiers in various dimensions $N$, i.e.,
\begin{align}\label{Mol}
\varGamma_\eps(x)=\eps^{-N}\varGamma(x/\eps),\ \eps\in(0,1),
\end{align}
where $\varGamma$ is a nonnegative smooth density function in $\mR^N$ with compact support in the unit ball.  For a function $f\in L^1_{\rm loc}(\mR^N)$, 
the mollifying approximation of $f$ is defined by
$$
f_\eps(x):=f*\varGamma_\eps(x)=\int_{\mR^N}f(x-y)\varGamma_\eps(y)\dif y.
$$
The dimension $N$ takes different values in different occasions, which should be clear from the respective context.

The following lemma is obvious by the definitions.
\bl
For any $f\in\wt\mL^\bbp_\x$, there is a constant $C=C(\bbp)>0$ such that for all $\eps\in(0,1)$,
\begin{align}\label{BZ130}
\nor f_\eps\nor_{\wt\mL^\bbp_\x}\leq C\nor f\nor_{\wt\mL^\bbp_\x},
\end{align}
and for any $R>0$,
\begin{align}\label{BZ30}
\lim_{\eps\to 0}\| (f_\eps-f)\chi^R_0\|_{\mL^\bbp_\x}=0.
\end{align}
\el
The local Hardy-Littlewood maximal function in $\mR^d$ is defined by
\begin{align*}
\cM f(x):=\sup_{r\in(0,1)}\frac{1}{|B^r_0|}\int_{B^r_0}f(x+y)\dif y.
\end{align*}
The following result is taken from Lemma 2.1 in \cite{Xi-Xi-Zh-Zh}.

\bl\label{Le21}
\begin{enumerate}[(i)]
\item There is a constant $C=C(d)>0$, such that for any $f\in L^\infty(\mR^d)$ with $\nabla f\in L^1_\text{\rm loc}(\mR^d)$, 
\begin{align}\label{Pre1:Mnf}
|f(x)-f(y)|\le C|x-y|\(\cM|\nabla f|(x)+\cM|\nabla f|(y)+\|f\|_\infty\)
\end{align}
for Lebesgue-almost all $x,y\in\mR^d$.
\item For any $(q,\bbp)\in(1,\infty)^{1+d}$, there is a $C=C(d,p,q)>0$ such that for all $f\in\wt\mL^q_T(\wt\mL^\bbp_\x)$,
\begin{align}\label{Pre1:Mix}
\nor \cM f\nor_{\wt\mL^q_T(\wt\mL^\bbp_\x)}\le C\nor f\nor_{\wt\mL^q_T(\wt\mL^\bbp_\x)}.
\end{align}
\end{enumerate}
\el
\subsection{A study of PDEs with mixed $L^\bbp$-drifts}
In this section we show the existence and uniqueness of strong solutions in the PDE sense to second order parabolic PDEs with drifts in mixed $L^p$-spaces.
For $t>0$, let $P_tf(x)=\mE f(x+W_t)$ be the Gaussian heat semigroup, i.e.,
$$
P_tf(x)=\int_{\mR^d} g_t(x-y)f(y)\dif y,
$$
where
$$
 g_t(x):=(2\pi t)^{-\frac d2}\e^{-\frac{|x|^2}{2t}}.
$$ 
First of all, we establish the following easy estimates about $P_t$.
\bl\label{Le23}
\begin{enumerate}[(i)]
\item For any $\bbp\in(1,\infty)^d$, $T>0$ and $\beta\geq 0$, there is a constant $C=C(T,\bbp,\beta, d)>0$ such that
 for all $f\in\mL^\bbp_\x$ and $t\in(0,T]$,
\begin{align}\label{SA1}
\|P_tf\|_{\cC^\beta}\leq Ct^{-\frac12(\beta+|\frac1\bbp|)}\|f\|_{\mL^\bbp_\x}.
\end{align}
\item For any $\bbq\geq\bbp$, there is a constant $C=C(\bbq,\bbp,d)>0$ such that 
for all $f\in\mL^\bbp_\x$ and $t>0$,
\begin{align}\label{SA2}
\|\nabla P_tf\|_{\mL^\bbq_\x}\leq Ct^{-\frac12(1+|\frac1\bbp|-|\frac1\bbq|)}\|f\|_{\mL^\bbp_\x}.
\end{align}
\end{enumerate}
\el
\begin{proof}
(i) Note that for $m=0,1,\cdots$,
$$
\nabla^m P_tf(x)=\int_{\mR^d}\nabla^m g_t(x-y)f(y)\dif y.
$$
For $\frac1\bbq+\frac1\bbp={\bf 1}$, by H\"older's inequality \eqref{HH1} and the scaling, we have
$$
\|\nabla^m P_tf\|_{\infty}
\leq \|\nabla^m g_t\|_{\mL^{\bbq}_\x}\|f\|_{\mL^{\bbp}_\x}=t^{-\frac12(m+|\frac1\bbp|)}\|\nabla^m g_1\|_{\mL^{\bbq}_\x}\|f\|_{\mL^{\bbp}_\x},
$$
where $\|\nabla^m g_1\|_{\mL^{\bbq}_\x}<\infty$.
Then estimate \eqref{SA1} follows by the interpolation theorem for H\"older spaces.

(ii) For $\bbr\in[1,\infty]^d$ with $\tfrac1\bbp+\tfrac1\bbr=\boldsymbol{1}+\tfrac1\bbq$, 
by Young's inequality \eqref{HH2}  and the scaling, 
we have
$$
\|\nabla P_tf\|_{\mL^\bbq_\x}
\leq \|\nabla g_t\|_{\mL^{\bbr}_\x}\|f\|_{\mL^{\bbp}_\x}=t^{-\frac12(1+d-|\frac1\bbr|)}\|\nabla g_1\|_{\mL^{\bbr}_\x}\|f\|_{\mL^{\bbp}_\x}.
$$
Then estimate \eqref{SA2} follows because $\|\nabla g_1\|_{\mL^{\bbr}_\x}<\infty$.
\end{proof}
We introduce the following index sets for later use:
\begin{align}\label{Ind}
\sI_m:=\Big\{(q,\bbp)\in(1,\infty)^{1+d}:\ |\tfrac1\bbp|+\tfrac2q<m\Big\},\ m=1,2.
\end{align}
\br\rm 
$\sI^o\subset\sI_1$, where $\sI^o$ is defined by \eqref{Indo}. For $(q,\bbp)\in\sI^o$, it holds that $(\frac q2,\frac\bbp2)\in\sI_2$.
\er
For $\lambda\geq 0$ and $f\in\mL^q_T(\mL^\bbp_\x)$, we define
$$
u(t,x):=\int^t_0\e^{-\lambda(t-s)}P_{t-s}f(s,x)\dif s,\ t>0,
$$
which solves the following non-homeogenous heat equation
$$
\p_t u=\tfrac12\Delta u-\lambda u+f,\ \ u(0)=0.
$$
\bl
\begin{enumerate}[(i)]
\item For any $T>0$, $(q,\bbp)\in\sI_2$ and $\beta\in[0,2-|\frac1\bbp|-\frac2q)$, 
there is a constant $C=C(T,d,q,\bbp,\beta)>0$ such that for all $\lambda\geq 0$,
\begin{align} \label{AC3}
(1\vee\lambda)^{\frac12(2-\beta-|\frac1\bbp|-\frac2q)}\|u\|_{\mL^\infty_T(\cC^\beta)}\leq C\|f\|_{\mL^q_T(\mL^\bbp_\x)}.
\end{align}
\item For any $T>0$, $(q,\bbp)\in\sI_2$ and $(q',\bbp')\geq(q,\bbp)$ with
$|\frac1\bbp|+\frac2q<|\frac1{\bbp'}|+\frac2{q'}+1$, there is a constant $C=C(T,d,q,\bbp,q',\bbp')>0$ such that for all $\lambda\geq 0$,
\begin{align}
(1\vee\lambda)^{\frac12(1+|\frac1{\bbp'}|+\frac2{q'}-|\frac1\bbp|-\frac2q)}\|\nabla u\|_{\mL^{q'}_T(\mL^{\bbp'}_\x)}\leq C \|f\|_{\mL^q_T(\mL^\bbp_\x)}.\label{BX9}
\end{align}
\item For any $T>0$, $(q,\bbp)\in\sI_1$ and $\lambda\geq 0$, there is a constant $C=C(\lambda,T,d,q,\bbp)>0$ such that for all $0\leq t_0<t_1\leq T$,
\begin{align}\label{AC4}
\|u(t_1)-u(t_0)\|_\infty\leq C(t_1-t_0)^{\frac12}\|f\|_{\mL^q_T(\mL^\bbp_\x)}.
\end{align}
\end{enumerate}
\el
\begin{proof}
(i) For $\beta\in[0,2-|\frac1\bbp|-\frac2q)$, by \eqref{SA1} and H\"older's inequality in the time variable, we have
\begin{align*}
\|u(t)\|_{\cC^\beta}&\lesssim \int^t_0\e^{-\lambda(t-s)}(t-s)^{-\frac12(\beta+|\frac1\bbp|)}\|f(s)\|_{\mL^\bbp_\x}\dif s\\
&\leq \left(\int^t_0\big(\e^{-\lambda s}s^{-\frac12(\beta+|\frac1\bbp|)}\big)^{\frac{q}{q-1}}\dif s\right)^{1-\frac{1}{q}}\|f\|_{\mL^q_T(\mL^\bbp_\x)}\\
&\lesssim (1\vee\lambda)^{-\frac12(2-\beta-|\frac1\bbp|-\frac2q)}\|f\|_{\mL^q_T(\mL^\bbp_\x)}.
\end{align*}
(ii)  For $(q',\bbp')\geq(q,\bbp)$ with $|\frac1\bbp|+\frac2q<|\frac1{\bbp'}|+\frac2{q'}+1$,
by \eqref{SA2} we have
$$
\|\nabla u(t)\|_{\mL^{\bbp'}_\x}
\lesssim \int^t_0\e^{-\lambda(t-s)}(t-s)^{-\frac12(1+|\frac1\bbp|-|\frac1{\bbp'}|)}\|f(s)\|_{\mL^\bbp_\x}\dif s.
$$
Let $r\geq 1$ be defined by $\frac1{r}=\frac1{q'}+1-\frac1q$. By Young's inequality we further have
\begin{align*}
\|\nabla u\|_{\mL^{q'}_T(\mL^{\bbp'}_\x)}
&\lesssim \left(\int^T_0\e^{-r\lambda s}s^{-\frac r2(1+|\frac1\bbp|-|\frac1{\bbp'}|)}\dif s\right)^{1/r}\|f\|_{\mL^q_T(\mL^\bbp_\x)}\\
&\lesssim (1\wedge\lambda)^{\frac1r-\frac 12(1+|\frac1\bbp|-|\frac1{\bbp'}|)}\|f\|_{\mL^q_T(\mL^\bbp_\x)}.
\end{align*}
(iii)
For $0\leq t_0<t_1\leq T$, by definition we have
\begin{align*}
u(t_1)-u(t_0)&=\int^{t_0}_0\e^{-\lambda(t_1-s)}(P_{t_1-s}-P_{t_0-s})f(s,x)\dif s\\
&\quad+(\e^{-\lambda(t_1-t_0)}-1)\int^{t_0}_0\e^{-\lambda(t_0-s)}P_{t_0-s}f(s,x)\dif s\\
&\quad+\int^{t_1}_{t_0}\e^{-\lambda(t_1-s)}P_{t_1-s}f(s,x)\dif s\\
&=:I_1+I_2+I_3.
\end{align*}
For $I_1$, noting that
$$
\|P_tf-f\|_\infty\leq\frac12\int^t_0\|\Delta P_sf\|_\infty\dif s
\lesssim\left(\int^t_0s^{-\frac12}\dif s\right)\|\nabla f\|_\infty\lesssim t^{\frac12}\|\nabla f\|_\infty,
$$
by \eqref{SA1} and H\"older's inequality, we have
\begin{align*}
\|I_1\|_\infty&\lesssim (t_1-t_0)^{\frac12}\int^{t_0}_0\|\nabla P_{t_0-s}f(s)\|_\infty\dif s\\
&\lesssim (t_1-t_0)^{\frac12}\int^{t_0}_0(t_0-s)^{-\frac12(1+|\frac1\bbp|)}\|f(s)\|_{\mL^\bbp_\x}\dif s\\
&\lesssim  (t_1-t_0)^{\frac12} t_0^{\frac12(1-\frac2q-|\frac1\bbp|)}\|f\|_{\mL^q_T(\mL^\bbp_\x)},
\end{align*}
and because $1-\e^{-\lambda(t_1-t_0)}\leq \lambda(t_1-t_0)$,
$$
\|I_2\|_\infty\lesssim  \lambda (t_1-t_0)\|f\|_{\mL^q_T(\mL^\bbp_\x)}.
$$
For $I_3$, as above, by \eqref{SA1} and H\"older's inequality, we have
$$
\|I_3\|_\infty\lesssim  \left(\int^{t_1-t_0}_0\big(\e^{-\lambda s}s^{-\frac12|\frac1\bbp|}\big)^{\frac{q}{q-1}}\dif s\right)^{1-\frac{1}{q}}\|f\|_{\mL^q_T(\mL^\bbp_\x)}
\leq (t_1-t_0)^{1-\frac12(\frac2q+|\frac1\bbp|)}\|f\|_{\mL^q_T(\mL^\bbp_\x)}.
$$
Combining the above estimates and because $\frac2q+|\frac1\bbp|<1$, we obtain \eqref{AC4}.
\end{proof}

Now we shall study the following second order parabolic PDE in $\mR_+\times\mR^d$:
\begin{align}\label{Pre3:PDE0}
\p_tu=\tr(a\cdot\nabla^2 u)+b\cdot\nabla u-\lambda u+f,\quad u(0)=0,
\end{align}
 where $\lambda\ge0$, $a:=\sigma\sigma^*/2$, $\sigma$ satisfies \eqref{SIG} and 
 $$
 b, f\in L^1_{\rm loc}(\mR_+\times\mR^d).
 $$
We introduce the following notion of solutions to PDE \eqref{Pre3:PDE0}.
\bd\label{Pre3:Def}
Let $T>0$ and $\sU_T\subset L^1_{\rm loc}(\mR_+\times\mR^d)$ be some subclass of locally integrable functions. We call $u\in\sU_T$ 
a solution of PDE \eqref{Pre3:PDE0} 
if for all $t\in[0,T]$ and $\varphi\in C_c(\mR^d)$,
\begin{align*}
\<u(t),\varphi\>=\int_0^t \<\tr(a\cdot\nabla^2 u)+b\cdot\nabla u,\varphi\>\dif s-\lambda \int_0^t \<u,\varphi\>\dif s+
\int_0^t 
\<f,\varphi\>\dif s,
\end{align*}
where we have implicitly assumed that $\nabla^2 u\in L^1_{\rm loc}$ and $\nabla u\in L^\infty_{\rm loc}$ so that the terms on the right hand side are well defined.
Here $\sU_T$ will be specified below in the respective cases.
\ed
We first show the following result for bounded drift $b$ (see \cite[Theorem 2.1]{LX21}).
\bt\label{Pre3:Sch}
Let $T>0$ and $(q,\bbp)\in(1,\infty)^{1+d}$. Suppose that \eqref{SIG} holds and $b$ is bounded measurable. 
Then for any $f\in\wt\mL^q_T(\wt\mL^\bbp_\x)$ and $\beta\in[0,2-|\frac1\bbp|-\frac2q)$,
there exists a unique solution $u\in\sU_T$ in the sense of Definition \ref{Pre3:Def}, where $\sU_T$ consists of all $u$ with 
\begin{align}\label{BB6}
(1\vee\lambda)^{\frac12(2-\beta-|\frac1\bbp|-\frac2q)}\|u\|_{\mL^\infty_T(\cC^\beta)}+
\nor\nabla^2 u\nor_{\wt\mL^q_T(\wt\mL^\bbp_\x)}\le C\nor f\nor_{\wt\mL^q_T(\wt\mL^\bbp_\x)}.
\end{align}
Here and below, the constant $C=C(T,\kappa_0, d,\bbp,q,\beta,\|b\|_{\mL^\infty_T})>0$ is independent of $\lambda$. Moreover, for any $(q',\bbp')\geq(q,\bbp)$ with
$|\frac1\bbp|+\frac2q<|\frac1{\bbp'}|+\frac2{q'}+1$, we also have
\begin{align}\label{BB61}
(1\vee\lambda)^{\frac12(1+|\frac1{\bbp'}|+\frac2{q'}-|\frac1\bbp|-\frac2q)}\nor\nabla u\nor_{\wt\mL^{q'}_T(\wt\mL^{\bbp'}_\x)}\le C\nor f\nor_{\wt\mL^q_T(\wt\mL^\bbp_\x)},
\end{align}
and for all $0\leq t_0\leq t_1\leq T$,
\begin{align}\label{BB62}
\|u(t_1)-u(t_0)\|_\infty\le C(\lambda)(t_1-t_0)^{\frac12}\nor f\nor_{\wt\mL^q_T(\wt\mL^\bbp_\x)}.
\end{align}
\et
\begin{proof}
We only prove the a priori estimates \eqref{BB6}, \eqref{BB61} and \eqref{BB62}. 
The existence is then standard by mollifying the coefficients and a compactness argument. Fix $r>0$.
Let $\chi^r_z$ be the cutoff function in \eqref{Cut} and  $w_z:=u\chi^r_z$. It is easy to see that
\begin{align}\label{BB8}
\p_tw_z=\tr(a\cdot\nabla^2 w_z)-\lambda w_z+g_z,\quad w_z(0)=0,
\end{align}
where
$$
g_z:=\tr(a\cdot\nabla^2 u)\chi^r_z-\tr(a\cdot\nabla^2 w_z)+(b\cdot\nabla u) \chi^r_z+f\chi^r_z.
$$
Let $(q,\bbp)\in(1,\infty)^{1+d}$. By \cite[Theorem 2.1]{LX21}, there is a constant $C=C(T,\kappa_0, d,\bbp,q)>0$ such that 
\begin{align}\label{BBS1}
\|w_z\|_{\mL^\infty_T(\mL^{\bbp}_\x)}+\|\nabla^2 w_z\|_{\mL^q_T(\mL^{\bbp}_\x)}\lesssim_C\|g_z\|_{\mL^q_T(\mL^{\bbp}_\x)}.
\end{align}
On the other hand, we can write \eqref{BB8} as
$$
\p_tw_z=\Delta w_z-\lambda w_z+\tr((a-\mI)\cdot\nabla^2 w_z)+g_z,\quad w_z(0)=0,
$$
and by Duhamel's formula,
$$
w_z(t,x)=\int^t_0\e^{-\lambda(t-s)}P_{t-s}(\tr((a-\mI)\cdot\nabla^2 w_z)+g_z)(s,x)\dif s.
$$
Note that by \eqref{BBS1},
\begin{align}\label{BB91}
\|\tr((a-\mI)\cdot\nabla^2 w_z)+g_z\|_{\mL^q_T(\mL^\bbp_\x)}\lesssim\|g_z\|_{\mL^q_T(\mL^\bbp_\x)}.
\end{align}
For $\beta\in[0,2-|\frac1\bbp|-\frac2q)$, by \eqref{AC3} and \eqref{BB91} we have
\begin{align}\label{BB00}
(1\vee\lambda)^{\frac12(2-\beta-|\frac1\bbp|-\frac2q)}\|w_z\|_{\mL^\infty_T(\cC^\beta)}
&\lesssim \|g_z\|_{\mL^q_T(\mL^\bbp_\x)}.
\end{align}
For $(q',\bbp')\geq(q,\bbp)$ with
$|\frac1\bbp|+\frac2q<|\frac1{\bbp'}|+\frac2{q'}+1$,
by \eqref{BX9} and \eqref{BB91} we have
\begin{align}
(1\vee\lambda)^{\frac12(1+|\frac1{\bbp'}|+\frac2{q'}-|\frac1\bbp|-\frac2q)}\|\nabla w_z\|_{\mL^{q'}_T(\mL^{\bbp'}_\x)}
\lesssim\|g_z\|_{\mL^q_T(\mL^\bbp_\x)}.\label{BX99}
\end{align}
For $0\leq t_0<t_1\leq T$, by \eqref{AC4} and \eqref{BB91} we have
\begin{align}
\|w_z(t_1)-w_z(t_0)\|_\infty\lesssim(t_1-t_0)^{\frac12}\|g_z\|_{\mL^q_T(\mL^\bbp_\x)}.\label{BX909}
\end{align}
Since $\chi^{2r}_z\nabla^j\chi^r_z=\nabla^j\chi^r_z$ for $j=0,1,2$, we have
\begin{align*}
\|g_z\|_{\mL^q_T(\mL^{\bbp}_\x)}
&\lesssim\|\nabla u\nabla\chi^r_z\|_{\mL^q_T(\mL^{\bbp}_\x)}+\|u\nabla^2\chi^r_z\|_{\mL^q_T(\mL^{\bbp}_\x)}
+\|b\|_{\mL^\infty_T}\|\nabla u\chi^r_z\|_{\mL^q_T(\mL^{\bbp}_\x)}\\
&\leq (\|\nabla\chi^r_z\|_\infty+\|b\|_{\mL^\infty_T})\|\nabla u\chi^{2r}_z\|_{\mL^q_T(\mL^{\bbp}_\x)}
+\|\nabla^2\chi^r_z\|_\infty\|u\chi^{2r}_z\|_{\mL^q_T(\mL^{\bbp}_\x)}.
\end{align*}
Substituting this into \eqref{BBS1}, \eqref{BB00}, \eqref{BX99} and \eqref{BX909} and
taking supremum in $z\in\mR^d$, we obtain 
\begin{align}\label{BB5}
\nor u\nor_{\wt\mL^\infty_T(\wt\mL^\bbp_\x)}+\nor\nabla^2 u\nor_{\wt\mL^q_T(\wt\mL^\bbp_\x)}
&\lesssim\nor f\nor_{\wt\mL^q_T(\wt\mL^\bbp_\x)}+\nor \nabla  u\nor_{\wt\mL^q_T(\wt\mL^\bbp_\x)}+\nor u\nor_{\wt\mL^q_T(\wt\mL^\bbp_\x)},
\end{align}
and for $\beta\in[0,2-|\frac1\bbp|-\frac2q)$,
\begin{align}\label{BB77}
(1\wedge\lambda)^{\frac12(2-\beta-|\frac1\bbp|-\frac2q)}\|u\|_{\mL^\infty_T(\cC^\beta)}
\lesssim\nor f\nor_{\wt\mL^q_T(\wt\mL^\bbp_\x)}+\nor \nabla  u\nor_{\wt\mL^q_T(\wt\mL^\bbp_\x)}+\nor u\nor_{\wt\mL^q_T(\wt\mL^\bbp_\x)},
\end{align}
and for $(q',\bbp')\geq(q,\bbp)$ with $|\frac1\bbp|+\frac2q<|\frac1{\bbp'}|+\frac2{q'}+1$,
\begin{align}\label{BB707}
(1\vee\lambda)^{\frac12(1+|\frac1{\bbp'}|+\frac2{q'}-|\frac1\bbp|-\frac2q)}\nor\nabla u\nor_{\wt\mL^{q'}_T(\wt\mL^{\bbp'}_\x)}
\lesssim\nor f\nor_{\wt\mL^q_T(\wt\mL^\bbp_\x)}+\nor \nabla  u\nor_{\wt\mL^q_T(\wt\mL^\bbp_\x)}+\nor u\nor_{\wt\mL^q_T(\wt\mL^\bbp_\x)},
\end{align}
and for all $0\leq t_0<t_1\leq T$,
\begin{align}\label{BB717}
\|u(t_1)-u(t_0)\|_\infty\lesssim(t_1-t_0)^{\frac12}
\Big(\nor f\nor_{\wt\mL^q_T(\wt\mL^\bbp_\x)}+\nor \nabla  u\nor_{\wt\mL^q_T(\wt\mL^\bbp_\x)}+\nor u\nor_{\wt\mL^q_T(\wt\mL^\bbp_\x)}\Big).
\end{align}
Note that by the interpolation inequality, for any $\eps\in(0,1)$,
$$
\nor \nabla  u\nor_{\wt\mL^q_T(\wt\mL^\bbp_\x)}\leq\eps\nor \nabla^2 u\nor_{\wt\mL^q_T(\wt\mL^\bbp_\x)}+C_\eps\nor u\nor_{\wt\mL^q_T(\wt\mL^\bbp_\x)}.
$$
Substituting this into \eqref{BB5} and choosing $\eps$ small enough, we derive that for any $t\in[0,T]$,
$$
\nor u(t)\nor_{\wt\mL^{\bbp}_\x}+\nor\nabla^2 u\nor_{\wt\mL^q_T(\wt\mL^\bbp_\x)}\lesssim\nor f\nor_{\wt\mL^q_T(\wt\mL^\bbp_\x)}
+\left(\int^t_0\nor u(s)\nor^q_{\wt\mL^{\bbp}_\x}\dif s\right)^{1/q}.
$$
By Gronwall's inequality, we get 
$$
\nor u\nor_{\mL^\infty_T(\wt\mL^{\bbp}_\x)}+\nor\nabla^2 u\nor_{\wt\mL^q_T(\wt\mL^\bbp_\x)}\leq C\nor f\nor_{\wt\mL^q_T(\wt\mL^\bbp_\x)},
$$
which together with \eqref{BB77}, \eqref{BB707} and \eqref{BB717} yields \eqref{BB6}, \eqref{BB61} and \eqref{BB62}.
\end{proof}
\br\rm
For any $T,\gamma>0$ and $(q,\bbp)\in\sI_2$, there is a $C=C(T,\gamma,d,q,\bbp)>0$ such that
\begin{align}\label{226}
\sup_x\mE\left(\int^T_0
h(s,x+W_{\gamma s})
\dif s\right)\leq C\nor h\nor_{\wt\mL^q_T(\wt\mL^\bbp_\x)}.
\end{align}
Indeed, let $a=\sqrt{\gamma/2}\mI$, $b=0$, $\lambda=0$ and $f(s,x)=h(T-s,x)$ in PDE \eqref{Pre3:PDE0}. By \eqref{BB6} we have
\begin{align*}
\mE\left(\int^T_0h(s,x+W_{\gamma s})\dif s\right)&=\int^T_0 P_{\gamma(T-s)} f(s,x)\dif s= u(T,x)\lesssim\nor f\nor_{\wt\mL^q_T(\wt\mL^\bbp_\x)}=\nor h\nor_{\wt\mL^q_T(\wt\mL^\bbp_\x)}.
\end{align*}
In particular, once we have the Gaussian type density estimate for SDEs, 
then by \eqref{226}, we immediately have the Krylov estimate as we shall see in Theorem \ref{Main1} below.
\er

Next we consider the drift $b$ being in the mixed $L^p$-space, where each component $b_i$ may lie in a different mixed $L^p$-space.
Thus the second order generalized derivative of $u$ stays in a direct sum space of mixed $L^p$-spaces. The following result seems to be new and is
the cornerstone of studying SDEs with singular mixed $L^p$-coefficients.
\bt\label{Pre3:Sch0}
Let $T>0$. Suppose \eqref{SIG} and for some  $(q_i,\bbp_i)\in\sI_1$ and $\x_i\in\sX$, $i=1,\cdots,d$,
\begin{align}\label{BB0}
\|b_1\|_{\wt\mL^{q_1}_T(\wt\mL^{\bbp_1}_{\x_1})}+\cdots+\|b_d\|_{\wt\mL^{q_d}_T(\wt\mL^{\bbp_d}_{\x_d})}\leq\kappa_1<\infty.
 \end{align}
 Let $\x_0\in\sX$ and $(q_0,\bbp_0)\in\sI_1$. Define
\begin{align}\label{BZ26}
 \vartheta:=1-\max_{i=0,\cdots,d}(|\tfrac1{\bbp_i}|+\tfrac2{q_i}).
 \end{align}
For any $f\in\wt\mL^{q_0}_T(\wt\mL^{\bbp_0}_{\x_0})$ and $\beta\in[0,\vartheta)$, 
there is a constant $C_0=C_0(T,\kappa_0,d,\bbp_i,q_i,\beta)\ge1$ so that for all $\lambda\ge C_0\kappa_1^{2/\vartheta}$, 
there exists a unique solution $u\in\sU_T$ to PDE \eqref{Pre3:PDE0} in the sense of Definition \ref{Pre3:Def}, where $\sU_T$ consists of all $u=u_0+u_1+\cdots+u_d$ with
\begin{align}\label{CC5}
\lambda^{\frac12(\vartheta-\beta)}\|u\|_{\mL^\infty_T(\cC^{1+\beta})}+\sum_{i=0}^d\nor\nabla^2 u_i\nor_{\widetilde{\mL}^{q_i}_T(\wt\mL^{\bbp_i}_{\x_i})}\le 
C_1\nor f\nor_{\wt\mL^{q_0}_T(\wt\mL^{\bbp_0}_{\x_0})},
\end{align}
where $C_1=C_1(T,\kappa_0, d,\bbp_i,q_i,\beta)>0$ is independent of $\lambda$ and $\kappa_1$.
Moreover, for all $0\leq t_0\leq t_1\leq T$,
\begin{align}\label{BB622}
\|u(t_1)-u(t_0)\|_\infty\le C(\lambda)(t_1-t_0)^{\frac12}\nor f\nor_{\wt\mL^{q_0}_T(\wt\mL^{\bbp_0}_{\x_0})}.
\end{align}
\et
\begin{proof}
Again we only show the a priori estimate \eqref{CC5} since then the existence can be shown by a compactness argument.
Let $u=u_0+u_1+\cdots+u_d$, where $u_0$ solves the following PDE:
$$
\p_tu_0=\tr(a\cdot\nabla^2 u_0)-\lambda u_0+f,\quad u_0(0)=0,
$$
and for each $i=1,\cdots,d$, $u_i$ solves
$$
\p_tu_i=\tr(a\cdot\nabla^2 u_i)+b_i\cdot\p_i u-\lambda u_i,\quad u_i(0)=0.
$$
Let $\lambda\geq 1$ and $\beta\in[0,\vartheta)$ with $\vartheta$ being defined by \eqref{BZ26}. By Theorem \ref{Pre3:Sch} with $b=0$, we have
$$
\lambda^{\frac12(1-|\frac1{\bbp_0}|-\frac2{q_0}-\beta)}\|u_0\|_{\mL^\infty_T(\cC^{1+\beta})}+
\nor \nabla^2 u_0\nor_{\wt\mL^{q_0}_T(\wt\mL^{\bbp_0}_{\x_0})}\lesssim \nor f\nor_{\wt\mL^{q_0}_T(\wt\mL^{\bbp_0}_{\x_0})},
$$
and
$$
\|u_0(t_1)-u_0(t_0)\|_\infty\le C(\lambda)(t_1-t_0)^{\frac12}\nor f\nor_{\wt\mL^{q_0}_T(\wt\mL^{\bbp_0}_{\x_0})},
$$
and for each $i=1,\cdots,d$,
\begin{align*}
&\lambda^{\frac12(1-|\frac1{\bbp_i}|-\frac2{q_i}-\beta)}\|u_i\|_{\mL^\infty_T(\cC^{1+\beta})}+
\nor \nabla^2 u_i\nor_{\wt\mL^{q_i}_T(\wt\mL^{\bbp_i}_{\x_i})}
\lesssim \nor b_i\cdot\p_i u\nor_{\wt\mL^{q_i}_T(\wt\mL^{\bbp_i}_{\x_i})}
\lesssim \nor b_i\nor_{\wt\mL^{q_i}_T(\wt\mL^{\bbp_i}_{\x_i})}\|\p_i u\|_{\mL^\infty_T},
\end{align*}
and
$$
\|u_i(t_1)-u_i(t_0)\|_\infty\le C(\lambda)(t_1-t_0)^{\frac12}\nor b_i\nor_{\wt\mL^{q_i}_T(\wt\mL^{\bbp_i}_{\x_i})}\|\p_i u\|_{\mL^\infty_T}.
$$
Summing up the above inequalities for $i$ from $0$ to $d$, we obtain
$$
\lambda^{\frac12(\vartheta-\beta)}\|u\|_{\mL^\infty_T(\cC^{1+\beta})}+
\sum_{i=0}^d\nor \nabla^2 u_i\nor_{\wt\mL^{q_i}_T(\wt\mL^{\bbp_i}_{\x_i})}
\leq C_1\nor f\nor_{\wt\mL^{q_0}_T(\wt\mL^{\bbp_0}_{\x_0})}+C_2\kappa_1\|\nabla u\|_{\mL^\infty_T},
$$
where $C_1, C_2$ only depend on $T,\kappa_0,d,\bbp_i,q_i,\beta$, and
$$
\|u(t_1)-u(t_0)\|_\infty\le C(\lambda)(t_1-t_0)^{\frac12}\Big(\kappa_1\|\nabla u\|_{\mL^\infty_T}+\nor f\nor_{\wt\mL^{q_0}_T(\wt\mL^{\bbp_0}_{\x_0})}\Big).
$$
Choosing $C_0=(C_2/2)^{2/\vartheta}\vee 1$, we obtain \eqref{CC5} and \eqref{BB622} for all $\lambda\geq C_0\kappa_1^{2/\vartheta}$.
\end{proof}

\section{SDEs with mixed $L^\bbp$-drifts}
In this section we first establish a priori Krylov estimates for any solution of SDEs with mixed drifts and for any index $(q,\bbp)\in\sI_1$, where $\sI_1$ is defined in \eqref{Ind}. 
Using this a priori estimates, one can perform the classical Zvonkin transformation (see \cite{Xi-Zh}),
and then establish the weak well-posedness under conditions \eqref{SIG} and \eqref{BB0}. Moreover, we also obtain the two-sided density estimates.
As a byproduct, one improves the Krylov estimate to any index $(q,\bbp)\in\sI_2$, which is crucial for the strong well-posedness and the propagation of chaos.

Let $\xi_t$ be a given $\mR^d$-valued measurable adapted process. We consider the following SDE:
\begin{align}\label{SDE0}
\dif X_t=[\xi_t+b(t,X_t)]\dif t+\sigma(t,X_t)\dif W_t,
\end{align}
where $b:\mR_+\times\mR^d\to\mR^d$ and $\sigma: \mR_+\times\mR^d\to\mR^d\otimes\mR^d$
are Borel measurable functions.
We first introduce the following notion of solutions, also called weak solutions.
\bd\label{Def21}
Let ${\frak U}:=(\Omega,\sF,\mP, (\sF_t)_{t\geq 0})$ be a stochastic basis, and $\xi_t$ be a given
$\mR^d$-valued measurable $\sF_t$-adapted process with $\int^t_0|\xi_s|\dif s<\infty$ a.s. for each $t>0$, 
and $(X, W)$ be a pair of continuous $\sF_t$-adapted processes.
Let $\mu_0\in\cP(\mR^d)$.
We call $(X,W,{\frak U})$ a solution of SDE \eqref{SDE0} with initial distribution $\mu_0$ if 
\begin{enumerate}[(i)]
\item $\mu_0=\mP\circ X^{-1}_0$ and $W$ is a standard Brownian motion on ${\frak U}$.

\item For all $t\geq 0$,
$$
\int^t_0|b(s,X_s)|\dif s+\int^t_0|\sigma(s,X_s)|^2\dif s<\infty,\ a.s.
$$
and
$$
X_t=X_0+\int^t_0[\xi_s+b(s,X_s)]\dif s+\int^t_0\sigma(s,X_s)\dif W_s,\ a.s.
$$
\end{enumerate}
\ed

By Theorem \ref{Pre3:Sch}, we can establish the following a priori Krylov estimate (see \cite{Xi-Zh}).
\bl\label{Le27}
Suppose that \eqref{SIG} and \eqref{BB0} hold.
Then for any  $(q,\bbp)\in\sI_1$, $\bx\in\sX$ and $T, \delta>0$, there is a constant $C_{T,\delta}=C_{T,\delta}(\Theta)>0$ such that for all $f\in \widetilde\mL^q_T(\wt\mL^\bbp_\x)$
and any solution $X$ of SDE \eqref{SDE0} in the sense of Definition \ref{Def21},
\begin{align}\label{Kry10}
\mE\left(\int^T_0f(s,X_s)\dif s\right)\leq \nor f\nor_{\widetilde\mL^q_T(\wt\mL^\bbp_\x)}
\left[C_{T,\delta}+\delta\mE\left(\int^T_0|\xi_s|\dif s\right)\right].
\end{align}
\el
\begin{proof}
By Theorem \ref{Pre3:Sch} and as in \cite[Lemma 5.5]{Xi-Zh}, for any $(q,\bbp)\in\sI_1$, $\bx\in\sX$  and $T,\delta>0$, there 
is a constant $C_{T,\delta}>0$ such that for any stopping time $\tau\leq T$ and $f\in\wt\mL^q_T(\wt\mL^{\bbp}_\x)$,
\begin{align}\label{BB11}
\mE\left(\int^{\tau}_0f(s,X_s)\dif s\right)
\leq\nor f\nor_{\wt\mL^q_T(\wt\mL^{\bbp}_\x)}\left[C_{T,\delta}
+\delta\mE\left(\int^{\tau}_0(|\xi_s|+|b(s, X_s)|)\dif s\right)\right].
\end{align}
Now for $n\in\mN$, define a stopping time
$$
\tau_n:=\inf\left\{t>0: \int^t_0|b(s,X_s)|\dif s\geq n\right\}\wedge T.
$$
Since $(q_i,\bbp_i)\in\sI_1$, by applying \eqref{BB11} with $f(s,x)=b_i(s,x)$, we obtain
\begin{align*}
\mE\left(\int^{\tau_n}_0|b_i(s,X_s)|\dif s\right)
\leq\nor b_i\nor_{\wt\mL^{q_i}_T(\wt\mL^{\bbp_i}_{\x_i})}\left[C_{T,\delta}+
\delta\mE\left(\int^{\tau_n}_0(|\xi_s|+|b(s, X_s)|)\dif s\right)\right].
\end{align*}
Summing up the above inequalities for $i$ from $1$ to $d$, we get
$$
\mE\left(\int^{\tau_n}_0|b(s,X_s)|\dif s\right)
\leq\kappa_1\left[C_{T,\delta}+\delta\mE\left(\int^{\tau_n}_0(|\xi_s|+|b(s, X_s)|)\dif s\right)\right].
$$
Letting $\delta=1/(2\kappa_1)$ and $n\to\infty$, we obtain
$$
\mE\left(\int^T_0|b(s,X_s)|\dif s\right)\leq \kappa_1C_{T,1/(2\kappa_1)}+\frac1{2}\mE\left(\int^T_0(|\xi_s|+|b(s, X_s)|)\dif s\right),
$$
which implies 
$$
\mE\left(\int^T_0|b(s,X_s)|\dif s\right)\leq 2\kappa_1C_{T,1/(2\kappa_1)}+\mE\left(\int^T_0|\xi_s|\dif s\right).
$$
Substituting this into \eqref{BB11} with $\tau=T$, we complete the proof.
\end{proof}

In the above lemma, the requirement of $(q,\bbp)\in\sI_1$ is too strong for applications. We need to improve it to $(q,\bbp)\in\sI_2$.
Firstly, we use Theorem \ref{Pre3:Sch0} and the above a priori Krylov estimate to construct the Zvonkin transformation.
For each $i=1,\cdots,d$, consider the following backward PDE:
\begin{align}\label{Pre3:PDE}
\p_t u_i+\tfrac12\tr((\sigma\sigma^*)\cdot\nabla^2u_i)+b\cdot\nabla u_i-\lambda u_i+b_i=0,\quad u_i(T)=0.
\end{align}
By reversing the time variable and by Theorem \ref{Pre3:Sch0}, there is a unique solution $u_i$ satisfying the following estimates: for any $\beta\in(0,\vartheta)$, 
where $\vartheta$ is defined in \eqref{BZ26},
there are $C_0, C_1\geq 1$ such that for all $\lambda\geq C_0\kappa_1^{2/\vartheta}$,
\begin{align}\label{CC0}
\lambda^{\frac12(\vartheta-\beta)}\|u_i\|_{\mL^\infty_T(\cC^{1+\beta})}+
\sum_{j=0}^d\nor\nabla^2 u_{ij}\nor_{\widetilde{\mL}^{q_{ij}}_T(\wt\mL^{\bbp_{ij}}_{\x_{ij}})}\le C_1\kappa_1,
\end{align}
and for all $0\leq t_0<t_1\leq T$,
\begin{align}\label{AC5}
\|u_i(t_1)-u_i(t_0)\|_\infty\leq C(\lambda)|t_1-t_0|^{1/2},
\end{align}
where 
\begin{align}\label{BZ11}
u_i=u_{i0}+u_{i1}+u_{i2}+\cdots+u_{id},
\end{align}
and
$$
q_{i0}=q_i,\ \bbp_{i0}=\bbp_i,\ \x_{i0}=\x_i,\ q_{ij}=q_j,\ \bbp_{ij}=\bbp_j,\ \x_{ij}=\x_j,\ j=1,\cdots,d.
$$
Below we set
$$
\u=(u_1,\cdots,u_d).
$$
By \eqref{CC0}, for any $\beta\in(0,\vartheta)$, we can choose $\lambda$ large enough so that
$$
\|\nabla \u\|_{\mL^\infty_T}\leq \|\u\|_{\mL^\infty_T(\cC^{1+\beta})}\leq\tfrac12.
$$
Once $\lambda$ is chosen, it shall be fixed below without further notice.
Now if we define
$$
\Phi(t,x):=x+\u(t,x),
$$
then for each $t$,
$$
x\mapsto\Phi(t,x)\mbox{ is a $C^1$-diffeomorphism},
$$
and
\begin{align}\label{CC3}
\|\nabla\Phi\|_{\mL^\infty_T}+\|\nabla\Phi^{-1}\|_{\mL^\infty_T}\leq 2,
\end{align}
and by \eqref{AC5}, for all $0\leq t_0<t_1\leq T$,
\begin{align}\label{AC6}
\|\Phi(t_1)-\Phi(t_0)\|_\infty\leq C(\lambda)(t_1-t_0)^{1/2}.
\end{align}
We have the following result (see \cite[Theorem 3.10]{Xi-Zh}).
\bl[Zvonkin's transformation]\label{Le28}
Under \eqref{SIG} and \eqref{BB0}, $Y_t:=\Phi(t,X_t)$ solves the following SDE
\begin{align}\label{SDE34}
Y_t=Y_0+\int^t_0\xi_s\cdot\nabla\Phi(s,\Phi^{-1}(s,Y_s))\dif s+\int^t_0\tilde b(s, Y_s)\dif s+\int^t_0\tilde\sigma(s,Y_s)\dif W_s,
\end{align}
where $Y_0:=\Phi(0,X_0)$ and
$$
\tilde b(s,y):=\lambda \u(s,\Phi^{-1}(s,y)),\ \ \tilde\sigma(s,y):=(\sigma^*\nabla \Phi)(s,\Phi^{-1}(s,y)).
$$
Moreover, for any $\beta\in(0,\vartheta\wedge\gamma_0)$, where $\vartheta$ is defined by \eqref{BZ26},
\begin{align}\label{BZ13}
\tilde b,\ \nabla\tilde b\in\mL^\infty_T,\ \tilde\sigma\in\mL^\infty_T(\cC^{\beta}),
\end{align}
and for some $\tilde\kappa_0\geq 1$,
\begin{align}\label{BZ104}
\tilde\kappa^{-1}_0|\eta|^2\leq |\tilde\sigma(s,y)\eta|^2\leq \tilde\kappa_0|\eta|^2,\ \ \eta\in\mR^d.
\end{align}
Vice versa, if $Y_t$ solves SDE \eqref{SDE34}, then $X_t:=\Phi^{-1}(t,Y_t)$ solves SDE \eqref{SDE0}.
\el
\begin{proof}
For each $n\in\mN$, define 
$$
\u^n(t,x):=(\u(t,\cdot)*\varGamma_{1/n})(x),\ \Phi^n(t,x):=x+\u^n(t,x).
$$
By It\^o's formula, we have
\begin{align*}
\Phi^n(t,X_t)=\Phi^n(0,X_0)+\int^t_0[\xi_s\cdot\nabla\Phi^n+\sL\Phi^n](s,X_s)\dif s+\int^t_0(\sigma^*\nabla\Phi^n)(s,X_s)\dif W_s,
\end{align*}
where
$$
\sL:=\p_s+\tr(a\cdot\nabla^2)+b\cdot\nabla,\ \ a:=(\sigma\sigma^*)/2.
$$
Since $x\mapsto \u(t,x)$ is $\cC^{1+\beta}$-continuous, it is easy to see that for each $t,x$,
$$
\lim_{n\to\infty}\nabla^j\Phi^n(t,x)=\nabla^j\Phi(t,x),\ \ j=0,1.
$$
Therefore, to show \eqref{SDE34}, it suffices to show that as $n\to\infty$,
$$
\int^t_0|\sL\Phi^n-\lambda \u|(s,X_s)\dif s\to 0,\ \ a.s.
$$
For $m\in\mN$, we define the stopping time
$$
\tau_m:=\inf\left\{t>0: |X_t|+\int^t_0|\xi_s|\dif s\geq m\right\}.
$$
Since $\tau_m\to\infty$ as $m\to\infty$, it suffices to show that for each fixed $m\in\mN$,
\begin{align} \label{CC01}
\mE\left(\int^{t\wedge\tau_m}_0|\sL\Phi^n-\lambda \u|(s,X_s)\dif s\right)=0.
\end{align}
Note that by definition,
\begin{align}\label{BZ12}
\begin{split}
\sL\Phi^n-\lambda \u&=\p_s\u^n+\tr(a\cdot\nabla^2\u^n)+b\cdot\nabla \u^n+b-\lambda \u\\
&=\big[\tr(a\cdot(\nabla^2\u)*\varGamma_{1/n})-\tr(a\cdot\nabla^2\u)*\varGamma_{1/n}\big]\\
&\quad+\big[b\cdot\nabla (\u*\varGamma_{1/n})-(b\cdot\nabla\u)*\varGamma_{1/n}\big]\\
&\quad+[b*\varGamma_{1/n}-b]+[\lambda(\u*\varGamma_{1/n}-\u)].
\end{split}
\end{align}
For each $i,j$, since $(q_{ij},\bbp_{ij})\in\sI_1$, by the Krylov estimates \eqref{Kry10} and \eqref{CC0}, we have
\begin{align*}
\mE\left(\int^{t\wedge\tau_m}_0|\nabla^2u^n_{ij}-\nabla^2u_{ij}|(s,X_s)\dif s\right)
&\leq C_m\|\nabla^2(u^n_{ij}-u_{ij})\b1_{B_m}\|_{\mL^{q_{ij}}_T(\mL^{\bbp_{ij}}_{\x_{ij}})}\to 0.
\end{align*}
From this and by \eqref{BZ11} and \eqref{BZ12}, it is easy to see that \eqref{CC01} holds. Moreover,
\eqref{BZ13} and \eqref{BZ104} directly follow by their definitions and \eqref{CC0}. 
On the other hand, if $Y_t$ solves SDE \eqref{SDE34}, then by similar calculations,  $X_t:=\Phi^{-1}(t,Y_t)$ solves SDE \eqref{SDE0}.
We omit the details here.
\end{proof}

\br\rm
Consider SDE \eqref{SDE0} with $\xi\equiv0$ and assume \eqref{SIG} and \eqref{BB0}.
An immediate consequence of Zvonkin's transformation together with \eqref{CC3} and \eqref{AC6} is that for any $p\geq1$ and $T>0$,
there is a constant $C=C(p,T,\Theta)>0$ such that
\begin{align}\label{Mom}
\mE|X_t-X_s|^{2p}\leq C|t-s|^p,\ \ t,s\in[0,T].
\end{align}
\er

 Now we show the following main result of this section.
\bt\label{Main1}
Suppose that \eqref{SIG} and \eqref{BB0} hold. For any $\mu_0\in\cP(\mR^d)$, there is a unique weak solution to SDE \eqref{SDE0} with $\xi\equiv0$ and initial distribution $\mu_0$
in the sense of Definition \ref{Def21}.
Moreover, we have :
\begin{enumerate}[(i)]
\item For each $t>0$, $X_t$ admits a density $\rho^X_t(y)$ with the following two-sided estimate:
for any $T>0$, there are $\delta_1, C_1\geq 1$ such that for all $t\in(0,T]$ and $y\in\mR^d$,
\begin{align}\label{DD1}
\frac{C^{-1}_1}{t^{d/2}}\int_{\mR^d}\e^{-\frac{\delta_1 |x-y|^2}{t}}\mu_0(\dif x)\leq
\rho^X_t(y)\leq \frac{C_1}{t^{d/2}}\int_{\mR^d}\e^{-\frac{|x-y|^2}{\delta_1t}}\mu_0(\dif x).
\end{align}
\item Let  $\vartheta$ be defined as in \eqref{BZ26}.
For any $\beta\in(0,\vartheta\wedge \gamma_0)$ and $T>0$, there are $\delta_2,C_2\geq 1$ such that for all $t\in(0,T]$ and $y,y'\in\mR^d$,
\begin{align}\label{Ho}
\frac{|\rho^X_t(y)-\rho^X_t(y')|}{|y-y'|^\beta}\leq C_2 t^{-\frac{d+\beta}2}
\left[\int_{\mR^d}\e^{-\frac{|x-y|^2}{\delta_2 t}}\mu_0(\dif x)+\int_{\mR^d}\e^{-\frac{|x-y'|^2}{\delta_2 t}}\mu_0(\dif x)\right].
\end{align}
\item For any $(q,\bbp)\in\sI_2$ and $T>0$, there is a constant $C_0>0$ such that for any $f\in\wt\mL^q_T(\wt\mL^\bbp_\x)$,
\begin{align}\label{Kry11}
\mE\left(\int^T_0f(s,X_s)\dif s\right)\leq C_0\nor f\nor_{\wt\mL^q_T(\wt\mL^\bbp_\x)}.
\end{align}
\end{enumerate}
\et
\begin{proof}
By \eqref{BZ13} and \eqref{BZ104}, it is well known that SDE \eqref{SDE34} with $\xi\equiv0$ admits a unique weak solution (cf. \cite{SV}).
The existence and uniqueness of weak solutions for the original SDE follow from Lemma \ref{Le28}.
 Next we shall prove \eqref{DD1}, \eqref{Ho} and \eqref{Kry11}. 
 
 (i) Let $\wt\sL$ be the generator of SDE \eqref{SDE34}, i.e.,
$$
\wt\sL:=\tr((\wt\sigma\wt\sigma^*)\cdot\nabla^2)/2+\wt b\cdot\nabla.
$$
By \eqref{BZ13}, \eqref{BZ104} and Theorems 1.1, 1.3 and 2.3  of \cite{CHXZ17},  there is a fundamental solution $p(s,x,t,y)$ associated with $\wt\sL$, which satisfies the following estimates: for all $0\leq s<t\leq T$ and $x,y\in\mR^d$,
$$
\frac{C^{-1}_0}{(t-s)^{d/2}}\e^{-\frac{\delta_0|x-y|^2}{t-s}}\leq p(s,x,t,y)\leq \frac{C_0}{(t-s)^{d/2}}\e^{-\frac{|x-y|^2}{\delta_0(t-s)}},
$$
and for any $\beta\in(0,\vartheta\wedge\gamma_0)$, and  for all $0\leq s<t\leq T$ and $x,y,y'\in\mR^d$,
$$
|p(s,x,t,y)-p(s,x,t,y')|\leq C_0|y-y'|^\beta(t-s)^{-\frac{d+\beta}2}\left[\e^{-\frac{|x-y|^2}{\delta_0(t-s)}}+\e^{-\frac{|x-y'|^2}{\delta_0(t-s)}}\right],
$$
where $\delta_0, C_0\geq 1$ only depend on $\Theta$ and  the bounds of $\tilde b$ and $\tilde\sigma$.
In particular, $p(0,x,t,y)$ is just the density of the solution of SDE \eqref{SDE0} starting from $x$ at time zero.
Note that the density $\rho^Y_t(y)$ of $Y_t$ starting from the initial distribution $\wt\mu_0=\mu_0\circ\Phi(0,\cdot)^{-1}$ is given by
$$
\rho^Y_t(y)=\int_{\mR^d}p(0,x,t,y)\wt\mu_0(\dif x).
$$
This can be shown by considering a smooth approximation and taking weak limits (see \cite[Section 5.1]{MZ21} for more details).
We thus have that for any $t\in(0,T]$ and all $y,y'\in\mR^d$,
$$
\frac{C^{-1}_0}{t^{d/2}}\int_{\mR^d}\e^{- \frac{\delta_0|x-y|^2}{t}}\wt\mu_0(\dif x)\leq
\rho^Y_t(y)\leq \frac{C_0}{t^{d/2}}\int_{\mR^d}\e^{-\frac{|x-y|^2}{\delta_0 t}}\wt\mu_0(\dif x)
$$
and
\begin{align}\label{AMG1}
|\rho^Y_t(y)-\rho^Y_t(y')|\leq \frac{C_0|y-y'|^\beta}{t^{(d+\beta)/2}}
\int_{\mR^d}\left[\e^{-\frac{|x-y|^2}{\delta_0t}}+\e^{-\frac{|x-y'|^2}{\delta_0t}}\right]\wt\mu_0(\dif x).
\end{align}
On the other hand, by change of variables, we have
\begin{align}\label{AM2}
\rho^X_t(y)=\rho^Y_t(\Phi(t,y))\det(\nabla\Phi(t,y)),
\end{align}
and for some $\tilde C_0\geq 1$,
$$
\frac{\tilde C^{-1}_0}{t^{d/2}}\int_{\mR^d}\e^{-\frac{\delta_0|\Phi(0,x)-\Phi(t,y)|^2}{ t}}\mu_0(\dif x)\leq\rho^X_t(y)\leq 
\frac{\tilde C^{-1}_0}{t^{d/2}}\int_{\mR^d}\e^{-\frac{|\Phi(0,x)-\Phi(t,y)|^2}{\delta_0t}}\mu_0(\dif x).
$$
which together with the following two estimates yields \eqref{DD1},
$$
|\Phi(0,x)-\Phi(t,y)|^2\geq\tfrac12|\Phi(t,x)-\Phi(t,y)|^2-|\Phi(0,x)-\Phi(t,x)|^2\stackrel{\eqref{CC3}\eqref{AC6}}{\geq} \tfrac1{8}|x-y|^2-Ct,
$$
and
$$
|\Phi(0,x)-\Phi(t,y)|^2\leq 2|\Phi(t,x)-\Phi(t,y)|^2+2|\Phi(0,x)-\Phi(t,x)|^2\stackrel{\eqref{CC3}\eqref{AC6}}{\leq} 8|x-y|^2-Ct.
$$

(ii) By \eqref{AM2} and \eqref{AMG1}, we have
\begin{align*}
|\rho^X_t(y)-\rho^X_t(y')|&\leq|\rho^Y_t(\Phi(t,y))-\rho^Y_t(\Phi(t,y'))|\det(\nabla\Phi(t,y))\\
&+\rho^Y_t(\Phi(t,y'))|\det(\nabla\Phi(t,y))-\det(\nabla\Phi(t,y'))|\\
&\lesssim\frac{|y-y'|^\beta}{t^{(d+\beta)/2}}
\int_{\mR^d}\left[\e^{-\frac{|x-\Phi(t,y)|^2}{\delta_0t}}+\e^{-\frac{|x-\Phi(t,y')|^2}{\delta_0t}}\right]\wt\mu_0(\dif x)\\
&+\frac{1}{t^{d/2}}\int_{\mR^d}\e^{-\frac{|x-\Phi(t,y)|^2}{\delta_0 t}}\wt\mu_0(\dif x)|\nabla\Phi(t,y)-\nabla\Phi(t,y')|,
\end{align*}
which in turn implies \eqref{Ho} by \eqref{CC0}.

(iii) For nonnegative $f\in \widetilde\mL^q_T(\wt\mL^\bbp_\x)$ with $(q,\bbp)\in\sI_2$, by \eqref{DD1} and \eqref{226}, we get
\begin{align*}
&\mE\left(\int^T_0f(s,X_s)\dif s\right)=\int^T_0\int_{\mR^d}f(s,y)\rho^X_s(y)\dif y\dif s\\
&\qquad\leq\int^T_0\int_{\mR^d}f(s,y)\left(\frac{C_2}{s^{d/2}}\int_{\mR^d}\e^{-\frac{|x-y|^2}{2\delta_1s}}\mu_0(\dif x)\right)\dif y\dif s\\
&\qquad=C_2(2\pi\delta_1)^{d/2}\int_{\mR^d}\left(\int^T_0\mE f(s,x-W_{\delta_1s})\dif s\right)\mu_0(\dif x)
\leq C_3\nor f\nor_{\widetilde\mL^q_T(\wt\mL^\bbp_\x)}.
\end{align*}
The proof is complete.
\end{proof}

As a corollary, we have the following important exponential integrability of singular functionals.
\bc(Khasminskii's estimate)
Let $X$ be the unique solution of SDE \eqref{SDE0} in Theorem \ref{Main1}. 
For any $T,\lambda>0$, $(q,\bbp)\in\sI_2$ and $\beta\in(0,2-|\frac1\bbp|-\frac2q)$, there is a constant 
$C_1>0$ 
depending only on $T,\lambda,d,\beta,\kappa_0,\kappa_1,q_i,\bbp_i, q,\bbp$ such that for all $f\in\widetilde\mL^q_T(\wt\mL^\bbp_\x)$,
\begin{align}\label{Kha1}
\mE\exp\left\{\lambda\int^T_{0}f(s,X_s)\dif s\right\}\leq \e^{C_1 \nor f\nor^{2/\beta}_{\widetilde\mL^q_T(\wt\mL^\bbp_\x)}}.
\end{align}
Moreover, if $b$ is bounded measurable, then for some $C_2=C_2(T,\lambda,d,\beta,\kappa_0,q,\bbp)>0$,
\begin{align}\label{Kha11}
\mE\exp\left\{\lambda\int^T_{0}f(s,X_s)\dif s\right\}\leq \e^{C_2\(\|b\|^2_{\mL^\infty_T}+ \nor f\nor^{2/\beta}_{\widetilde\mL^q_T(\wt\mL^\bbp_\x)}\)}.
\end{align}
\ec
\begin{proof}
Let $\beta\in(0,2-|\frac1\bbp|-\frac2q)$. For \eqref{Kha1}, by \cite[Lemma 3.5]{Xi-Zh}, it suffices to show that for any $0\leq t_0<t_1\leq T$,
\begin{align}\label{AA1}
\mE\left(\int^{t_1}_{t_0}f(s,X_s)\dif s\Big|\sF_{t_0}\right)\leq C_0(t_1-t_0)^{\frac\beta2}\nor f\nor_{\widetilde\mL^q_T(\wt\mL^\bbp_\x)}.
\end{align}
Let $\frac1{q'}=\frac1q+\frac\beta2$. Since $\beta\in(0,2-|\frac1\bbp|-\frac2q)$, we have $(q',\bbp)\in\sI_2$. By  \eqref{Kry11} and H\"older's inequality,
$$
\mE\left(\int^{t_1-t_0}_0f(s,X_s)\dif s\right)\leq C_0\nor f\nor_{\widetilde\mL^{q'}_{t_1-t_0}(\wt\mL^\bbp_\x)}\leq C_0(t_1-t_0)^{\frac\beta2}\nor f\nor_{\widetilde\mL^{q}_T(\wt\mL^\bbp_\x)}.
$$
By the Markov property of $X_t$, we get \eqref{AA1}. \eqref{Kha11} follows by Girsanov's theorem.
\end{proof}

\bt\label{Th217}
{\bf (Strong well-posedness)} In addition to the assumptions, in Theorem \ref{Main1}, we also assume \eqref{SIG0} and that $(q_i,\bbp_i)\in\sI^o$ in \eqref{BB0}. Then 
there is a unique strong solution to SDE \eqref{SDE0} with $\xi\equiv0$.
\et
\begin{proof}
By Yamada-Watanabe's theorem, it suffices to show the pathwise uniqueness.
 But this follows by Zvonkin's transformation (see Lemma \ref{Le28}),  Lemma \ref{Le21} and \eqref{Kha1} (see \cite[Theorem 3.9]{Xi-Zh} for more details).
\end{proof}

\section{Well-posedness of dDDSDEs with mixed $L^\bbp$-drifts}
We consider the following {distribution density-distribution dependent SDE} (abbreviated as {dDDSDE}):
\begin{align}\label{SDE2}
\dif X_t=b(t,X_t,\rho_{t}(X_t),\mu_{X_t})\dif t+\sigma(t,X_t)\dif W_t,
\end{align}
where $\rho_{t}(x)$ is the density of $X_t$ and $b(t,x,r,\mu):\mR_+\times\mR^d\times\mR_+\times\cP(\mR^d)\to\mR^d$ is a measurable function.
As in Definition \ref{Def21}, we introduce the following notion of solutions to the above SDE.
\bd\label{Def26}
Let ${\frak U}:=(\Omega,\sF,\mP, (\sF_t)_{t\geq 0})$ be a stochastic basis and $(X,W)$ be a pair of continuous $\sF_t$-adapted processes.
Let $\mu_0\in\cP(\mR^d)$.
We call $(X,W,{\frak U})$ a solution of dDDSDE \eqref{SDE2} with initial distribution $\mu_0$ if 

(i) $\mu_0=\mP\circ X^{-1}_0$ and $W$ is a standard Brownian motion on ${\frak U}$.

(ii) For each $t>0$, the distribution $\mu_{X_t}$ of $X_t$ admits a density $\rho_{t}$,
$$
\int^t_0|b(s,X_s,\rho_{s}(X_s),\mu_{X_s})|\dif s+\int^t_0|\sigma(s,X_s)|^2\dif s<\infty,\ a.s.,
$$
and
$$
X_t=X_0+\int^t_0b(s,X_s,\rho_{s}(X_s),\mu_{X_s})\dif s+\int^t_0\sigma(s,X_s)\dif W_s,\ a.s.
$$
\ed
Let $T>0$ and $\mC_T$ be the space of all continuous functions from $[0,T]$ to $\mR^d$. We use 
$\omega$ to
denote a path in $\mC_T$ 
and by $w_t(\omega)=\omega_t$ to denote the coordinate process. Let $\cB_t:=\sigma\{w_s, s\leq t\}$ be the natural filtration.
We also introduce the following notion of martingale solutions to dDDSDE \eqref{SDE2}.
\bd
Let $\mu_0\in\cP(\mR^d)$. A probability measure $\mP\in\cP(\mC_T)$ is called a martingale solution of dDDSDE \eqref{SDE2} with initial distribution $\mu_0$
if $\mP\circ w_0^{-1}=\mu_0$ and for any $f\in C^2_c(\mR^d)$, the process
\begin{align}\label{BZ24}
M^f_t(\omega):=f(w_t)-f(w_0)-\int^t_0\Big(\tfrac12\tr((\sigma\sigma^*)(s,w_s)\cdot\nabla^2)+b(s,w_s,\rho_s(w_s),\mu_s)\cdot\nabla \Big)f(w_s)\dif s
\end{align}
is a $\cB_t$-martingale, where $\mu_s:=\mP\circ w^{-1}_s$ has a density $\rho_s(x)$. 
We shall use $\sM^{\sigma,b}_{\mu_0}$ to denote the set of all martingale 
solutions of dDDSDE \eqref{SDE2} associated with $\sigma,b$ and initial distribution ${\mu_0}$.
\ed

\br\rm\label{Re43}
It is well known that weak solutions are equivalent to the martingale solutions (see \cite{SV}). If we consider the {classical SDE}, i.e., 
$b$ only depends on $(t,x)$, and if for each starting point $(s,x)$, there is a unique martingale solution starting from $(s,x)$, then as usual, we say the martingale problem is well-posed.
\er

\subsection{Relative entropy}
In this subsection we recall the notion and some basic facts about the relative entropy.
Let $E$ be a Polish space and $\mu,\nu$ be two probability measures on $E$. The relative entropy $\cH(\mu|\nu)$ is defined by
\begin{align}\label{Rela}
\cH(\mu|\nu):=
\left\{
\begin{aligned}
&\int_E \frac{\dif\mu}{\dif\nu}\log\Big(\frac{\dif\mu}{\dif\nu}\Big)\dif\nu,&\mu\ll\nu,\\
&\infty,&\mbox{otherwise.}
\end{aligned}
\right.
\end{align}
Since $x\mapsto x\log x$ is convex on $[0,\infty)$, by Jensen's inequality, we have $\cH(\mu|\nu)\geq 0$.

The following theorem contains the main tools used below (see \cite[Theorem 2.1(ii)]{BV05}, \cite[Lemma 1.4.3(a)]{DE97} and \cite[Lemma 3.9]{DM01}).
\bt
\begin{enumerate}[(i)]
\item (Pinsker's inequality) For any two probability measures $\mu,\nu$,
\begin{align}\label{Pin0}
\|\mu-\nu\|^2_{\rm var}\leq 2\cH(\mu|\nu).
\end{align}
\item (The weighted Pinsker inequality) For any  $\mu,\nu\in\cP(E)$ and Borel measurable function $f$, 
\begin{align}\label{Pin1}
|\<\mu-\nu, f\>|^2\leq 2\left(1+\log\int_E \e^{f(x)^2}\nu(\dif x)\right)\cH(\mu|\nu).
\end{align}
\item (Variational representation of the relative entropy) For any $\mu,\nu\in\cP(E)$,
\begin{align}\label{Var}
\cH(\mu|\nu)=\sup_{\psi\in\cB_b(E)}\left(\int_E\psi\dif\mu-\log\int_E\e^\psi\dif\nu\right),
\end{align}
where $\cB_b(E)$ is the set of all bounded Borel measurable functions.
\item (Dimensional bounds on entropy) Let $\mu^N$ be a symmetric probability measure on $E^N$ and $\mu\in\cP(E)$. Then for any $k\leq N$, 
\begin{align}\label{BB4}
\cH\(\mu^{N,k}|\mu^{\otimes k}\)\leq\frac{2k}{N}\cH\(\mu^N|\mu^{\otimes N}\),
\end{align}
where $\mu^{N,k}$ is the marginal distribution of the first $k$-component of $\mu^N$.
\end{enumerate}
\et

We recall the following entropy formula for the martingale solutions of classical SDEs,
which is a consequence of Girsanov's theorem (see \cite[Lemma 4.4 and Remark 4.5]{La21} for the most general form).
\bl\label{Lem34}
For $i=1,2$, let $b^i:\mR_+\times\mR^d\to\mR^d$ be two measurable functions. Suppose that the martingale problem associated with 
$(\sigma,b^2)$ is well-posed (see Remark \ref{Re43}). Let $\mu^1_0,\mu^2_0\in\cP(\mR^d)$ be two initial distributions.
For any two martingale solutions $\mP_i\in\cM^{\sigma,b^i}_{\mu^i_0}$, $i=1,2$, and 
any $t\geq 0$, if we let $\mu^i_t:=\mP_i\circ w_t^{-1}$ be the marginal distributions, then
$$
\cH(\mu^1_t|\mu^2_t)\leq\cH(\mu^1_0|\mu^2_0)+\frac12\mE^{\mP_2}\left(\int^t_0|\sigma^{-1}(s,w_s)(b^1(s,w_s)-b^2(s,w_s))|^2\dif s\right).
$$
\el

\subsection{Stability of density}
In this section we prepare a stability result about the density of classical SDEs.
Our starting point is the associated Fokker-Planck equation.
Fix $z\in\mR^d$. Let 
$$
A^z_{s,t}:=\int^t_s A(r,z)\dif r\mbox{ with $A(r,z)=(a_{ij}(r,z))=((\sigma\sigma^*)_{ij}(r,z))/2$}.
$$
Let $P_{s,t}^z f(x)$ be the Gaussian heat kernel associated with $A^z_{s,t}$, i.e.,
$$
P_{s,t}^z f(x)=\int_{\mR^d} h_{A^z_{s,t}}(x-y)f(y)\dif y,
$$
where for a symmetric positive definite matrix $A$,
$$
h_A(x):=\frac{\e^{-\<A^{-1}x,x\>/2}}{\sqrt{(2\pi)^d\det(A)}}.
$$
\bl\label{Le46}
Let $\beta\in[0,1]$, $k\in\mN_0$, $\bbp\in[1,\infty]^d$ and $\x\in\sX$. Under \eqref{SIG},
for any $T>0$, there is a constant $C=C(T, d,\beta,k,\bbp,\kappa_0)>0$ such that for all $0\leq s<t\leq T$ and $0\leq f\in\wt\mL^\bbp_\x$,
\begin{align*}
|\nabla^kP_{s,t}^z(|\cdot|^\beta f)(0)|\leq C (t-s)^{\frac12(\beta-k-|\frac1\bbp|)} \nor f\nor_{\wt\mL^\bbp_\x}.
\end{align*}
\el
\begin{proof}
First of all, by definition and  \eqref{SIG}, it is easy to see that for some $\lambda>0$,
$$
|\nabla^kh_{A^z_{s,t}}(x)|\lesssim (t-s)^{-\frac{k+d}2}\e^{-\frac{|x|^2}{\lambda(t-s)}}=(t-s)^{-\frac k2}(2\pi\lambda)^{\frac{d}2}g_{\lambda(t-s)}(x),
$$
and for some $\lambda'>\lambda$,
\begin{align*}
|\nabla^kP_{s,t}^z(|\cdot|^\beta f)(0)|
&\lesssim(t-s)^{-\frac k2}\int_{\mR^d}g_{\lambda (t-s)}(y)|y|^\beta f(y)\dif y\no\\
&\lesssim(t-s)^{\frac{\beta}{2}-\frac k2}\int_{\mR^d}g_{\lambda' (t-s)}(y) f(y)\dif y.
\end{align*}
Let $\bbp'\in(1,\infty)^d$ be defined by $\frac1{\bbp}+\frac1{\bbp'}={\bf 1}$. Fix $r>0$.
By H\"older's inequality we have
\begin{align}
\int_{\mR^d}g_{\lambda' (t-s)}(y) f(y)\dif y
&=\frac{1}{|B^r_0|}\int_{\mR^d}\int_{\mR^d}g_{\lambda'(t-s)}(y)\b1_{B^r_z}(y)f(y)\dif y\dif z\no\\
&\leq\frac{1}{|B^r_0|}\int_{\mR^d}\|\b1_{B^r_z}g_{\lambda' (t-s)}\|_{\mL^{\bbp'}_\x}\|\b1_{B^r_z}f\|_{\mL^\bbp_\x}\dif z\no\\ 
&\leq\frac{1}{|B^r_0|}\left(\int_{\mR^d}\|\b1_{B^r_z}g_{\lambda' (t-s)}\|_{\mL^{\bbp'}_\x}\dif z\right)\nor f\nor_{\wt\mL^\bbp_\x}\label{BZ401}.
\end{align}
Below, without loss of generality, we suppose $s=0$. By a change of variables, we have
\begin{align*}
\int_{\mR^d}\|\b1_{B^r_{z}}g_{\lambda t}\|_{\mL^{\bbp'}_\x}\dif z
&= (2\pi\lambda t)^{-\frac{d}2}\int_{\mR^d}\|\b1_{B^r_{z}}\e^{-\frac{|\cdot|^2}{\lambda t}}\|_{\mL^{\bbp'}_\x}\dif z\\
&\lesssim t^{-\frac{d}2}\prod_{i=1}^d\int_{\mR}\left(\int_{|y_i-z_i|\leq r}\e^{-\frac{p'_i|y_i|^2}{\lambda t}}\dif y_i\right)^{\frac{1}{p'_i}}\dif z_i
=:t^{-\frac{d}2}\prod_{i=1}^d \cJ_i.
\end{align*}
For each $i$, we have
\begin{align*}
\cJ_i&=\int_{|z_i|\leq 2r}\left(\int_{|y_i-z_i|\leq r}\e^{-\frac{p'_i|y_i|^2}{\lambda t}}\dif y_i\right)^{\frac{1}{p'_i}}\dif z_i
+\int_{|z_i|>2r}\left(\int_{|y_i-z_i|\leq r}\e^{-\frac{p'_i|y_i|^2}{\lambda t}}\dif y_i\right)^{\frac{1}{p'_i}}\dif z_i\\
&\leq \int_{|z_i|\leq 2r}\left(\int_\mR\e^{-\frac{p'_i|y_i|^2}{\lambda t}}\dif y_i\right)^{\frac{1}{p'_i}}\dif z_i+\int_{|z_i|>2r}\e^{-\frac{p_i'(|z_i|-r)^2}{\lambda t}}\left(\int_{|y_i-z_i|\leq r}\dif y_i\right)^{\frac{1}{p'_i}}\dif z_i\\
&\lesssim \left(\int_\mR\e^{-\frac{p'_i|y_i|^2}{\lambda t}}\dif y_i\right)^{\frac{1}{p'_i}}+\int_{\mR}\e^{-\frac{p'_i|z_i|^2}{2\lambda t}}\dif z_i
\lesssim t^{\frac1{2p'_i}}+t^{\frac12}\lesssim t^{\frac1{2p'_i}}=t^{\frac12(1-\frac1{p_i})}.
\end{align*}
Hence,
$$
\int_{\mR^d}\|\b1_{B^r_{z}}g_{\lambda t}\|_{\mL^{\bbp'}_\x}\dif z\lesssim t^{-\frac{d}2}\prod_{i=1}^dt^{\frac12(1-\frac1{p_i})}= t^{-|\frac1{\bbp}|/2}.
$$
Combining the above estimates, we obtain the desired estimate.
\end{proof}

The following stability result shall be used below to show the existence and uniqueness.
\bl\label{Le47}
Let $b_0, b_1$ be two Borel measurable functions satisfying \eqref{BB0} and for $k=0,1$, $\mu_k(\dif x):=\rho^o_k(x)\dif x$ with 
$\rho^o_k\in\mL^\infty$.
Let $\mP_k\in\cM^{\sigma,b_k}_{\mu_k}$ be the unique martingale solution and $\rho_k(t,x)$ be the density of the coordinated process
$w_t$ under $\mP_k$. Then for any $T>0$, there is a constant $C=C(T,\Theta)>0$ such that for all $t\in[0,T]$,
\begin{align}\label{CS2}
\|\rho_0(t)-\rho_1(t)\|_{\mL^\infty}\lesssim_C \|\rho^o_0-\rho^o_1\|_{\mL^\infty}+
\sum_{i=1}^d\int^t_0 (t-s)^{-\frac12(1+|\frac1{\bbp_i}|)}\nor b^i_0(s)-b^i_1(s)\nor_{\wt\mL^{\bbp_i}_{\x_i}}\dif s.
\end{align}
\el
\begin{proof}
First of all, by the heat kernel estimate \eqref{DD1}, we have for all $t,y$,
\begin{align}\label{CS1}
\rho_k(t,y)\leq \frac{C_1}{t^{d/2}}\int_{\mR^d}\e^{-\frac{|x-y|^2}{2\delta_1t}}\rho^o_k(x)\dif x\lesssim\|\rho^o_k\|_{\mL^\infty},\ k=0,1.
\end{align}
Note that $\rho_k$ solves the following Fokker-Planck equation in the distributional sense:
$$
\p_t\rho_k=\p_i\p_j(a_{ij}\rho_k)+\div(b_k\rho_k),\ k=0,1,
$$
where $a=\sigma\sigma^*/2$ and we use the Einstein convention for summation.
Below we use the freezing technique to show our result. Fix $z\in\mR^d$. For a function $f$, we set
$$
\tau_z f(x):=f(x+z),\ \ \ell(t,x):=\rho_0(t,x)-\rho_1(t,x).
$$
By the invariance of shifting the spatial variable $x$, we have
\begin{align*}
\p_t\tau_z\ell&=\p_i\p_j(\tau_z a_{ij}\tau_z\ell)+\div(\tau_zb_0\tau_z\ell)+\div(\tau_z(b_0-b_1)\tau_z\rho_1)\\
&=a_{ij}(t,z)\p_i\p_j\tau_z\ell+\p_i\p_j((\tau_za_{ij}-a_{ij}(t,z))\tau_z\ell)\\
&\quad+\div(\tau_z b_0\tau_z\ell)+\div(\tau_z(b_0-b_1)\tau_z\rho_1).
\end{align*}
By Duhamel's formula we have
\begin{align*}
\tau_z\ell(t,x)&=P^z_{0,t}\tau_z\ell(0,x)+\int^t_0 P^z_{s,t}(\p_i\p_j((\tau_za_{ij}-a_{ij}(s,z))\tau_z\ell))(s,x)\dif s\\
&\quad+\int^t_0 P^z_{s,t}\div(\tau_z b_0\tau_z\ell)(s,x)\dif s+\int^t_0 P^z_{s,t}\div(\tau_z(b_0-b_1)\tau_z\rho_1)(s,x)\dif s.
\end{align*}
By \eqref{SIG} and Lemma \ref{Le46} we have
\begin{align*}
|\tau_z\ell(t,0)|
&\lesssim |P^z_{0,t}\tau_z\ell(0,0)|+\int^t_0 (t-s)^{\frac{\gamma_0} 2-1}\|\tau_z\ell\|_{\mL^\infty}\dif s\\
&\quad+\sum_{i=1}^d\int^t_0 (t-s)^{-\frac12(1+|\frac1{\bbp_i}|)}\nor\tau_z b^i_0\tau_z\ell\nor_{\wt\mL^{\bbp_i}_{\x_i}}\dif s\\
&\quad+\sum_{i=1}^d\int^t_0 (t-s)^{-\frac12(1+|\frac1{\bbp_i}|)}\nor\tau_z(b^i_0-b^i_1)\tau_z\rho_1\nor_{\wt\mL^{\bbp_i}_{\x_i}}\dif s.
\end{align*}
Noting that
\begin{align*}
\nor \tau_zb^i_0\tau_z\ell\nor_{\wt\mL^{\bbp_i}_{\x_i}}
\leq \nor \tau_zb^i_0\nor_{\wt\mL^{\bbp_i}_{\x_i}}\|\tau_z\ell\|_{\mL^\infty}\leq \kappa_1\|\ell\|_{\mL^\infty},
\end{align*}
and by \eqref{CS1},
\begin{align*}
\nor \tau_z(b^i_0-b^i_1)\tau_z\rho_1\nor_{\wt\mL^{\bbp_i}_{\x_i}}
&\leq\nor\tau_z(b^i_0-b^i_1)\nor_{\wt\mL^{\bbp_i}_{\x_i}}\|\tau_z\rho_1\|_{\mL^\infty}\\
&=\nor b^i_0-b^i_1\nor_{\wt\mL^{\bbp_i}_{\x_i}}\|\rho_1\|_{\mL^\infty}\lesssim\nor b^i_0-b^i_1\nor_{\wt\mL^{\bbp_i}_{\x_i}}\|\rho^o_1\|_{\mL^\infty},
\end{align*}
we further have
\begin{align*}
\|\ell(t)\|_{\mL^\infty}=\sup_z|\tau_z\ell(t,0)|
&\lesssim \|\ell(0)\|_{\mL^\infty}+\int^t_0 (t-s)^{\frac{\gamma_0} 2-1}\|\ell(s)\|_{\mL^\infty}\dif s\\
&\quad+\sum_{i=1}^d\int^t_0 (t-s)^{-\frac12(1+|\frac1{\bbp_i}|)}\|\ell(s)\|_{\mL^\infty}\dif s\\
&\quad+\sum_{i=1}^d\int^t_0 (t-s)^{-\frac12(1+|\frac1{\bbp_i}|)}\nor b^i_0-b^i_1\nor_{\wt\mL^{\bbp_i}_{\x_i}}\dif s.
\end{align*}
By Gronwall's inequality of Volterra's type (see \cite[Lemma 2.2]{Zh10}), we obtain the desired estimate.
\end{proof}

\subsection{Well-posedness of dDDSDEs}
Now we are ready to prove the main result of this section.
\bt\label{Th215}
{\bf (Weak well-posedness)} Suppose that \eqref{SIG} holds and for any $T>0$ and $i=1,\cdots,d$, there are 
indices $(q_i,\bbp_i)\in\sI^o$ and $\bx_i\in\sX$ such that
\begin{align}\label{CQ1}
\sup_{\mu\in C([0,T];\cP(\mR^d))}\nor\sup_{r\geq 0} |b^i(\cdot,\cdot,r,\mu_\cdot)|\nor_{\wt\mL^{q_i}_T(\wt\mL^{\bbp_i}_{\x_i})}\leq \kappa_1,
\end{align}
and for some $h_i\in\mL^{q_i}_T(\wt\mL^{\bbp_i}_{\x_i})$ and for all $t,x\in[0,T]\times\mR^d$, $r,r'\geq 0$ and $\mu,\nu\in\cP(\mR^d)$,
\begin{align}\label{BZ25}
|b^i(t,x,r,\mu)-b^i(t,x,r',\nu)|\leq h_i(t,x)\big(|r-r'|+\|\mu-\nu\|_{\rm var}\big).
\end{align}
Then for any probability measure $\mu_0(\dif x)=\rho_0(x)\dif x$ with $\rho_0\in\mL^\infty$, 
there is a unique weak solution $(X,W,{\frak U})$, or equivalently, a martingale solution to dDDSDE \eqref{SDE2} with initial distribution $\mu_0$.
\et

\begin{proof}
We divide the proof into three steps.

{\bf (Step 1).} Let $\mu^0_t\equiv\mu_0$ for any $t\geq 0$. We consider the following Picard iteration: for $n\in\mN$,
\begin{align}\label{SDE5}
\dif X^n_t=b_n(t,X^n_t)\dif t+\sigma(t,X^n_t)\dif W_t,\ \ X^n_0\stackrel{(d)}{=}\mu_0,
\end{align}
where 
$$
b_n(t,x):=b(t,x,\rho^{n-1}_t(x),\mu^{n-1}_t),
$$
and
\begin{align}\label{DDD3}
\mu^{n-1}_t\mbox{ is the marginal distribution of $X^{n-1}_t$, which has a density $\rho^{n-1}_t$.}
\end{align}
By \eqref{CQ1}, one sees that for each $i=1,\cdots, d$,
\begin{align}\label{DD3}
\sup_n\nor b^i_n\nor_{\wt\mL^{q_i}_T(\wt\mL^{\bbp_i}_{\x_i})}\leq \kappa_1.
\end{align}
Thus, by Theorem \ref{Main1}, for each $n\in\mN$, there is a unique weak solution 
$(X^n,W^n,{\frak U}^n)$ to SDE \eqref{SDE5}, where
$$
{\frak U}^n:=(\Omega^n,\sF^n,\bP_n;(\sF^n_t)_{t\geq 0}),
$$  
and for each $t>0$, $X^n_t$ admits a density $\rho^n_t$ satisfying the following estimate: for all $(t,y)\in[0,T]\times\mR^d$,
\begin{align}\label{BOU}
\rho^n_t(y)\leq \frac{C_1}{t^{d/2}}\int_{\mR^d}\e^{-\frac{|x-y|^2}{\delta_1t}}\rho_0(x)\dif x\lesssim\|\rho_0\|_\infty.
\end{align}
Moreover,  for any $T>0$,  by  \eqref{Mom},
there is a constant $C>0$ such that
$$
\sup_n\mE^{\bP_n}|X^n_t-X^n_s|^4\leq C|t-s|^2,\ \ s,t\in[0,T],
$$
and by \eqref{Kry11}, for any $(q_0,\bbp_0)\in\sI_2$,
there is a constant $C>0$ such that for all $f\in\wt\mL^{q_0}_T(\widetilde\mL^{\bbp_0}_\x)$,
\begin{align}\label{Kry13}
\sup_n\mE^{\bP_n}\left(\int^T_0f(s,X^n_s)\dif s\right)\leq C\nor f\nor_{\wt\mL^{q_0}_T(\widetilde\mL^{\bbp_0}_\x)}.
\end{align}
In particular, by Kolmogorov's criterion,
\begin{align}\label{BZ17}
\mbox{the laws $\mP_n$ of $X^n_\cdot$ in $\mC_T$ are tight.}
\end{align}

{\bf (Step 2).} For simplicity of notations, we write
$$
\Gamma_{n,m}(t):=\|\rho^{n}_t-\rho^{m}_t\|_{\mL^\infty}+\|\rho^{n}_t-\rho^{m}_t\|_{\mL^1}.
$$
Noting that by \eqref{BZ25} and \eqref{DDD3},
$$
|b^i_n(s,x)-b^i_m(s,x)|\leq h_i(s,x)\big(|\rho^{n-1}_s(x)-\rho^{m-1}_s(x)|+\|\mu^{n-1}_s-\mu^{m-1}_s\|_{\rm var}\big)\leq h_i(s,x)\Gamma_{n-1,m-1}(s),
$$
we have
\begin{align}\label{BOU1}
\nor b^i_n(s)-b^i_m(s)\nor_{\wt\mL^{\bbp_i}_{\x_i}}\leq \nor h_i(s)\nor_{\wt\mL^{\bbp_i}_{\x_i}}\Gamma_{n-1,m-1}(s)=:\ell_i(s) \Gamma_{n-1,m-1}(s).
\end{align}
Since $(\frac{q_i}2,\frac{\bbp_i}2)\in\sI_2$, by Lemma \ref{Lem34} and \eqref{Kry13}, \eqref{BOU1}, we have
\begin{align*}
\cH(\mu^n_t|\mu^m_t)&\leq\frac12\mE^{\mP_m}\left(\int^t_0|\sigma^{-1}(s,w_s)(b_n(s,w_s)-b_m(s,w_s))|^2\dif s\right)\\
&\leq\frac{\|\sigma^{-1}\|^2_\infty}{2}\mE^{\mP_m}\left(\int^t_0|b_n(s,w_s)-b_m(s,w_s)|^2\dif s\right)\\
&\lesssim \sum_{i=1}^d\left(\int^t_0\nor |b^i_n(s)-b^i_m(s)|^2\nor^{q_i/2}_{\wt\mL^{\bbp_i/2}_{\x_i}}\dif s\right)^{\frac2{q_i}}\\
&=\sum_{i=1}^d\left(\int^t_0\nor b^i_n(s)-b^i_m(s)\nor^{q_i}_{\wt\mL^{\bbp_i}_{\x_i}}\dif s\right)^{\frac2{q_i}}\\
&\lesssim \sum_{i=1}^d\left(\int^t_0\ell^{q_i}_i(s)\Gamma^{q_i}_{n-1,m-1}(s)\dif s\right)^{\frac2{q_i}}.
\end{align*}
By Pinsker's inequality \eqref{Pin0}, we get
\begin{align}\label{CS3}
\|\rho^n_t-\rho^m_t\|_{\mL^1}=\|\mu^{n}_t-\mu^{m}_t\|_{\rm var}
\lesssim \sum_{i=1}^d\left(\int^t_0\ell^{q_i}_i(s)\Gamma^{q_i}_{n-1,m-1}(s)\dif s\right)^{\frac{1}{q_i}}.
\end{align}
On the other hand, by \eqref{CS2}, \eqref{BOU1} and H\"older's inequality, for $q'_i=\frac{q_i}{q_i-1}$, we have
\begin{align*}
\|\rho^{n}_t-\rho^{m}_t\|_{\mL^\infty}
&\lesssim 
\sum_{i=1}^d\int^t_0 (t-s)^{-\frac12(1+|\frac1{\bbp_i}|)}\ell_i(s)\Gamma_{n-1,m-1}(s)\dif s\\
&\lesssim 
\sum_{i=1}^d\left(\int^t_0 (t-s)^{-\frac{q_i'}2(1+|\frac1{\bbp_i}|)}\dif s\right)^{\frac1{q_i'}}
\left(\int^t_0\ell^{q_i}_i(s)\Gamma^{q_i}_{n-1,m-1}(s)\dif s\right)^{\frac{1}{q_i}}\\
&\lesssim \sum_{i=1}^d \left(\int^t_0\ell^{q_i}_i(s)\Gamma^{q_i}_{n-1,m-1}(s)\dif s\right)^{\frac{1}{q_i}},
\end{align*}
which together with \eqref{CS3} yields
\begin{align*}
\Gamma_{n,m}(t)\lesssim \sum_{i=1}^d \left(\int^t_0\ell^{q_i}_i(s)\Gamma^{q_i}_{n-1,m-1}(s)\dif s\right)^{\frac{1}{q_i}}.
\end{align*}
Let $q=q_1\vee \cdots\vee q_d$. By H\"older's inequality with respect to $\ell^{q_i}_i(s)\dif s$,  we get
\begin{align*}
\Gamma_{n,m}^q(t)&\lesssim\sum_{i=1}^d\left(\int^t_0\ell^{q_i}_i(s)\Gamma_{n-1,m-1}^q(s)\dif s\right)\left(\int^t_0\ell^{q_i}_i(s)\dif s\right)^{\frac{q}{q_i}-1}\\
&\lesssim\int^t_0\sum_{i=1}^d\ell^{q_i}_i(s)\Gamma_{n-1,m-1}^q(s)\dif s.
\end{align*}
Therefore, by \eqref{BOU} and Fatou's lemma,
$$
\varlimsup_{n,m\to\infty}\Gamma_{n,m}^q(t)\lesssim \int^t_0\sum_{i=1}^d\ell^{q_i}_i(s)\varlimsup_{n,m\to\infty}\Gamma_{n-1,m-1}^q(s)\dif s,
$$
which implies by Gronwall's inequality  that  for each $t\in[0,T]$,
\begin{align}\label{LL1}
\varlimsup_{n,m\to\infty}\big(\|\rho^{n}_t-\rho^{m}_t\|_{\mL^\infty}+\|\rho^{n}_t-\rho^{m}_t\|_{\mL^1}\big)=\varlimsup_{n,m\to\infty}\Gamma_{n,m}^q(t)=0.
\end{align}
Now by \eqref{BZ17}, there is a subsequence $n_k$ such that as $k\to\infty$,
 $$
 \mP_{n_k}\mbox{ weakly converges to some $\mP\in\cP(\mC_T)$},
 $$
 and by \eqref{LL1}, $\mP\circ w^{-1}_t(\dif x)=\mu_t(\dif x)=\rho_t(x)\dif x$ and for each $t\in(0,T]$,
\begin{align}\label{LL2}
\varlimsup_{n\to\infty}\big(\|\rho^{n}_t-\rho_t\|_{\mL^\infty}+\|\rho^{n}_t-\rho_t\|_{\mL^1}\big)=0.
\end{align}

{\bf (Step 3).} In this step we show $\mP\in\cM^{\sigma,b}_{\mu_0}$. More precisely,
we want to show that for fixed $f\in C^2_c(\mR^d)$, the process $M^f_t$ defined by \eqref{BZ24}
is a $\cB_t$-martingale under $\mP$, that is, for any $t_0<t_1$ and every bounded $\cB_{t_0}$-measurable continuous function $\eta$,
\begin{align}\label{BZ20}
\mE \((M^f_{t_1}-M^f_{t_0})\eta\)=0.
\end{align}
Note that for each $k\in\mN$, by SDE \eqref{SDE5} and It\^o's formula,
$$
\mE^{\mP_{n_k}} \((M^k_{t_1}-M^k_{t_0})\eta\)=0,
$$
where
$$
M^k_t:=f(w_t)-f(w_0)-\int^t_0\Big(\tr(a_{n_k}\cdot\nabla^2f)+b_{n_k}\cdot\nabla f \Big)(s,w_s)\dif s.
$$
Since $x\mapsto a_{n_k}(s,x)$ is continuous, to show \eqref{BZ20}, the key point is to prove the following:
$$
\lim_{k\to\infty}\mE^{\mP_{n_k}}\left(\eta\int^{t_1}_{t_0}b_{n_k}(s,w_s)\cdot\nabla f(s,w_s)\dif s\right)
=\mE\left(\eta\int^{t_1}_{t_0}b(s,w_s,\rho_s(w_s),\mu_s)\cdot\nabla f(s,w_s)\dif s\right),
$$
which follows from:
\begin{align}\label{BZ22}
\lim_{m\to\infty}\sup_k\mE^{\mP_{n_k}}\left(\int^{t_1}_{t_0}|b_{n_m}(s,w_s)-b(s,w_s,\rho_s(w_s),\mu_s)|\dif s\right)=0,
\end{align}
together with
\begin{align}\label{BZ23}
\lim_{k\to\infty}\mE^{\mP_{n_k}}\left(\eta\int^{t_1}_{t_0}b_{n_m}(s,w_s)\cdot\nabla f(w_s)\dif s\right)=
\mE^{\mP}\left(\eta\int^{t_1}_{t_0}b_{n_m}(s,w_s)\cdot\nabla f(w_s)\dif s\right)
\end{align}
for each $m\in\mN$.
The first limit \eqref{BZ22} follows by the Krylov estimates \eqref{Kry13}, \eqref{BZ25} and \eqref{LL2}.
For the second, let
$$
b^\eps_{n_m}(s,x):=b_{n_m}(s,\cdot)*\varGamma_\eps(x),\ \ \eps\in(0,1),
$$
where $\varGamma_\eps$ is the mollifiers in \eqref{Mol}. For each $\eps\in(0,1)$, since $x\mapsto b^\eps_{n_m}(s,x)$ is bounded continuous, by the weak convergence of $\mP_{n_k}$, we have
\begin{align}\label{BZ231}
\lim_{k\to\infty}\mE^{\mP_{n_k}}\left(\eta\int^{t_1}_{t_0}b^\eps_{n_m}(s,w_s)\cdot\nabla f(w_s)\dif s\right)=
\mE^{\mP}\left(\eta\int^{t_1}_{t_0}b^\eps_{n_m}(s,w_s)\cdot\nabla f(w_s)\dif s\right).
\end{align}
Moreover, for each $m\in\mN$ and $R>0$, by the Krylov estimate \eqref{Kry13}, we also have
\begin{align}\label{BZ232}
\begin{split}
&\lim_{\eps\to0}\sup_k\mE^{\mP_{n_k}}\left(\int^{t_1}_{t_0}|b^\eps_{n_m}-b_{n_m}|(s,w_s)|\b1_{|w_s|\leq R}\dif s\right)\\
&\qquad\lesssim \lim_{\eps\to0}\sum_{i=1}^d\|(b^\eps_{n_m}-b_{n_m})^i\b1_{B^R_0}\|_{\mL^{q_i}_T(\mL^{\bbp_i}_{\x_i})}=0,
\end{split}
\end{align}
and
\begin{align}\label{BZ233}
\lim_{R\to\infty}\sup_{k,\eps}\mE^{\mP_{n_k}}\left(\int^{t_1}_{t_0}|b^\eps_{n_m}-b_{n_m}|(s,w_s)|\b1_{|w_s|\geq R}\dif s\right)=0.
\end{align}
Combining \eqref{BZ231},  \eqref{BZ232} and  \eqref{BZ233}, we obtain \eqref{BZ23}. Thus we  complete the proof of existence. On the other hand, by the same calculations as in \eqref{LL1}, one can show 
that any two weak solutions have the same marginal distribution. Then by Theorem \ref{Main1}, we get the weak uniqueness.
\end{proof}
\br\rm
If $b$ does not depend on the density variable $r$, then we can drop the assumption $\mu_0(\dif x)=\rho_0(x)\dif x$.
In this case, we can only use \eqref{CS3} to show that $\mu^n_t$ is a Cauchy sequence. We note that a similar 
result has been shown in \cite{Wa21}. However, even in the non-mixed norm case,
the results in \cite{Wa21} do not cover our case since we are using the total variational norm in \eqref{BZ25}. 
Moreover, our proofs are based on the Fokker-Planck equation, and Wang's proofs are based on the backward Kolmogorov equation.
\er
\bt\label{Th216}
{\bf (Strong  well-posedness)} In addition to the assumptions in Theorem \ref{Th215}, we also assume \eqref{SIG0}. Then 
there is a unique strong solution.
\et
\begin{proof}
This is a direct consequence of Theorems \ref{Th215} and  \ref{Th217}. 
\end{proof}

\section{Weak convergence of propagation of chaos}

Throughout this section we assume {\bf (H$^\sigma$)} and {\bf (H$^b$)}. 
Let
$$
\bX^N_t:=(X^{N,1}_t,\cdots, X^{N,N}_t),\ \ \bW^N_t:=(W^1,\cdots, W^N),
$$
and for $\x=(x^1,\cdots,x^N)$, define
\begin{align}\label{BB1}
B(t,\x):=\left(F\Big(t,x^1,\frac1N\sum_{j=1}^N\phi_t(x^1,x^j)\Big),\cdots,F\Big(t,x^N,\frac1N\sum_{j=1}^N\phi_t(x^N,x^j)\Big)\right),
\end{align}
and a $(dN)\times (dN)$-matrix ${\boldsymbol\sigma}$ by
\begin{align}\label{BB2}
{\boldsymbol\sigma}(t,\x):={\rm diag}_N(\sigma(t, x^1),\cdots,\sigma(t, x^N)).
\end{align}
Then the particle system \eqref{MV2} can be written as an SDE in $\mR^{dN}$:
$$
\dif\bX^N_t=B(t,\bX^N_t)\dif t+{\boldsymbol \sigma}(t, \bX^N_t)\dif \bW^N_t.
$$
Noting that by {\bf (H$^b$)},
$$
|B_i(t,\x)|\leq h(t,x^i)+\frac{\kappa_1}N\sum_{j=1}^N|\phi_t(x^i,x^j)|,
$$
we have for $\vec{\bbp}=(\infty,\cdots,\infty,\bbp)\in [1,\infty]^{dN}$ and for $\x_i=(\cdots, x^{i-1}, x^{i+1},\cdots, x^N, x^i)$,
\begin{align*}
\nor B_i\nor_{\mL^q_T(\wt\mL^{\vec{\bbp}}_{\x_i})}\leq\nor h\nor_{\mL^q_T(\wt\mL^{\bbp}_{\x})}
+\kappa_1\left[\int^T_0\sup_{y\in\mR^d}\nor\phi_t(\cdot,y)\nor^{q}_{\wt\mL^{\bbp}_{\x}}\dif t\right]^{\frac1q}<\infty.
\end{align*}
Then, by Theorem \ref{Th217}, for any initial value $\bX^N_0$, there is a unique strong solution to the above SDE.
In particular, there is a measurable functional $\Phi: \mR^{dN}\times\mC_T^{N}\to\mC_T^{N}$ such that
\begin{align}\label{STR}
\bX^N_t=\Phi(\bX^N_0,\bW^N_\cdot)(t),\ \ t\in[0,T].
\end{align}
\subsection{Martingale approach}
In this section we use the classical martingale approach to show the following qualitative result of weak convergence.
\bt\label{Th31}
For any $N\in\mN$,
let $\xi^N_1,\cdots,\xi^N_N$ be N-random variables and $\mu_0\in\cP(\mR^d)$. Suppose that
the law of $(\xi^N_1,\cdots,\xi^N_N)$ is invariant under any permutation of $\{1,\cdots,N\}$, and for any $k\leq N$,
\begin{align}\label{CC94}
\mP\circ\(\xi^N_1,\cdots,\xi^N_k\)^{-1}\to \mu_0^{\otimes k},\ \ N\to\infty.
\end{align}
Then for any $k\leq N$ and $T>0$,
\begin{align}\label{CC14}
\mP\circ\(X^{N,1}_{[0,T]},\cdots,X^{N,k}_{[0,T]}\)^{-1}\to \mu_{[0,T]}^{\otimes k}, \ \ N\to\infty,
\end{align}
where $\mu_{[0,T]}$ is the law of the unique solution of dDDSDE \eqref{MV1} with initial distribution $\mu_0$ on $\mC_T$.
\et

First of all, we use the partial Girsanov transform as used in \cite{JTT18,To20} 
to show some uniform Krylov estimate for particle system \eqref{MV2}.
Let $\{\W^i_t, i\in\mN\}$ be a sequence of independent $d$-dimensional standard Brownian motions.
For each $x\in\mR^d$, let $Z_t(x)$ be the unique strong solution of the following SDE starting from $x$:
$$
\dif Z_t=\sigma\(t,Z_t\)\dif\W^1_t,\ \ Z_0=x.
$$
For each ${\boldsymbol z}=(z^2,\cdots,z^N)\in\mR^{(N-1)d}$, 
let ${\boldsymbol Z}^N_t({\boldsymbol z}):={\boldsymbol Z}^N_t:=(Z^{N,2}_t,\cdots,Z^{N,N}_t)$ be the unique strong solution of the following SDE starting from ${\boldsymbol z}$:
$$
\dif Z^{N,k}_t=b\(t,Z^{N,k}_t,\eta_{{\boldsymbol Z}^N_t}\)\dif t+\sigma\(t,Z^{N,k}_t\)\dif\W^k_t,\ \ Z^{N,k}_0=z^k,
$$
where $k=2,\cdots,N$ and
$$
\eta_{{\boldsymbol z}}(\dif y):=\frac1N\sum_{j=2}^N\delta_{z^j}(\dif y).
$$
In particular, as Brownian functionals of $\wt W^1$ and $(\wt W^2,\cdots,\wt W^N)$ respectively,
\begin{align}\label{13}
Z_\cdot(\cdot)\mbox{ is independent of }{\boldsymbol Z}^N_\cdot(\cdot),
\end{align}
and by the notion of strong solution of SDEs (see \eqref{STR}),
\begin{align}\label{23}
\X^{N,1}_t:=Z_t(\xi^N_1),\ \ (\X^{N,2}_t,\cdots,\X^{N,N}_t):={\boldsymbol Z}^N_t(\xi^N_2,\cdots,\xi^N_N)=:{\boldsymbol Y}^N_t,
\end{align}
solves the following SDE:
\begin{align}\label{NN3}
\left\{
\begin{aligned}
&\dif \X^{N,1}_t=\sigma\(t,\X^{N,1}_t\)\dif\W^1_t,\ \ \X^{N,1}_0=\xi^N_1,\\
&\quad\mbox{ and for each $k=2,\cdots,N,$}\\
&\dif \X^{N,k}_t=b\(t,\X^{N,k}_t,\eta_{{\boldsymbol Y}^N_t}\)\dif t+\sigma\(t,\X^{N,k}_t\)\dif\W^k_t,\ \ \X^{N,k}_0=\xi^N_k,
\end{aligned}
\right.
\end{align}
where
$$
\eta_{{\boldsymbol Y}^N_t}:=\frac1{N}\sum_{j=2}^N\delta_{\X^{N,j}_t}(\dif y).
$$
Now let us define
$$
\eta_{\wt\bX^{N}_t}(\dif y):=\frac1{N}\sum_{j=1}^N\delta_{\X^{N,j}_t}(\dif y),\ \ H^{N,1}_t:=\sigma\(t,\X^{N,1}_t\)^{-1}b\(t,\X^{N,1}_t,\eta_{\wt\bX^{N}_t}\),
$$
and for $k=2,\cdots,N$,
$$
H^{N,k}_t:=\sigma\(t,\X^{N,k}_t\)^{-1}\Big[b\(t,\X^{N,k}_t,\eta_{\wt\bX^{N}_t}\)-b\(t,\X^{N,k}_t,\eta_{{\boldsymbol Y}^N_t}\)\Big].
$$
By the above definition, we clearly have for each $i=1,\cdots,N$,
\begin{align}\label{MV22}
\dif \X^{N,i}_t=b\(t,\X^{N,i}_t,\eta_{\wt\bX^{N}_t}\)\dif t+\sigma\(t,\X^{N,i}_t\)\(\dif\W^i_t-H^{N,i}_t\dif t\).
\end{align}

The following uniform estimate is the key step for performing the Girsanov transform to derive the Krylov estimate for the particle system,
whose proof strongly depends on the independence in \eqref{13} and the strong uniqueness used in {\eqref{NN3}}.
\bl\label{Le52}
For any $\gamma, T>0$, 
\begin{align}\label{111}
\sup_N\mE\exp\left\{\gamma\sum_{i=1}^N\int^T_0|H^{N,i}_t|^2\dif t\right\}<\infty.
\end{align}
\el
\begin{proof}
For $x\in\mR^d$ and ${\boldsymbol y}=(y^2,\cdots,y^N)\in\mR^{(N-1)d}$, let us
write $\eta_{{\boldsymbol y}}:=\frac1N\sum_{j=2}^N\delta_{y^j}$ and define
$$
\Gamma_1(t,x,{\boldsymbol y}):=\sigma(t,x)^{-1}b\left(t,x,\frac{\delta_x}{N}+\eta_{{\boldsymbol y}}\right),
$$
and for $k=2,\cdots,N$,
$$
\Gamma_k(t,x,{\boldsymbol y}):=\sigma(t,y^k)^{-1}\left[b\left(t,y^k,\frac{\delta_{x}}{N}+\eta_{{\boldsymbol y}}\right)
-b\left(t,y^k, \eta_{{\boldsymbol y}}\right)\right].
$$
From the very definition, one sees that for each $i=1,\cdots,N$,
$$
H^{N,i}_s=\Gamma_i\(s,\X^{N,1}_s,{\boldsymbol Y}^N_s\),
$$
and by \eqref{23} and \eqref{13},
\begin{align}
&\mE\exp\left\{\gamma\sum_{i=1}^N\int^T_0|H^{N,i}_s|^2\dif s\right\}
=\mE\exp\left\{\gamma\int^T_0\sum_{i=1}^N|\Gamma_i\(s,Z_s(\xi^N_1),{\boldsymbol Y}^N_s\)|^2\dif s\right\}\no\\
&\qquad=\mE\left(\mE\exp\left\{\gamma\int^T_0\sum_{i=1}^N|\Gamma_i\(s,Z_s(x),{\boldsymbol y}_s\)|^2\dif s\right\}
\Big|_{(x,{\boldsymbol y_\cdot})=(\xi^N_1,{\boldsymbol Y}^N_\cdot)}\right)\no\\
&\qquad\leq\sup_{x,{\boldsymbol y_\cdot}}\mE\exp\left\{\gamma\int^T_0\sum_{i=1}^N|\Gamma_i\(s,Z_s(x),{\boldsymbol y}_s\)|^2\dif s\right\}\no\\
&\qquad=\sup_{x,{\boldsymbol y_\cdot}}\mE\exp\left\{\gamma\int^T_0f_{{\boldsymbol y}}\(s,Z_s(x)\)\dif s\right\},\label{AS66}
\end{align}
where for ${\boldsymbol y}=({\boldsymbol y}_s)_{s\in[0,T]}$,
$$
f_{{\boldsymbol y}}(s,x):=\sum_{i=1}^N|\Gamma_i\(s,x,{\boldsymbol y}_s\)|^2.
$$
Note that by \eqref{CC1} and because $\phi_t(x,x)=0$,
\begin{align*}
|\Gamma_1(t,x,{\boldsymbol y})|&=\left|\sigma(t,x)^{-1}F\left(t,x,\frac1{N}\Big(\phi_t(x,x)+\sum_{j=2}^N|\phi_t(x,y^j)|\Big)\right)\right|\\
&\leq 
\|\sigma^{-1}\|_\infty\left(h(t,x)+\frac{\kappa_1}{N}\sum_{j=2}^N\phi_t(x,y^j)\right),
\end{align*}
and
\begin{align*}
|\Gamma_k(t,x,{\boldsymbol y})|
&\leq\frac{\kappa_1\|\sigma^{-1}\|_\infty}{N}|\phi_t(y^k,x)|,
\end{align*}
and by \eqref{CC2},
\begin{align*}
\left(\int^T_0\sup_{\boldsymbol y}\nor\Gamma_1(t,\cdot,{\boldsymbol y})\nor_{\wt\mL^\bbp_\x}^{q}\dif t\right)^{1/q}\leq \|\sigma^{-1}\|_\infty(\kappa_1+\kappa^2_1)
\end{align*}
and
\begin{align*}
\left(\int^T_0\sup_{\boldsymbol y}\nor\Gamma_k(t,\cdot,{\boldsymbol y})\nor_{\wt\mL^\bbp_\x}^{q}\dif t\right)^{1/q}\leq \frac{\kappa^2_1\|\sigma^{-1}\|_\infty}{N-1}.
\end{align*}
From these two estimates, by Minkowskii's inequality, we derive 
\begin{align*}
\left(\int^T_0\sup_{\boldsymbol y}\nor f_{{\boldsymbol y}_\cdot}(s,\cdot)\nor_{\wt\mL^{\bbp/2}_\x}^{q/2}\dif t\right)^{2/q}
&\leq \sum_{i=1}^N\left(\int^T_0\sup_{\boldsymbol y}\nor|\Gamma_i(t,\cdot,{\boldsymbol y})|^2\nor_{\wt\mL^{\bbp/2}_\x}^{q/2}\dif t\right)^{2/q}\\
&=\sum_{i=1}^N\left(\int^T_0\sup_{\boldsymbol y}\nor\Gamma_i(t,\cdot,{\boldsymbol y})\nor_{\wt\mL^\bbp_\x}^{q}\dif t\right)^{2/q}\\
&\leq \|\sigma^{-1}\|^2_\infty\left((\kappa_1+\kappa^2_1)^2+\frac{\kappa^4_1}{N}\right).
\end{align*}
Thus, because $(\frac q2,\frac\bbp 2)\in\sI_2$, by \eqref{Kha1} we have
\begin{align*}
\sup_{x,{\boldsymbol y}}\mE\exp\left\{\gamma\int^T_0f_{{\boldsymbol y}}\(s,Z_s(x)\)\dif s\right\}\leq C,
\end{align*}
which together with \eqref{AS66} yields \eqref{111}.  
\end{proof}

Now if we define 
$$
\sE^N_t:=\exp\left\{\sum_{i=1}^N\int^t_0H^{N,i}_s\dif \W^i_s-\frac12\sum_{i=1}^N\int^t_0|H^{N,i}_s|^2\dif s\right\},
$$
then by \eqref{111} and Novikov's criterion, $t\mapsto Z_t$ is an exponential martingale  and
$$
\sE^N_t=1+\sum_{i=1}^N\int^t_0H^{N,i}_s\sE^N_s\dif \W^i_s.
$$
Thus, by Girsanov's theorem, $\big(\W^i_t-\int^t_0H^{N,i}_s\dif s\big)_{t\in[0,T]}^{i=1,\cdots,N}$ are $N$-independent standard Brownian motions under the new probability measure
$$
\mQ:=\sE^N_T\mP.
$$
Moreover, by \eqref{MV22} and the weak uniqueness for SDE \eqref{MV2}, we have
\begin{align}\label{NN1}
\mQ\circ\big(\bX^{N}_{[0,T]}\big)^{-1}=\mP\circ\big(\bX^{N}_{[0,T]}\big)^{-1},
\end{align}
and for any $\gamma\in\mR$, by \eqref{111} it is standard to derive that
\begin{align}\label{NN2}
\sup_N\mE\left(\sup_{t\in[0,T]}|\sE^N_t|^\gamma\right)<\infty.
\end{align}
From these, we can derive the following crucial Krylov estimate for the particle system.
\bl\label{Le33}
\begin{enumerate}[(i)]
\item The law of $(X^{N,1}_t)_{t\in[0,T]}$, $N\in\mN$, in $\mC_T$ is tight.

\item For any $T>0$, $(q,\bbp)\in\sI_2$ and $\x\in\sX$, there is a constant $C_1=C_1(T,\Theta)>0$ such that for any $f\in\wt\mL^q_T(\wt\mL^\bbp_\x)$,
\begin{align}\label{CC9}
\sup_N\mE\left(\int^T_0f\(t,X^{N,1}_t\)\dif t\right)\leq C_1\nor f\nor_{\wt\mL^q_T(\wt\mL^\bbp_\x)},
\end{align}
and for any $\lambda>0$ and $\beta\in(0,2-|\frac1\bbp|-\frac2q)$, there is a $C_2=C_2(T,\Theta,\lambda,\beta)>0$ such that for any $f\in\wt\mL^q_T(\wt\mL^\bbp_\x)$,
\begin{align}\label{CC98}
\sup_N\mE\exp\left\{\lambda\int^T_0f\(t,X^{N,1}_t\)\dif t\right\}\leq\e^{C_2 \nor f\nor^{2/\beta}_{\widetilde\mL^q_T(\wt\mL^\bbp_\x)}}.
\end{align}
\item Let $\bbp_1,\bbp_2\in(1,\infty)^d$ and let $q\in(1,\infty)$ with $|\frac{1}{\bbp_1}|+|\frac{1}{\bbp_2}|+\frac{2}{q}<2$
and $\x_1,\x_2\in\sX$.
Then for any $T>0$, it holds that for some $C_3=C_3(T,\Theta)>0$,
\begin{align}\label{CC10}
\sup_N\mE\left(\int^T_0f\(t,X^{N,1}_t, X^{N,2}_t\)\dif t\right)\leq C_3\nor f\nor_{\wt\mL^q_T(\wt\mL^{\bbp_1}_{\x_1}(\wt\mL^{\bbp_2}_{\x_2}))},
\end{align}
where $\wt\mL^q_T(\wt\mL^{\bbp_1}_{\x_1}(\wt\mL^{\bbp_2}_{\x_2}))$ is the localization of $\mL^q_T(\mL^{\bbp_1}_{x_1}(\mL^{\bbp_2}_{\x_2}))$ as in \eqref{LOC1}.
\end{enumerate}
\el
\begin{proof}
(i) By \eqref{NN1}, H\"older's inequality, \eqref{NN2}  and \eqref{NN3}, there is a constant $C>0$ such that for all $0\leq s<t\leq T$ and $N\in\mN$,
\begin{align*}
\mE|X^{N,1}_t-X^{N,1}_s|^4&=\mE\(\sE^N_T|\X^{N,1}_t-\X^{N,1}_s|^{4}\)\\
&\leq \(\mE (\sE^N_T)^2\)^{1/2}\(\mE|\X^{N,1}_t-\X^{N,1}_s|^8\)^{1/2}\leq C|t-s|^2,
\end{align*}
which, together with \eqref{CC94}, implies the tightness by Kolmogorov's criterion.

(ii) Let $\gamma>1$ be such that $(\frac{q}\gamma,\frac\bbp\gamma)\in\sI_2$. 
By \eqref{NN1}, H\"older's inequality, \eqref{NN2} and \eqref{Kry10}, we have
\begin{align*}
\mE\left(\int^T_0 f(t,X^{N,1}_t)\dif t\right)&=\mE\left(\sE^N_T\int^T_0 f(t,\X^{N,1}_t)\dif t\right)\\
&\leq\Big[\mE (\sE^N_T)^{\frac{\gamma}{\gamma-1}}\Big]^{1-1/\gamma}\left[\mE\left(\int^T_0 |f(t,\X^{N,1}_t)|^\gamma\dif t\right)\right]^{1/\gamma}\\
&\leq C\nor f^\gamma\nor^{1/\gamma}_{\wt\mL^{q/\gamma}_T(\wt\mL^{\bbp/\gamma}_\x)}= C\nor f\nor_{\wt\mL^q_T(\wt\mL^\bbp_\x)}.
\end{align*}
\eqref{CC98} follows by the same method and \eqref{Kha1}.

(iii) Let $\gamma\in(1,\min_i (p_{1i},p_{2i})\wedge q)$ be such that $|\frac{1}{\bbp_1/\gamma}|+|\frac{1}{\bbp_2/\gamma}|+\frac2{q/\gamma}<2$. 
By \eqref{NN1}, H\"older's inequality and \eqref{NN2}, we have
\begin{align*}
\mE\left(\int^T_0 f(t,X^{N,1}_t,X^{N,2}_t)\dif t\right)
&=\mE\left(\sE^N_T\int^T_0 f(t,\X^{N,1}_t,\X^{N,2}_t)\dif t\right)\\
&\leq\Big[\mE (\sE^N_T)^{\frac{\gamma}{\gamma-1}}\Big]^{\frac{\gamma-1}\gamma}\left[\mE\left(\int^T_0 |f(t,\X^{N,1}_t,\X^{N,2}_t)|^\gamma\dif t\right)\right]^{\frac 1\gamma}\\
&\lesssim \sup_{x}\left[\mE\left(\int^T_0 |f(t,Z_t(x),\X^{N,2}_t)|^\gamma\dif t\right)\right]^{\frac1\gamma}.
\end{align*}
By $|\frac{1}{\bbp_1/\gamma}|+|\frac{1}{\bbp_2/\gamma}|+\frac2{q/\gamma}<2$, one can choose $q_1,q_2>\gamma$ so that
$\frac1{q_1/\gamma}+\frac{1}{q_2/\gamma}=1+\frac1{q/\gamma}$ and $(q_i/\gamma,\bbp_i/\gamma)\in\sI_2$, $i=1,2$.
Since $Z_\cdot(x)$ and $\X^{N,2}$ are independent by \eqref{13} and \eqref{23}, and satisfy the Krylov estimate \eqref{CC9}, the desired estimate now follows by using \cite[Lemma 2.6]{RZ21}.
\end{proof}

In the following, in order to take weak limits, we need to mollify the coefficients. For $\eps\in(0,1)$ and $k\in\mN$, we define
\begin{align}\label{CC13}
b_{\eps,k}(t,x,\mu):=F_\eps(t,x,(\phi^k_t\circledast\mu)(x)),
\end{align}
where
$$
F_\eps(t,x,r):=(-\eps^{-1})\vee \((F(t,\cdot,r)*\varGamma_\eps)(x)\)\wedge\eps^{-1}
$$
and
$$
\phi^k_t(x,y):=(-k)\vee \((\phi_t*\varGamma_{1/k})(x,y)\)\wedge k.
$$

We have the following properties for the above approximation.
\bl
(i) $b_{\eps,k}\in L^\infty_T(C_b(\mR^d\times\cP(\mR^d)))$ and
$$
|b_{\eps,k}(t,x,\mu)|\leq h_t*\varGamma_\eps(x)+\kappa_0(\phi^k_t\circledast\mu)*\varGamma_\eps(x)
$$
and
$$
|b-b_{\eps,k}|\(t,x,\mu\)\leq\sup_{|r|\leq k}|F_\eps-F|(t,x,r)
+\kappa_0|(\phi^k_t-\phi_t)\otimes\mu|(x).
$$
(ii) For any $T>0$,
\begin{align}\label{CC11}
\lim_{k\to\infty }\lim_{\eps\to 0}\sup_N\mE\left(\int^T_0|b-b_{\eps,k}|\(s, X^{N,1}_s,\eta_{\bX^N_s}\)\dif s\right)=0.
\end{align}
\el
\begin{proof}
(i) is obvious by definition and the assumptions. We now show (ii). Note that
\begin{align}\label{CC12}
|b-b_{\eps,k}|\(s, X^{N,1}_s,\eta_{\bX^N_s}\)\leq\sup_{|r|\leq k}|F_\eps-F|(s, X^{N,1}_s,r)
+\frac{\kappa_0}N\sum_{j=1}^N|\phi^k_s-\phi_s|(X^{N,1}_s,X^{N,j}_s).
\end{align}
We first show that for fixed $r\in\mR^m$,
\begin{align}\label{AD9}
\lim_{\eps\to 0}\sup_N\mE\left(\int^T_0|F_\eps-F|(s, X^{N,1}_s,r)\dif s\right)=0.
\end{align}
Let $R>0$. Since $(\frac q2,\frac\bbp2)\in\sI_2$ and $\nor F(\cdot,r)\nor^2_{\wt\mL^{q}_T(\wt\mL^{\bbp}_\x)}<\infty$ by \eqref{CC1} and \eqref{CC2}, 
by H\"older's inequality and \eqref{CC9}, \eqref{BZ130}, we have
\begin{align}
&\mE\left(\int^T_0|F_\eps-F|(s, X^{N,1}_s,r)\b1_{|X^{N,1}_s|>R}\dif s\right)\no\\
&\quad\leq\left[\mE\left(\int^T_0|F_\eps-F|^2(s, X^{N,1}_s,r)\dif s\right)\right]^{\frac12}\left[\int^T_0\mP(|X^{N,1}_s|>R)\dif s\right]^{\frac12}\no\\
&\quad\lesssim\nor|F_\eps-F|^2(\cdot,r)\nor^{1/2}_{\wt\mL^{q/2}_T(\wt\mL^{\bbp/2}_\x)}
\left[\int^T_0\Big(\mP\big(|X^{N,1}_s-X^{N,1}_0|>\tfrac R2\big)+\mP\big(|X^{N,1}_0|>\tfrac R2\big)\Big)\dif s\right]^{\frac12}\no\\
&\quad\lesssim\nor F(\cdot,r)\nor_{\wt\mL^{q}_T(\wt\mL^{\bbp}_\x)}
\left[\int^T_0\left(\frac{\mE|X^{N,1}_s-X^{N,1}_0|}{R}+\mP\big(|X^{N,1}_0|>\tfrac R2\big)\right)\dif s\right]^{\frac12}\no\\
&\quad\lesssim\nor F(\cdot,r)\nor_{\wt\mL^{q}_T(\wt\mL^{\bbp}_\x)}\left[\frac{C}{R}+\mP\big(|\xi^N_1|>\tfrac R2\big)\right]^{\frac12}\to 0, \ \ R\to\infty.\label{AD8}
\end{align}
On the other hand, for each $R>0$, by \eqref{CC9} again, we have
\begin{align*}
&\mE\left(\int^T_0|F_\eps-F|(s, X^{N,1}_s,r)\b1_{|X^{N,1}_s|\leq R}\dif s\right)\lesssim \nor (F_\eps-F)(\cdot,r)\b1_{B_R}\nor_{\wt\mL^{q}_T(\wt\mL^{\bbp}_\x)}\stackrel{\eqref{BZ30}}{\to 0},\ \ \eps\to 0,
 \end{align*}
which together with \eqref{AD8} yields \eqref{AD9}. 

Since $|F_\eps(t,x,r)-F_\eps(t,x,r')|\leq \kappa_0|r-r'|$, by \eqref{AD9} and a finite covering technique, for each $k\in\mN$, we further have
\begin{align}\label{AD10}
\lim_{\eps\to 0}\sup_N\mE\left(\int^T_0\sup_{|r|\leq k}|F_\eps-F|(s, X^{N,1}_s,r)\dif s\right)=0.
\end{align}
Indeed, for any given $\delta>0$, one can find $M$-balls in $\mR^m$ with centers in $\{r_i, i=1,\cdots,M\}$ and radius $\delta$ such that
$$
\big\{r: |r|\leq k\big\}\subset \cup_{i=1,\cdots,M}B_{\delta}(r_i).
$$
Thus,
\begin{align*}
\mE\left(\int^T_0\sup_{|r|\leq k}|F_\eps-F|(s, X^{N,1}_s,r)\dif s\right)
\leq\sum_{i=1}^M\mE\left(\int^T_0|F_\eps-F|(s, X^{N,1}_s,r_i)\dif s\right)+\kappa_0\delta.
\end{align*}
By \eqref{AD9} and firstly letting $\eps\to 0$ and then $\delta\to0$, we get \eqref{AD10}.

Moreover, for $j\not=1$, since
$$
\mE\left(\int^T_0|\phi^k_s-\phi_s|(X^{N,1}_s,X^{N,j}_s)\dif s\right)=\mE\left(\int^T_0|\phi^k_s-\phi_s|(X^{N,1}_s,X^{N,2}_s)\dif s\right),
$$
as in proving \eqref{AD9} and by \eqref{CC2} and \eqref{CC10}, we also have
$$
\lim_{k\to\infty}\sup_N\mE\left(\int^T_0|\phi^k_s-\phi_s|(X^{N,1}_s,X^{N,2}_s)\dif s\right)=0,
$$
and because $\phi_s(x,x)=0$ and \eqref{CC9},
$$
\lim_{k\to\infty}\sup_N\mE\left(\int^T_0|\phi^k_s|(X^{N,1}_s,X^{N,1}_s)\dif s\right)=0.
$$
Hence,
$$
\lim_{k\to\infty}\sup_N\sup_{j=1,\cdots,N}\mE\left(\int^T_0|\phi^k_s-\phi_s|(X^{N,1}_s,X^{N,j}_s)\dif s\right)=0,
$$
which together with \eqref{AD10} and \eqref{CC12} yields \eqref{CC11}.
\end{proof}

Now we are ready to give the
\begin{proof}[Proof of Theorem \ref{Th31}]
Consider the following random measure with values in $\cP(\mC_T)$,
$$
\omega\to\Pi_N(\omega,\dif w):=\frac1{N}\sum_{i=1}^N\delta_{X^{N,i}_\cdot(\omega)}(\dif w).
$$
By (i) of Lemma \ref{Le33} and \cite[(ii) of Proposition 2.2]{Sz}, the laws of $\Pi_N$, $ N\in\mN$, are tight in $\cP(\cP(\mC_T))$.
Without loss of generality, we assume that the laws of $\Pi_N$ weakly converge to some $\Pi_\infty\in \cP(\cP(\mC_T))$.
By \eqref{CC9} and \eqref{CC10}, it is standard to derive that for any $(q,\bbp)\in\sI_2$ and $f\in\wt\mL^{q}_T(\wt\mL^{\bbp}_\x)$ (see \cite[Remark 3.4]{Xi-Zh}),
\begin{align}\label{CC15}
\left|\int_{\cP(\mC_T)}\int_{\mC_T}\left(\int^T_0f(s, w_s)\dif s\right)\nu(\dif w)\Pi_\infty(\dif\nu)\right|\leq C\nor f\nor_{\wt\mL^{q}_T(\wt\mL^{\bbp}_\x)},
\end{align}
and for any $\bbp_1,\bbp_2\in(1,\infty)^d$ and $q\in(1,\infty)$ with $|\frac{1}{\bbp_1}|+|\frac{1}{\bbp_2}|+\frac{2}{q}<2$,
and $\x_1,\x_2\in\sX$, $f\in\wt\mL^q_T(\wt\mL^{\bbp_1}_{\x_1}(\wt\mL^{\bbp_2}_{\x_2}))$,
\begin{align}\label{CC16}
\left|\int_{\cP(\mC_T)}\int_{\mC_T}\int_{\mC_T}\left(\int^T_0f(s, w_s,w'_s)\dif s\right)\nu(\dif w')\nu(\dif w)\Pi_\infty(\dif\nu)\right|\leq 
C\nor f\nor_{\wt\mL^q_T(\wt\mL^{\bbp_1}_{\x_1}(\wt\mL^{\bbp_2}_{\x_2}))}.
\end{align}
Our aim below is to show that $\Pi_\infty$ is a Dirac measure, i.e.,
$$
\Pi_\infty(\dif\nu)=\delta_{\mu}(\dif\nu),\ \ \Pi_\infty-a.s.,
$$
where $\mu\in\cM^{\sigma,b}_{\mu_0}$ is the unique martingale solution of dDDSDE with initial distribution $\mu_0$.

We divide the proofs into two steps.

{\bf (Step 1)} For given $f\in C^2_0(\mR^d)$ and $\nu\in\cP(\mC_T)$, we define a functional on $\mC_T$ by
$$
M^{\sigma,b}_{f,\nu}(t,w):=f(w_t)-f(w_0)-\int^t_0\sL^{\sigma,b}_{\nu} f(s,w_s)\dif s,\ t\in[0,T],
$$
where 
$$
\sL^{\sigma,b}_{\nu}f(s,x):=\tfrac12\tr(\sigma\sigma^*\cdot\nabla^2 f)(s,x)+b(s,x,\nu_s)\cdot\nabla f(x),
$$
and 
$$
\nu_s:=\nu\circ w_s^{-1}\mbox{ is the marginal distribution of $\nu$ at time $s$}.
$$
Fix $n\in\mN$. For given  $g\in C_0(\mR^{nd})$ and $0\leq s_1<\cdots<s_n\leq s$, we also introduce a functional $\Xi^g_f$ on $\cP(\mC_T)$ by
$$
\Xi^g_f(\nu):=\int_{\mC_T} \(M^{\sigma,b}_{f,\nu}(t,w)-M^{\sigma,b}_{f,\nu}(s,w)\)g(w_{s_1},\cdots, w_{s_n})\nu(\dif w).
$$
In particular,
\begin{align}\label{BZ33}
\Xi^g_f(\Pi_N)=\frac{1}{N}\sum_{i=1}^N\(M^{\sigma,b}_{f,\Pi_N}(t,X^{N,i}_\cdot)-M^{\sigma,b}_{f,\Pi_N}(s,X^{N,i}_\cdot)\)
g\(X^{N,i}_{s_1},\cdots, X^{N,i}_{s_n}\)
\end{align}
and
$$
\Pi_N\circ w^{-1}_s=\eta_{\bX^N_s}.
$$
Noting that by It\^o's formula,
\begin{align*}
M^{\sigma,b}_{f,\Pi_N}(t,X^{N,i}_\cdot)&=f(X^{N,i}_t)-f(X^{N,i}_0)-\int^t_0\sL^{\sigma,b}_{\Pi_N} f(s,X^{N,i}_s)\dif s=\int^t_0(\sigma^*\cdot\nabla f)\(s, X_s^{N,i}\)\dif W^i_s,
\end{align*}
by \eqref{BZ33} and the It\^o isometry for stochastic integrals, we have
\begin{align}
\mE|\Xi^g_f(\Pi_N)|^2&=\frac{1}{N^2}\mE\left|\sum_{i=1}^N\int^t_s(\sigma^*\cdot\nabla f)\(r, X_r^{N,i}\)g\(X^{N,i}_{s_1},\cdots, X^{N,i}_{s_n}\)\dif W^i_r\right|^2\no\\
&=\frac{1}{N^2}\sum_{i=1}^N\int^t_s\mE\big|(\sigma^*\cdot\nabla f)\(r, X_r^{N,i}\)g\(X^{N,i}_{s_1},\cdots, X^{N,i}_{s_n}\)\big|^2\dif r\no\\
&\leq\frac{1}{N}(t-s)\|\sigma^*\cdot\nabla f\|^2_\infty\|g\|^2_\infty.\label{AA9}
\end{align}
Suppose that we have proven
\begin{align}\label{AA8}
\lim_{N\to\infty}\mE|\Xi^g_f(\Pi_N)|=\int_{\cP(\mC_T)}|\Xi^g_f(\nu)|\Pi_\infty(\dif \nu).
\end{align}
Then by \eqref{AA9} and \eqref{AA8}, for each $f\in C^2_0(\mR^d)$ and $n\in\mN$, $g\in C_0(\mR^{nd})$,
$$
\int_{\cP(\mC_T)}|\Xi^g_f(\nu)|\Pi_\infty(\dif \nu)=0\Rightarrow \Xi^g_f(\nu)=0\mbox{ for $\Pi_\infty$-a.s.  $\nu\in\cP(\mC_T)$}.
$$
Since $C^2_0(\mR^d)$ and $C_0(\mR^{nd})$ are separable, one can find a common $\Pi_\infty$-null set $\cN\subset\cP(\mC_T)$
such that for all $\nu\notin\cN$ and for all $0\leq s<t\leq T$, $f\in C_0^2(\mR^{d})$  and $n\in\mN$, $g\in C_0(\mR^{nd})$,
$$
\Xi^g_f(\nu)=\int_{\mC_T} \(M^{\sigma,b}_{f,\nu}(t,w)-M^{\sigma,b}_{f,\nu}(s,w)\)g(w_{s_1},\cdots, w_{s_n})\nu(\dif w)=0.
$$
Moreover, by \eqref{CC94} and \eqref{LAR}, we also have
$$
\Pi_\infty\{\nu\in\cP(\mC_T): \nu_0=\mu_0\}=1.
$$
Hence, for $\Pi_\infty$-almost all $\nu$,
$$
\nu\in\cM^{\sigma,b}_{\mu_0}.
$$
Since $\cM^{\sigma,b}_{\mu_0}$ only contains one point by uniqueness (see Theorem \ref{Th215}), all the points $\nu\notin\cN$ are the same.
Hence, $\Pi_N$ weakly converges to a one-point measure. By \cite[(ii) of Proposition 2.2]{Sz}, we conclude \eqref{CC14}. Thus it remains to show \eqref{AA8}. 

{\bf (Step 2)} Let $b_{\eps,k}$ be defined by \eqref{CC13} and define
$$
\Xi_{\eps,k}(\nu):=\int_{\mC_T} \(M^{\sigma,b_{\eps,k}}_{f,\nu}(t,w)-M^{\sigma,b_{\eps,k}}_{f,\nu}(s,w)\)g(w_{s_1},\cdots, w_{s_n})\nu(\dif w).
$$
By $b_{\eps,k}\in L^\infty_T(C_b(\mR^d\times\cP(\mR^d)))$, we have
\begin{align}\label{AA0}
\Xi_{\eps,k}\in C_b(\cP(\mC_T)),\ \ \forall\eps>0, k\in\mN.
\end{align}
Indeed, note that
\begin{align*}
\Xi_{\eps,k}(\nu)&=\int_{\mC_T} \left(f(w_t)-f(w_s)+\frac12\int^t_s\tr(\sigma\sigma^*\cdot\nabla^2 f)(r,w_r)\dif r\right)g(w_{s_1},\cdots, w_{s_n})\nu(\dif w)\\
&+\int_{\mC_T} \left(\int^t_s(b_{\eps,k}\cdot\nabla f)(r,w_r,\nu_r)\dif r\right)g(w_{s_1},\cdots, w_{s_n})\nu(\dif w)=:\Xi_{\eps,k}^{(1)}(\nu)+\Xi_{\eps,k}^{(2)}(\nu).
\end{align*}
Since $f\in C^2_b$ and $\sigma, g$ are bounded continuous, 
we have $\Xi_{\eps,k}^{(1)}\in C_b(\cP(\mC_T))$. 
For $\Xi_{\eps,k}^{(2)}$, since it is a non-linear functional of $\nu$, we have to take some care for the continuity of $\nu\mapsto \Xi_{\eps,k}^{(2)}(\nu)$.
Suppose that $\nu_m\in\cP(\mC_T)$ weakly converges to $\nu\in\cP(\mC_T)$. By definition, we have
\begin{align*}
|\Xi_{\eps,k}^{(2)}(\nu_m)-\Xi_{\eps,k}^{(2)}(\nu)|
&\leq
\left|\int_{\mC_T}\left(\int^t_s(b_{\eps,k}\cdot\nabla f)(r,w_r,\nu_r)\dif r\right)g(w_{s_1},\cdots, w_{s_n})(\nu_m-\nu)(\dif w)\right|\\
&+\kappa_0\|\nabla f\|_\infty\|g\|_\infty
\int_{\mC_T}\left(\int^t_s|\phi^k_r\otimes(\nu_m-\nu)_r|(w_r)\dif r\right)\nu_m(\dif w)
=:I^{(1)}_m+I^{(2)}_m,
\end{align*}
where we have used that
$$
|F_{\eps}(r,x,s_1)-F_{\eps}(r,x,s_2)|\leq\kappa_0|s_1-s_2|.
$$
For $I^{(1)}_m$, we clearly have
$$
\lim_{m\to\infty}I^{(1)}_m=0.
$$
For $I^{(2)}_m$, by the dominated convergence theorem, it suffices to show that for each $r\in[s,t]$,
$$
\lim_{m\to\infty}\int_{\mC_T}|\phi^k_r\otimes(\nu_m-\nu)_r|(w_r)\nu_m(\dif w)=0,
$$
which follows by noting that (see the proof of \eqref{AD10})
$$
\lim_{m\to\infty}|\phi^k_r\otimes(\nu_m-\nu)_r|(x)=0,\ x\in\mR^d,
$$
and
$$
\lim_{|x-y|\to 0}\sup_m|(\phi^k_r\otimes\nu_{m,r})(x)-(\phi^k_r\otimes\nu_{m,r})(y)|=0.
$$
Thus we get \eqref{AA0}, and so,
$$
\lim_{N\to\infty}\mE|\Xi_{\eps,k}(\Pi_N)|=\int_{\cP(\mC_T)}|\Xi_{\eps,k}(\nu)|\Pi_\infty(\dif \nu).
$$
On the other hand, we note that
$$
\Xi_{\eps,k}(\nu)-\Xi^g_f(\nu)=\int_{\mC_T}\left(\int^t_s(b-b_{\eps,k})(r, w_r,\nu_r)\cdot\nabla f(w_r)\dif r\right) g(w_{s_1},\cdots, w_{s_n})\nu(\dif w),
$$
and
$$
\Xi_{\eps,k}(\Pi_N)-\Xi^g_f(\Pi_N)=
\frac{1}{N}\sum_{i=1}^N\left(\int^t_s((b-b_{\eps,k})\cdot\nabla f)\(r, X^{N,i}_r,\eta_{\bX^N_r}\)\dif r\right) g\(X^{N,i}_{s_1},\cdots, X^{N,i}_{s_n}\).
$$
By \eqref{CC11}, we have
\begin{align*}
&\lim_{k\to\infty}\lim_{\eps\to 0}\sup_N\mE|\Xi_{\eps,k}(\Pi_N)-\Xi^g_f(\Pi_N)|\\
&\qquad\leq\|\nabla f\|_\infty\|g\|_\infty\lim_{k\to\infty}\lim_{\eps\to 0}\sup_N\mE\left(\int^t_s|b-b_{\eps,k}|\(r, X^{N,1}_r,\eta_{\bX^N_r}\)\dif r\right)=0,
\end{align*}
and by \eqref{CC15} and \eqref{CC16}, as in showing \eqref{CC11},
\begin{align*}
&\lim_{k\to\infty}\lim_{\eps\to 0}\int_{\cP(\mC_T)}|\Xi_{\eps,k}(\nu)-\Xi^g_f(\nu)|\Pi_\infty(\dif\nu)\\
&\quad\leq\|\nabla f\|_\infty\|g\|_\infty\lim_{k\to\infty}\lim_{\eps\to 0}\int_{\cP(\mC_T)}
\int_{\mC_T}\left(\int^T_0|b-b_{\eps,k}|(s, w_s,\nu_s)\dif s\right)\nu(\dif w)\Pi_\infty(\dif\nu)=0.
\end{align*}
Thus we obtain \eqref{AA8} and the proof is complete.
\end{proof}

\subsection{Entropy method}\label{55}
In this section we recall the entropy method used in \cite{JW18} to show a quantitative result for weak convergence when
the interaction kernel is bounded measurable, which is essentially contained in \cite{JW18}. For the completeness of the paper,
we provide a detailed proof. 
We first prepare the following lemma.
\bl\label{Le55}
Let $\phi:\mR^d\times\mR^d\to\mR$ be a bounded measurable function with $\phi(x,x)=0$
and $\boldsymbol{\xi}:=(\xi_1,\cdots, \xi_N)$ be a sequence of independent identical distributed random variables. Set
$$
\bar\phi(x,y):=\phi(x,y)-(\phi\circledast\mu)(x).
$$
Then for any $\lambda\leq\frac{1}{16\e^2\|\phi\|_\infty^2}$,
$$
\mE\e^{\lambda N|(\bar\phi\circledast\eta_{\boldsymbol{\xi}})(\xi_1)|^2}\leq 6,
$$
where $\eta_{\boldsymbol{\xi}}(\dif y):=\frac1N\sum_{i=1}^N\delta_{\xi_i}(\dif y)$.
\el
\begin{proof}
Note that by Taylor's expansion,
\begin{align*}
\e^{\lambda N|(\bar\phi\circledast\eta_{\boldsymbol{\xi}})(\xi_1)|^2}
&=\sum_{m=0}^\infty\frac{\lambda^m N^m}{m!}|(\bar\phi\circledast\eta_{\boldsymbol{\xi}})(\xi_1)|^{2m}
=\sum_{m=0}^\infty\frac{\lambda^m }{m! N^m}\Big|\sum_{j=1}^N\bar\phi(\xi_1,\xi_j)\Big|^{2m}\\
&\leq\sum_{m=0}^\infty\frac{\lambda^m }{m! N^m}2^{2m}\Big(|\bar\phi(\xi_1,\xi_1)|^{2m}+\Big|\sum_{j=2}^N\bar\phi(\xi_1,\xi_j)\Big|^{2m}\Big)\\
&\leq\sum_{m=0}^\infty\frac{(4\lambda)^m }{m! N^m}\left(\|\bar\phi\|^{2m}_\infty+\sum_{j_1,\cdots,j_{2m}=2}^N\bar\phi(\xi_1,\xi_{j_1})\cdots\bar\phi(\xi_1,\xi_{j_{2m}})\right).
\end{align*}
Let $\bJ$ be the set of all indices $(j_1,\cdots,j_{2m})\in\{2,\cdots,N\}^{2m}$ such that there is at least one index $j_k$ different from all others.
Since for $j\in\{2,\cdots,N\}$ and $x\in\mR^d$,
$$
\mE\bar\phi(x,\xi_j)=0,
$$
by the independence of the components of $\boldsymbol{\xi}$, we have for any $(j_1,\cdots,j_{2m})\in\bJ$,
$$
\mE\Big[\bar\phi(\xi_1,\xi_{j_1})\cdots\bar\phi(\xi_1,\xi_{j_{2m}})\Big]=\mE\Big[\mE\big[\bar\phi(x,\xi_{j_1})\cdots\bar\phi(x,\xi_{j_{2m}})\big]|_{x=\xi_1}\Big]=0.
$$
Hence,
$$
\mE\e^{\lambda N|(\bar\phi\circledast\eta_{\boldsymbol{\xi}})(\xi_1)|^2}\leq
\sum_{m=0}^\infty\frac{(4\lambda)^m }{m! N^m}\|\bar\phi\|_\infty^{2m}(1+\sharp \bJ^c),
$$
where $\sharp \bJ^c$ stands for the cardinality of the complement set $\bJ^c$. 

Suppose $2m\leq N$. It is easy to see that
$(j_1,\cdots,j_{2m})\in\bJ^c$ if and only if each $j_k$ appears at least twice and there are at most $m$-distinct $j_k$.
Thus one has
$$
\bJ^c=\cup_{n=1}^m\bJ_n,
$$
where $\bJ_n$ is the set of $(j_1,\cdots,j_{2m})$ such that each $j_k$ appears at least twice and exactly $n$-integers appear.
Clearly, by Stirling's formula $n^n\leq\e^nn!\leq \e^{2n} n^n$, we have
$$
\sharp\bJ_n\leq\left(
\begin{array}{c}
N-1\\
n
\end{array}
\right) n^{2m}=\frac{(N-1)^n}{n!} n^{2m}
\leq \frac{\e^n (N-1)^n}{n^n} n^{2m}\leq (N\e)^n n^m.
$$
Thus, for $2m\leq N$,
$$
\sharp\bJ^c\leq\sum_{n=1}^m(N\e)^n n^m\leq 2(N\e)^{m} m^m\leq2(N\e)^{m}\e^m m!.
$$
Moreover, for $2m>N$, we obviously have
$$
\sharp\bJ^c\leq N^{2m}\leq N^m (2m)^m\leq N^m(2\e)^{m} m!.
$$
So, for $\lambda\leq\frac{1}{16\e^2\|\phi\|_\infty^2}$,
$$
\mE\e^{\lambda N|(\bar\phi\circledast\eta_{\boldsymbol{\xi}})(\xi_1)|^2}\leq
\sum_{m=0}^\infty(4\lambda)^m \|\bar\phi\|_\infty^{2m}\Big(\frac{1}{m! N^m}+(2\e)^{m}\Big)
\leq 2\sum_{m=0}^\infty2^{-m}=6.
$$
The proof is complete.
\end{proof}
Now we can use the entropy formula in Lemma \ref{Lem34} to show the following result.
\bt\label{Th56}
Suppose that {\bf (H$^\sigma$)} and {\bf (H$^b$)} hold and $\phi$ is bounded measurable. Let $\mu^N_t$ be the law of $\bX^N_t$ in $\mR^{dN}$
and $\mu_t$ be the law of $X_t$ in $\mR^d$. Then there is a constant $C=C(\kappa_0,\kappa_1)>0$ independent of $\phi$ such that for any $t>0$,
\begin{align*}
\cH\big(\mu^{N}_t|\mu_t^{\otimes N}\big)\leq 
\e^{C\|\phi\|^2_\infty t}\Big(\cH\big(\mu^{N}_0|\mu_0^{\otimes N}\big)+C\|\phi\|^2_\infty t\Big).
\end{align*}
\et
\begin{proof}
Let $\eta_{\bw_s}:=\frac{1}{N}\sum_{i=1}^N\delta_{w^i_s}$ and $B, {\boldsymbol\sigma}$ be defined by \eqref{BB1} and \eqref{BB2}, respectively. 
By Lemma \ref{Lem34} and \eqref{CC1}, we have
\begin{align*}
\cH\big(\mu^{N}_t|\mu_t^{\otimes N}\big)
&\leq\cH\big(\mu^{N}_0|\mu_0^{\otimes N}\big)
+\frac{1}2\int^t_0\mE^{\mu^N_s}|{\boldsymbol\sigma}(s,\bw_s)^{-1}(B(s,\bw_s,\mu_s)-B(s,\bw_s,\eta_{\bw_s}))|^2\dif s\\
&\leq\cH\big(\mu^{N}_0|\mu_0^{\otimes N}\big)
+\frac{\kappa_0}{2}\int^t_0\mE^{\mu^N_s}|B(s,\bw_s,\mu_s)-B(s,\bw_s,\eta_{\bw_s})|^2\dif s\\
&\leq\cH\big(\mu^{N}_0|\mu_0^{\otimes N}\big)
+\frac{\kappa_0\kappa_1}{2}\sum_{i=1}^N\int^t_0\mE^{\mu^N_s}|(\phi_s\circledast\mu_s)(w^i_s)-(\phi_s\circledast\eta_{\bw_s})(w^i_s)|^2\dif s\\
&=\cH\big(\mu^{N}_0|\mu_0^{\otimes N}\big)
+\frac{\kappa_0\kappa_1}{2}\int^t_0N\mE^{\mu^N_s}|(\bar\phi_s\circledast\eta_{\bw_s})(w^1_s)|^2\dif s.
\end{align*}
Now by the variational representation \eqref{Var} and Lemma \ref{Le55} with $\lambda=\frac{1}{16\e^2\|\phi\|_\infty^2}$, we further have
\begin{align*}
\cH\big(\mu^{N}_t|\mu_t^{\otimes N}\big)
&\leq\cH\big(\mu^{N}_0|\mu_0^{\otimes N}\big)
+\frac{\kappa_0\kappa_1}{2\lambda}\int^t_0\Big[\cH\big(\mu^{N}_s|\mu_s^{\otimes N}\big)
+\ln\mE^{\mu^{\otimes N}_s}\e^{\lambda N|(\bar\phi_s\circledast\eta_{\bw_s})(w^1_s)|^2}\Big]\dif s\\
&\leq\cH\big(\mu^{N}_0|\mu_0^{\otimes N}\big)
+C\|\phi\|^2_\infty\int^t_0\Big[\cH\big(\mu^{N}_s|\mu_s^{\otimes N}\big)+\ln 6\Big]\dif s,
\end{align*}
which yields the desired estimate by Gronwall's inequality.
\end{proof}
\br\rm
By the Pinsker inequalities \eqref{Pin0} and \eqref{BB4}, we have for any $k\leq N$,
$$
\|\mu^{N,k}_t-\mu^{\otimes k}_t\|_{\rm var}\leq \sqrt{2\cH\(\mu^{N,k}_t|\mu^{\otimes k}_t\)}\leq 
\sqrt{\frac{\e^{C\|\phi\|^2_\infty t}k}{N}}\Big(\cH\big(\mu^{N}_0|\mu_0^{\otimes N}\big)+C\|\phi\|^2_\infty t\Big).
$$
Note that when $F(t,x,r)=r$ is linear and $\cH\(\mu^{N,k}_t|\mu^{\otimes k}_t\)\leq C_0k^2/N^2$, 
by a delicate analysis of the BBGKY hierarchy, the following sharp estimate is obtained by Lacker 
(see  Theorem 2.10 of  \cite{La21}):
$$
\|\mu^{N,k}_t-\mu^{\otimes k}_t\|_{\rm var}\leq \sqrt{2\cH\(\mu^{N,k}_t|\mu^{\otimes k}_t\)}\leq Ck/N.
$$
\er

\section{From weak convergence to strong convergence: Proof of Theorem \ref{Th11}}

In this section we show how to use the previous weak convergence result to derive the strong convergence of the particle system.
The following lemma is the key point.
\bl\label{Le61}
Let $\phi:\mR_+\times\mR^d\times\mR^d\to\mR$ be a measurable function. Set
$$
\bar \phi_t(x,y):=\phi_t(x,y)-(\phi_t\circledast\mu_{X_t})(x).
$$
(i) If $\phi$ is bounded measurable, then there is a constant $C=C(\kappa_0,\kappa_1)>0$ such that for all $t>0$,
\begin{align}\label{WX1}
\mE|(\bar \phi_t\circledast\eta_{\bX^N_t})(X^{N,1}_t)|^2\leq C\|\phi\|^2_\infty\e^{C\|\phi\|^2_\infty t}
\Big(\cH\big(\mu^{N}_0|\mu_0^{\otimes N}\big)+1\Big)/N.
\end{align}
(ii) If $\phi$ satisfies \eqref{CC2}, then for any $T>0$,
\begin{align}\label{WX2}
\lim_{N\to\infty}\mE\left(\int^T_0|(\bar \phi_t\circledast\eta_{\bX^N_t})(X^{N,1}_t)|^2\dif t\right)=0.
\end{align}
\el
\begin{proof}
(i) By the variational representation \eqref{Var}, for any $\eps>0$, we have
\begin{align*}
&\eps N\mE|(\bar \phi_t\circledast\eta_{\bX^N_t})(X^{N,1}_t)|^2=\eps N\mE^{\mu^N_{t}}|\bar \phi_t(w^1_t,\eta_{\bw_t})|^2
\leq \cH(\mu^N_{t}|\mu^{\otimes N}_t)
+\log\mE^{\mu^{\otimes N}_t}\e^{\eps N|\bar \phi_t(w^1_t,\eta_{\bw_t})|^2},
\end{align*}
which in turn implies \eqref{WX1} by Lemma \ref{Le55} with $\eps=\frac{1}{16\e^2\|\phi\|_\infty^2}$ and Theorem \ref{Th56}.

(ii)  By definition we have
\begin{align}\label{AA2}
\mE\left(\int^T_0|(\bar \phi_t\circledast\eta_{\bX^N_t})(X^{N,1}_t)|^2\dif t\right)
=\frac1{N^2}\sum_{j,k=1}^N\mE\left(\int^T_0\Gamma_t\(X^{N,1}_t,X^{N,j}_t,X^{N,k}_t\)\dif t\right),
\end{align}
where
$$
\Gamma_t(x,y,z):=\bar \phi_t(x,y)\bar \phi_t(x,z).
$$
Let $\phi^\eps_t(x,y):=(\phi_t*\varGamma_\eps)(x,y)$ be the mollifying approximation of $\phi_t$ and
$$
\bar \phi^\eps_t(x,y):=\phi^\eps_t(x,y)-(\phi^\eps_t\circledast\mu_{X_t})(x),
$$
and
$$
\Gamma_t^\eps(x,y,z):=\bar \phi^\eps_t(x,y)\bar \phi^\eps_t(x,z).
$$
Noting that
$$
(\Gamma_t-\Gamma_t^\eps)(x,y,z)=(\bar \phi_t-\bar \phi^\eps_t)(x,y)\bar \phi^\eps_t(x,z)+
\bar \phi_t(x,y)(\bar \phi_t-\bar \phi^\eps_t)(x,z),
$$
by H\"older's inequality, we have
\begin{align*}
I^{N}_{j,k}(\eps)&:=\left|\mE\left(\int^T_0(\Gamma_t-\Gamma_t^\eps)\(X^{N,1}_t,X^{N,j}_t,X^{N,k}_t\)\dif t\right)\right|\\
&\leq \left(\mE\int^T_0(\bar \phi_t-\bar \phi^\eps_t)^2\(X^{N,1}_t,X^{N,j}_t\)\dif t\right)^{1/2}
\left(\mE\int^T_0\bar \phi^\eps_t\(X^{N,1}_t,X^{N,k}_t\)^2\dif t\right)^{1/2}\\
&+\left(\mE\int^T_0\bar \phi_t\(X^{N,1}_t,X^{N,j}_t\)^2\dif t\right)^{1/2}
\left(\mE\int^T_0(\bar \phi_t-\bar \phi^\eps_t)^2\(X^{N,1}_t,X^{N,k}_t\)\dif t\right)^{1/2}.
\end{align*}
Using the Krylov estimate \eqref{CC10} and as in showing \eqref{AD9}, we get
\begin{align}
\lim_{\eps\to 0}\sup_N\sup_{j,k}I^{N}_{j,k}(\eps)=0.\label{AA4}
\end{align}
On the other hand, for fixed $\eps$, by \eqref{CC14} we have
\begin{align}
&\lim_{N\to\infty}\sup_{j\not=k\not=1}\mE\left(\int^T_0\Gamma_t^\eps\(X^{N,1}_t,X^{N,j}_t,X^{N,k}_t\)\dif t\right)\no\\
&\quad=\lim_{N\to\infty}\mE\left(\int^T_0\Gamma_t^\eps\(X^{N,1}_t,X^{N,2}_t,X^{N,3}_t\)\dif t\right)\no\\
&\quad=\mE\left(\int^T_0\Gamma_t^\eps\(X^{1}_t,X^{2}_t,X^{3}_t\)\dif t\right)=0,\label{AA404}
\end{align}
where the last step is due to the fact that
$$
\mE\Gamma^\eps_t\(X^{1}_t,X^2_t,X^3_t\)=\mE\left[\mE\bar\phi^\eps_t(x,X^2_t)\mE\bar\phi^\eps_t(x,X^3_t); x=X^{1}_t \right]=0.
$$
Thus by \eqref{AA4}  and \eqref{AA404},
\begin{align}\label{AA22}
\lim_{N\to\infty}\sup_{j\not=k\not=1}\mE\left(\int^T_0\Gamma_t(X^{N,1}_t,X^{N,j}_t,X^{N,k}_t\)\dif t\right)=0.
\end{align}
Moreover, by the Krylov estimate \eqref{CC10} we also have
\begin{align}\label{AA202}
&\sup_{j,k}\mE\left(\int^T_0\Gamma_t\(X^{N,1}_t,X^{N,j}_t,X^{N,k}_t\)\dif t\right)\no\\
&\quad\leq \sup_{j,k}\mE\left(\int^T_0\bar\phi_t\(X^{N,1}_t,X^{N,j}_t\)^2\dif t\right)^{\frac12}
\mE\left(\int^T_0\bar\phi_t\(X^{N,1}_t,X^{N,k}_t\)^2\dif t\right)^{\frac12}\no\\
&\quad=\sup_j \mE\left(\int^T_0\bar\phi_t\(X^{N,1}_t,X^{N,j}_t\)^2\dif t\right)<\infty.
\end{align}
By \eqref{AA2}, \eqref{AA22} and \eqref{AA202}, we obtain \eqref{WX2}.
\end{proof}

Now we can give the
\begin{proof}[Proof of Theorem \ref{Th11}]
Let $X_t$ be the unique strong solution of dDDSDE \eqref{MV1} starting from $X_0$ (see Theorem \ref{Th215}). Define
$$
\bar b(t,x):=b(t,x,\mu_{X_t})=F(t,x,(\phi_t\circledast\mu_{X_t})(x)).
$$
By {\bf (H$^b$)}, it is easy to see that
$$
\flat:=\nor \bar b\nor_{\mL^q_T(\wt\mL^\bbp_\x)}<\infty.
$$
Consider the following backward PDE
$$
\p_t \u+\tfrac12\tr(\sigma\sigma^*\cdot\nabla^2 \u)+\bar b\cdot\nabla \u-\lambda \u+\bar b=0,\ \ \u(T)=0.
$$
By reversing the time variable and Theorem \ref{Pre3:Sch0}, there is a unique solution $\u$ satisfying the following estimate: for any $\beta\in(0,\vartheta)$,
where $\vartheta:=1-|\frac1\bbp|-\frac2q$,
there is a constant $C_0=C_0(T,\kappa_0, d,\bbp,q,\beta)\geq 1$ such that for all $\lambda\geq C_0\flat^{2/\vartheta}$,
\begin{align}\label{CC09}
\lambda^{\frac12(\vartheta-\beta)}\|\u\|_{\mL^\infty_T(\cC^{1+\beta})}+
\nor\nabla^2 \u\nor_{\widetilde{\mL}^{q}_T(\wt\mL^{\bbp}_{\x})}\le C_0\flat.
\end{align}
In particular, one can choose $\lambda=(2C_0\flat)^{2/\vartheta}$ so that
\begin{align}\label{CC091}
\|\nabla \u\|_{\mL^\infty_T}\leq\tfrac12.
\end{align}
Now if we define
$$
\Phi(t,x):=x+\u(t,x),
$$
then for each $t$,
$$
x\mapsto\Phi(t,x)\mbox{ is a $C^1$-diffeomorphism on $\mR^d$},
$$
and
\begin{align}\label{CC93}
\|\nabla\Phi\|_{\mL^\infty_T}+\|\nabla\Phi^{-1}\|_{\mL^\infty_T}\leq 2.
\end{align}
Define
$$
Y_t:=\Phi(t,X_t),\ \ Y^{N,1}_t:=\Phi(t,X^{N,1}_t).
$$
By It\^o's formula (see the proof in Lemma \ref{Le28}), we have
\begin{align*}
\dif Y_t=\lambda \u(t,X_t)\dif t+\tilde\sigma(t,X_t)\dif W^1_t
\end{align*}
and
\begin{align*}
\dif Y^{N,1}_t=\lambda \u(t,X^{N,1}_t)\dif t+\(B\cdot\nabla \u\)\(t,X^{N,1}_t\)\dif t+\tilde\sigma(t,X^{N,1}_t)\dif W^1_t,
\end{align*}
where $\tilde\sigma:=\sigma^*\nabla\Phi$ and
$$
B(t,x):=b(t,x,\eta_{\bX^N_t})-b(t,x,\mu_{X_t}).
$$
In particular, we have
\begin{align*}
Y^{N,1}_t-Y_t&=\Phi(0,X^{N,1}_0)-\Phi(0,X_0)+\lambda\int^t_0\Big[\u(s,X^{N,1}_s)-\u(s,X_s)\Big]\dif s\\
&\quad+\int^t_0\(B\cdot\nabla \u\)\(s,X^{N,1}_s\)\dif s+\int^t_0\Big[\tilde\sigma(s,X^{N,1}_s)-\tilde\sigma(s,X_s)\Big]\dif W^1_s.
\end{align*}
By It\^o's formula and \eqref{CC091},  \eqref{CC93}, we further have
\begin{align}\label{AX8}
\begin{split}
|Y^{N,1}_t-Y_t|^2&\leq 4|X^{N,1}_0-X_0|^2+\int^t_0|Y^{N,1}_s-Y_s|\(\lambda|X^{N,1}_s-X_s|+|B\(s,X^{N,1}_s\)|\)\dif s\\
&\qquad +\int^t_0|\tilde\sigma(s,X^{N,1}_s)-\tilde\sigma(s,X_s)|^2\dif s+M_t,
\end{split}
\end{align}
where $M_t$ is a continuous local martingale.
Note that by \eqref{Pre1:Mnf},
$$
|\tilde\sigma(s,X^{N,1}_s)-\tilde\sigma(s,X_s)|^2\leq 2\ell_{N,0}(s)|X^{N,1}_s-X_s|^2,
$$
where
$$
\ell_{N,\lambda}(s):=\cM|\nabla\tilde\sigma(s,\cdot)|^2(X^{N,1}_s)+\cM|\nabla\tilde\sigma(s,\cdot)|^2(X_s)+\|\tilde\sigma\|^2_\infty+\lambda+1.
$$
Thus, by \eqref{AX8} and \eqref{CC93} we have
\begin{align}\label{AM5}
\begin{split}
|X^{N,1}_t-X_t|^2&\leq C \Bigg(|X^{N,1}_0-X_0|^2+\int^t_0\ell_{N,\lambda}(s)|X^{N,1}_s-X_s|^2\dif s\\
&\qquad\quad+\int^t_0|B\(s,X^{N,1}_s\)|^2\dif s\Bigg)+M_t,
\end{split}
\end{align}
where $C>0$ is an absolute constant.
By the chain rule, we have
$$
\cM|\nabla\tilde\sigma|^2\leq 4\cM|\nabla\sigma|^2+\|\sigma\|_\infty^2\cM|\nabla^2\u|^2.
$$
By \eqref{Pre1:Mix} and \eqref{SIG0}, we have
\begin{align*}
\nor\cM|\nabla\sigma|^2\nor_{\mL^{q_0/2}_T(\wt\mL^{\bbp_0/2}_\x)}\lesssim
\nor|\nabla\sigma|^2\nor_{\mL^{q_0/2}_T(\wt\mL^{\bbp_0/2}_\x)}
=\nor \nabla\sigma\nor_{\mL^{q_0}_T(\wt\mL^{\bbp_0}_\x)}^2\leq\kappa_0,
\end{align*}
and by \eqref{CC09},
\begin{align*}
\nor\cM|\nabla^2\u|^2\nor_{\mL^{q/2}_T(\wt\mL^{\bbp/2}_\x)}\lesssim
\nor|\nabla^2\u|^2\nor_{\mL^{q/2}_T(\wt\mL^{\bbp/2}_\x)}
=\nor \nabla^2\u\nor_{\mL^{q}_T(\wt\mL^{\bbp}_\x)}^2\leq (C_0\flat)^2.
\end{align*}
Since $(\frac{q_0}2,\frac{\bbp_0}2), (\frac{q}2,\frac{\bbp}2)\in\sI_2$, by \eqref{CC98} and \eqref{Kha1} we have for any $\gamma>0$,
$$
A_\gamma:=\sup_N\mE\exp\left\{\gamma\int^T_0\ell_{N,\lambda}(s)\dif s\right\}<\infty.
$$
Thus by \eqref{AM5} and the stochastic Gronwall inequality (cf. \cite{Sc} or \cite[Lemma 3.7]{Xi-Zh}), we get for any $\gamma\in(0,1)$,
\begin{align}\label{AM1}
\mE\left(\sup_{t\in[0,T]}|X^{N,1}_t-X_t|^{2\gamma}\right)\leq C_\gamma A_{\frac{\gamma+1}{\gamma-1}}\left(\mE|X^{N,1}_0-X_0|^2
+\mE\int^T_0|B\(s,X^{N,1}_s\)|^2\dif s\right)^\gamma.
\end{align}
Noting that by \eqref{CC1},
$$
|B(t,x)|\leq\kappa_1|(\phi_t\circledast\eta_{\bX^N_t})(x)-(\phi_t\circledast \mu_{X_t})(x)|=\kappa_1|(\bar \phi_t\circledast\eta_{\bX^N_t})(x)|,
$$
where 
$$
\bar\phi_t(x,y):=\phi_t(x,y)-(\phi_t\circledast \mu_{X_t})(x),
$$
we further have for any $\gamma\in(0,1)$,
$$
\mE\left(\sup_{t\in[0,T]}|X^{N,1}_t-X_t|^{2\gamma}\right)\leq C_\gamma A_{\frac{\gamma+1}{\gamma-1}}\left(\mE|X^{N,1}_0-X_0|^2
+\kappa^2_1\mE\int^T_0|(\bar \phi_s\circledast\eta_{\bX^N_s})(X^{N,1}_s)|^2\dif s\right)^{\gamma}.
$$

Now, (i) follows by \eqref{WX2} and the above estimate.

(ii)  When $h$ and $\phi$ are bounded, by \eqref{CC1} one has
$$
|\bar b(t,x)|\leq \|h\|_\infty+\kappa_1\|\phi\|_\infty.
$$
Thus for any $\delta>2$, one can choose $q,\bbp$ in \eqref{CC09} close to $\infty$ so that $\vartheta=\frac{2}{\delta}=1-\frac2q-|\frac1\bbp|$ and
$$
\flat:=\nor \bar b\nor_{\mL^q_T(\wt\mL^\bbp_\x)}\leq C(1+\|\phi\|_\infty).
$$
By \eqref{CC98},  \eqref{Kha11} and for $\lambda=(2C_0\flat)^{2/\vartheta}$,  we have
$$
A_\gamma=\sup_N\mE\exp\left\{\gamma\int^T_0\ell_{N,\lambda}(s)\dif s\right\}\leq C\e^{C\flat^{2/\vartheta}}\leq C\e^{C\|\phi\|_\infty^{2/\vartheta}}.
$$
Estimate \eqref{LIM2} now follows by the above estimates and \eqref{WX1}.
\end{proof}

\section{Moderate interacting particle system: Proof of Theorem \ref{Th12}}

We consider the following McKean-Vlasov type approximation for density-dependent SDE \eqref{MV3}:
$$
\dif X^\eps_t=F(t,X^\eps_t,(\phi_{\eps}*\rho^\eps_t)(X^\eps_t))\dif t+\sigma(t,X^\eps_t)\dif W^1_t,\ X^\eps_0=X_0,
$$
where $\phi_\eps(x)=\eps^{-d}\phi(x/\eps)$, and $\phi$ is a bounded probability density function with support in the unit ball, $F$ is bounded measurable and
$\rho^\eps_t$ is the density of $X^\eps_t$.

We first show the following lemma.
\bl\label{Le71}
For any $T>0$, $\beta\in(0,\gamma_0)$ and $\gamma\in(0,1)$, there is a constant $C=C(T,\beta,\gamma,\Theta)>0$ such that for all $\eps\in(0,1)$,
$$
\mE\left(\sup_{t\in[0,T]}|X^\eps_t-X_t|^{2\gamma}\right)\leq C \eps^{2\beta \gamma}.
$$
\el
\begin{proof}
Let $X_t$ be the unique strong solution of DDSDE \eqref{MV3} starting from $X_0$. Define
$$
\bar b(t,x):=F(t,x,\rho_t(x)).
$$
By assumption we have
$$
\|\bar b\|_{\mL^\infty_T}\leq\|F\|_{\mL^\infty_T}.
$$
Consider the following backward PDE
$$
\p_t \u+\tfrac12\tr(\sigma\sigma^*\cdot\nabla^2 \u)+\bar b\cdot\nabla \u-\lambda \u+\bar b=0,\ \ \u(T)=0.
$$
As in the proof of Theorem \ref{Th11} we construct a $C^1$-diffeomorphism
$$
\Phi(t,x):=x+\u(t,x),
$$
and define
$$
Y^\eps_t:=\Phi(t,X^\eps_t),\ \ Y_t:=\Phi(t,X_t).
$$
By the generalized It\^o formula, we have
\begin{align*}
\dif  Y_t=\lambda \u(t, X_t)\dif t+\tilde\sigma(t, X_t)\dif W^1_t
\end{align*}
 and
\begin{align*}
\dif Y^\eps_t=\lambda \u(t,X^\eps_t)\dif t+\(B_\eps\cdot\nabla \u\)\(t,X^\eps_t\)\dif t+\tilde\sigma(t,X^\eps_t)\dif W^1_t,
\end{align*}
where $\tilde\sigma=\sigma^*\nabla\Phi$ and
$$
B_\eps(t,x):=F(t,x,(\phi_\eps*\rho^\eps_t)(x))-F(t,x,\rho_t(x)).
$$
In particular, we have
\begin{align*}
Y^\eps_t- Y_t&=\lambda\int^t_0\Big[\u(s,X^\eps_s)-\u(s, X_s)\Big]\dif s+\int^t_0\(B_\eps\cdot\nabla \u\)\(s,X^\eps_s\)\dif s\\
&\quad+\int^t_0\Big[\tilde\sigma(s,X^\eps_s)-\tilde\sigma(s, X_s)\Big]\dif W^1_s.
\end{align*}
By It\^o's formula and \eqref{CC091}, we further have
\begin{align}\label{AX88}
\begin{split}
|Y^\eps_t- Y_t|^2&\leq \int^t_0|Y^\eps_s- Y_s| \(\lambda |X^\eps_s- X_s|+|B_\eps\(s,X^\eps_s\)|\)\dif s\\
&\quad+\int^t_0|\tilde\sigma(s,X^\eps_s)-\tilde\sigma(s, X_s)|^2\dif s+M_t,
\end{split}
\end{align}
where $M_t$ is a continuous local martingale. Completely the same way as in proving \eqref{AM1}, we have
\begin{align}\label{GA1}
\mE|X^{N,1}_t- X_t|^{2\gamma}\lesssim\left(\mE\int^T_0|B_\eps\(s,X^\eps_s\)|^2\dif s\right)^{\gamma}.
\end{align}
On the other hand, for any $p>d$, by Lemma \ref{Le47} we have
\begin{align*}
\|\rho^\eps_t-\rho_t\|_{\mL^\infty}
\lesssim_C \int^t_0 (t-s)^{-\frac12(1+\frac dp)}\nor B_\eps(s)\nor_{\wt\mL^p}\dif s.
\end{align*}
By the Lipschitz assumption on $F$ in $r$, we have
\begin{align*}
\nor B_\eps(s)\nor_{\wt\mL^p}\leq\|B_\eps(s)\nor_{\mL^\infty}
&\lesssim\| \phi_{\eps}*\rho^\eps_s-\rho_s\|_{\mL^\infty}
\leq\|\rho^\eps_s-\rho_s\|_{\mL^\infty}+\| \phi_{\eps}*\rho_s-\rho_s\|_{\mL^\infty}.
\end{align*}
For any $\beta\in(0,\gamma_0)$, noting that by \eqref{Ho},
$$
\|\rho_s(\cdot+y)-\rho_s\|_{\mL^\infty}\leq C\|\rho_0\|_\infty|y|^\beta s^{-\beta/2},
$$
we have
\begin{align}
\|\phi_{\eps}*\rho_s-\rho_s\|_{\mL^\infty}
&\leq\int_{\mR^d}\|\rho_s(\cdot+y)-\rho_s\|_{\mL^\infty}\cdot|\phi_\eps(y)|\dif y\no\\
&\lesssim s^{-\beta/2}\int_{\mR^d}|y|^\beta\cdot|\phi_\eps(y)|\dif y\lesssim s^{-\beta/2}\eps^\beta.\label{ZA1}
\end{align}
Hence,
\begin{align*}
\|\rho^\eps_t-\rho_t\|_{\mL^\infty}
\lesssim_C 
\int^t_0 (t-s)^{-\frac12(1+\frac dp)}(\|\rho^\eps_s-\rho_s\|_{\mL^\infty}
+s^{-\frac\beta2}\eps^\beta)\dif s.
\end{align*}
By Gronwall's inequality of Volterra's type, we have
\begin{align}\label{ZA2}
\|\rho^\eps_t-\rho_t\|_{\mL^\infty}\leq C t^{\frac12-\frac{d}{2p}-\frac{\beta}{2}} \eps^\beta\le C t^{-\frac{\beta}{2}} \eps^\beta .
\end{align}
Note that by \eqref{CC1}, \eqref{ZA1} and \eqref{ZA2},
\begin{align*}
\mE|B_\eps\(s,X^\eps_s\)|^2&\leq\kappa_1^2\int_{\mR^d}|\phi_\eps*\rho^\eps_s(x)-\rho_s(x)|^2\rho_s^\eps(x)\dif x\\
&\leq\kappa_1^2\|\phi_\eps*\rho^\eps_s-\rho_s\|^2_{\mL^\infty}\leq Cs^{-\beta}\eps^{2\beta}.
\end{align*}
Substituting this into \eqref{GA1}, we obtain the desired estimate.
\end{proof}

Now we can give the
\begin{proof}[Proof of Theorem \ref{Th12}]This is a direct combination of Lemma \ref{Le71} and (ii) of Theorem \ref{Th11}.
\end{proof}

\end{document}